\documentclass[]{scrartcl}

\usepackage[utf8]{inputenc}

\usepackage[USenglish]{babel}
\usepackage{csquotes}

\usepackage[a4paper,top=27mm,bottom=20mm,inner=25mm,outer=20mm]{geometry}

\usepackage[%
  backend=bibtex,bibencoding=ascii,
  style=numeric-comp,
  giveninits=true, uniquename=init, 
  natbib=true,
  url=true,
  doi=true,
  isbn=false,
  backref=false,
  maxnames=99,
  ]{biblatex}
\usepackage{amsmath}
\allowdisplaybreaks
\numberwithin{equation}{section}
\usepackage{amssymb}
\usepackage{commath}
\usepackage{mathtools}
\usepackage{bbm}
\usepackage{nicefrac}
\usepackage{subdepth}

\usepackage{siunitx}
\usepackage{algorithm}
\usepackage{algpseudocode}
\usepackage{todonotes}

\usepackage{amsthm}
\usepackage{thmtools}
\usepackage{etoolbox}
\makeatletter
\patchcmd{\thmt@setheadstyle}
 {\bgroup\thmt@space}
 {\thmt@space}
 {}{}
\patchcmd{\thmt@setheadstyle}
 {\egroup\fi}
 {\fi}
 {}{}
\makeatother
\declaretheoremstyle[
  bodyfont=\normalfont\itshape,
  headformat=\NAME\ \NUMBER\NOTE,
]{myplain}
\declaretheoremstyle[
  headformat=\NAME\ \NUMBER\NOTE,
]{mydefinition}
\newcommand{\envqed}{{\lower-0.3ex\hbox{$\triangleleft$}}}
\declaretheorem[style=myplain,numberwithin=section]{theorem}

\declaretheorem[style=mydefinition,numberlike=theorem,qed=\envqed]{definition}
\declaretheorem[style=mydefinition,numberlike=theorem,qed=\envqed]{remark}

\usepackage[plainpages=false,pdfpagelabels,hidelinks,unicode]{hyperref}

\usepackage{color}
\usepackage{graphicx}
\usepackage[small]{caption}
\usepackage{subfigure,placeins,paralist}

\begingroup\expandafter\expandafter\expandafter\endgroup
\expandafter\ifx\csname pdfsuppresswarningpagegroup\endcsname\relax
\else
\fi

\usepackage{booktabs}
\usepackage{rotating}
\usepackage{multirow}

\usepackage{ifluatex}
\ifluatex
  \usepackage[no-math]{fontspec}
\else
  \usepackage[T1]{fontenc}
\fi
\usepackage{newpxtext,newpxmath}
\makeatletter
\newcommand*\dashline{\rotatebox[origin=c]{90}{$\dabar@\dabar@\dabar@$}}
\makeatother
\usepackage{trimclip}

\usepackage{xparse}

\let\epsilon\varepsilon
\let\phi\varphi
\let\rho\varrho

\newsavebox{\DelimiterBox}
\newlength{\DelimiterHeight}
\newlength{\DelimiterDepth}
\newsavebox{\ArgumentBox}
\newlength{\ArgumentHeight}
\newlength{\ArgumentDepth}
\newlength{\ResizedDelimiterHeight}
\newlength{\ResizedDelimiterDepth}

\newcommand{\encloseby}[3]{%
  \savebox{\ArgumentBox}{$\displaystyle #1$}%
  \settoheight{\ArgumentHeight}{\usebox{\ArgumentBox}}%
  \settodepth{\ArgumentDepth}{\usebox{\ArgumentBox}}%
  \savebox{\DelimiterBox}{#2}%
  \settoheight{\DelimiterHeight}{\usebox{\DelimiterBox}}%
  \settodepth{\DelimiterDepth}{\usebox{\DelimiterBox}}%
  \setlength{\ResizedDelimiterHeight}{%
    \maxof{1.2\ArgumentHeight}{\DelimiterHeight}%
  }
  \setlength{\ResizedDelimiterDepth}{%
    \maxof{1.2\ArgumentDepth}{\DelimiterDepth}%
  }
  \raisebox{-\ResizedDelimiterDepth}{%
    \resizebox{\width}{\ResizedDelimiterHeight+\ResizedDelimiterDepth}{%
      \raisebox{\DelimiterDepth}{#2}%
    }%
  }
  #1
  \raisebox{-\ResizedDelimiterDepth}{%
    \resizebox{\width}{\ResizedDelimiterHeight+\ResizedDelimiterDepth}{%
      \raisebox{\DelimiterDepth}{#3}%
    }%
  }
}

\ifluatex
  \newcommand{\mean}[1]{\encloseby{#1}{$\{\mkern-5mu\{$}{$\}\mkern-5mu\}$}}
  
  \newcommand{\jump}[1]{\encloseby{#1}{$[\mkern-4mu[$}{$]\mkern-4mu]$}}
\else
  \newcommand{\mean}[1]{\encloseby{#1}{$\{\mkern-6mu\{$}{$\}\mkern-6mu\}$}}
  
  \newcommand{\jump}[1]{\encloseby{#1}{$[\mkern-3mu[$}{$]\mkern-3mu]$}}
\fi

\newcommand{\R}{\mathbb{R}}
\newcommand{\N}{\mathbb{N}}
\newcommand{\dd}{\textup{d}}
\newcommand{\bb}[1]{\boldsymbol{#1}}
\newcommand{\dbl}{\left\llbracket}
\newcommand{\dbr}{\right\rrbracket}

\newcommand{\orcid}[1]{ORCID:~\href{https://orcid.org/#1}{#1}}
\usepackage{authblk}

\newenvironment{keywords}{\par\textbf{Key words.}}{\par}
\newenvironment{AMS}{\par\textbf{AMS subject classification.}}{\par}

\title{Computing Radially-Symmetric Solutions of the Ultra-Relativistic Euler Equations with Entropy-Stable Discontinuous Galerkin Methods}

\author[1]{Ferdinand~Thein\thanks{\orcid{0000-0002-0170-8284}}}
\affil[1]{Institute of Mathematics, Johannes Gutenberg University Mainz, Staudingerweg 9, 55128 Mainz, Germany}

\author[1]{Hendrik~Ranocha\thanks{\orcid{0000-0002-3456-2277}}}

\date{February 27, 2026} 

\makeatletter
\hypersetup{pdfauthor={Ferdinand~Thein, Hendrik~Ranocha}} 
\hypersetup{pdftitle={Computing Radially-Symmetric Solutions of the Ultra-Relativistic Euler Equations with Entropy-Stable Discontinuous Galerkin Methods}} 
\makeatother

\begin{document}

\maketitle

\begin{abstract}
\noindent
The ultra--relativistic Euler equations describe gases in the relativistic case when the thermal energy dominates.
These equations for an ideal gas are given in terms of the pressure, the spatial part of the dimensionless four-velocity, and the particle density.
Kunik et al.\ (2024, \url{https://doi.org/10.1016/j.jcp.2024.113330}) proposed genuine multi--dimensional benchmark problems for the ultra--relativistic Euler equations.
In particular, they compared full two-dimensional discontinuous Galerkin simulations for radially symmetric problems with solutions computed using a specific one-dimensional scheme.
Of particular interest in the solutions are the formation of shock waves and a pressure blow-up.
In the present work we derive an entropy-stable flux for the ultra--relativistic Euler equations.
Therefore, we derive the main field (or entropy variables) and the corresponding potentials.
We then present the entropy-stable flux and conclude with simulation results for different test cases both in 2D and in 3D.
\end{abstract}



\section{Introduction}\label{sec:intro}
Relativistic fluid dynamics is an area of ongoing study in recent research literature, covering various scales from ions to neutron stars \cite{Andersson2021}. The commonality is that the energy or velocity of the fluid in question is so significant that it surpasses the confines of classical Newtonian mechanics. For a recent review of relativistic fluid dynamics, we refer to \cite{Disconzi2024}. In this work, we examine relativistic fluid dynamics in the context of hyperbolic partial differential equations (PDEs). The study of the structural properties of the governing equations, based on their symmetric structure, can be found in \cite{Ruggeri2021a,Ruggeri2021b,Ruggeri2021,Pennisi2019}. For works exploring diverse aspects of the relativistic Euler equations, we refer to \cite{Steinhardt1982,Smoller1993,Chen1997}. Additional recent results can be found in \cite{Groah2007,Smoller2012a,Lai2019,Freistuehler2019,Freistuehler2024}. Here, we are particularly interested in the ultra-relativistic Euler equations, i.e., the scenario where thermal energy predominates. Results on this specific system can be located in the aforementioned literature. For a general introduction to the mathematical theory of hyperbolic conservation laws, see Smoller \cite{Smoller2012} and Dafermos \cite{Dafermos2016}.

In this work, we focus on the numerical aspects of the relativistic flows under consideration. For a recent treatment of numerical methods for relativistic fluid dynamics, we refer to \cite{Zanotti2015,Gaburro2021,Dumbser2024} and the references therein. Previous results on the numerical treatment of the ultra-relativistic Euler equations are presented by Abdelrahman et al.\ \cite{Abdelrahman2014}, who propose a front tracking scheme; for kinetic schemes, we cite Kunik et al.\ \cite{Kunik2003,Kunik2003a,Kunik2004}. In Kunik et al.\ \cite{Kunik2021}, the ultra-relativistic Euler equations are studied, and a one-dimensional radially symmetric solver that can provide non-trivial reference solutions for genuinely three-dimensional calculations is developed. Recently, in \cite{kunik2024radially}, this work continues, and multi--dimensional benchmark problems are introduced. These benchmarks are solved using full-dimensional simulations, and the results are compared to the specific solver provided in \cite{Kunik2021,kunik2024radially}.
In this work, we extend the previous works by developing an entropy-stable method for the ultra-relativistic equations and applying it not only to two-dimensional but also to more challenging three-dimensional problems.
The entropy-stable numerical methods developed in this paper are based on the framework pioneered by Tadmor \cite{tadmor1987numerical,tadmor2003entropy}. A key ingredient of the method is an entropy-conservative two-point flux, which we derive for the first time for the ultra-relativistic Euler equations. This flux is not based on integral averages like the first entropy-conservative numerical fluxes proposed by Tadmor, but rather on differential averages similar to the affordable numerical fluxes developed for the compressible Euler equations of an ideal gas in the non-relativistic setting \cite{ismail2009affordable,chandrashekar2013kinetic,ranocha2018comparison}. The second-order method is extended to a higher-order scheme using flux differencing \cite{lefloch2002fully,fisher2013high,ranocha2018comparison,chen2017entropy} combined with summation-by-parts (SBP) operators \cite{svard2014review,fernandez2014review}. While many numerical schemes, including finite differences \cite{kreiss1974finite,strand1994summation}, finite volumes \cite{nordstrom2001finite,nordstrom2003finite}, continuous finite elements \cite{hicken2016multidimensional,hicken2020entropy,abgrall2020analysisI}, flux reconstruction \cite{huynh2007flux,vincent2011newclass,ranocha2016summation}, active flux \cite{barsukow2025stability}, and meshless methods \cite{hicken2024constructing} can be formulated using SBP operators, we focus specifically on discontinuous Galerkin (DG) methods \cite{gassner2013skew,carpenter2014entropy,chan2018discretely}. Finally, we combine the entropy-stable baseline schemes with shock capturing using finite volume subcells \cite{hennemann2021provably} and adaptive mesh refinement implemented in the open-source software Trixi.jl \cite{ranocha2022adaptive,schlottkelakemper2021purely}. In particular, we provide all Julia source code and data needed to reproduce the numerical results presented in this paper in the accompanying reproducibility repository \cite{thein2025computingRepro}.

The outline of this paper is as follows. In Section~\ref{sec:eqns}, the governing equations relevant to this study are provided. In Section~\ref{sec:deriv_consflux}, an entropy-conservative numerical flux is derived.
Section~\ref{sec:num_meth} summarizes the numerical methods used. In Section~\ref{sec:num_res}, several numerical tests are provided, and the simulation results are presented. Concluding remarks are found in Section~\ref{sec:concl}. Appendix~\ref{app:entropy} offers detailed and supplementary calculations for the readers' convenience.
\section{Conservative Formulations of the Equations}\label{sec:eqns}
In this section we specify the governing equations precisely and provide necessary details on the conservative, primitive as well as entropy variables.
Here we consider the flat Minkowski metric given by
\begin{equation}
    g_{\alpha \beta} = \begin{cases}
        +1, &\alpha = \beta = 0,\\
        -1, &\alpha = \beta = 1,2,3,\\
        \hphantom{-}0, &\alpha \neq \beta.
    \end{cases}
\end{equation}
It is of further interest to study other space--times which then would result in curvi--linear coordinates, see for example \cite{Banyuls1997,DelZanna2007}. A precise study of such cases is beyond the scope of this work and will be considered in the future.
As mentioned before we study the ultra-relativistic equations for a perfect fluid in Minkowski space-time $\mathbf{x}=(x_1,x_2,x_3)$, $t=x_0$, namely
\begin{equation}\label{div}
    \sum_{\beta=0}^{3}\frac{\partial T_{\alpha \beta}}{\partial x_{\beta}} = 0
\end{equation}
with $\alpha,\beta\in\{0,1,2,3\}$ and the energy-momentum tensor
\begin{equation}
    T_{\alpha\beta} =-p g_{\alpha\beta}+4pu_\alpha u_\beta.
\end{equation}
Here $p$ represents the pressure, $\mathbf{u} \in \R^3$ is the spatial part of the four-velocity vector $(u_0,u_1,u_2,u_3)=(\sqrt{1+|\mathbf{u}|^2},\mathbf{u})$.
Note that in this work we consider the dimensionless velocity $\mathbf{u}$ and we have the following relation to the usual velocity $\mathbf{v}$
\[
  \mathbf{u} = \frac{\Gamma}{c}\mathbf{v}
\]
where $c$ is the speed of light and we have the Lorentz factor
\[
  \Gamma = \dfrac{1}{\sqrt{1 - \dfrac{|\mathbf{v}|^2}{c^2}}}.
\]
Using the dimensionless velocity we hence have $\Gamma = \sqrt{1 + |\mathbf{u}|^2}$ and we further scale to $c=1$.
For the ideal ultra-relativistic gas the thermal energy dominates and thus the local energy density $e$ at rest is linked to the pressure $p$ by
\begin{align}
    e = 3p.\label{ultra_rel_gas_eos}
\end{align}
Note that \eqref{ultra_rel_gas_eos} is universal for an ideal gas in the ultra-relativistic regime.
For the physical background we refer to Weinberg \cite[Part I, pp 47-52]{Weinberg1972}, further details can be found in Kunik \cite[Chapter 3.9]{Kunik2005} and in the works by Ruggeri and co-authors \cite{Pennisi2019,Ruggeri2021b,Ruggeri2021}.\\
The conservative multi--dimensional $d + 1$ form  of the ultra--relativistic Euler equations \eqref{div} is presented in the following.
We consider either $d = 2$ or $d = 3$ space dimensions. Then the unknown quantities $p$ and $\mathbf{u} = (u_1,\ldots,u_d) \in \R^d$
satisfying \eqref{div}
depend on time $t\geq 0$ and position $\mathbf{x} = (x_1,\ldots,x_d) \in \R^d$.\\
Putting $\alpha = 0$ in \eqref{div} gives the conservation of energy
\begin{equation}\label{energy_general}
    \frac{\partial}{\partial t}\left(3p + 4p|\mathbf{u}|^2 \right) + \sum\limits_{k=1}^d\frac{\partial}{\partial x_k}\left(4pu_k\sqrt{1 + |\mathbf{u}|^2} \right)=0,
\end{equation}
whereas for $\alpha = j = 1,\dots,d$ we obtain the conservation of momentum
\begin{equation}\label{momentum_general}
    \frac{\partial}{\partial t}\left(4pu_j \sqrt{1 + |\mathbf{u}|^2} \right) + \sum\limits_{k=1}^d\frac{\partial}{\partial x_k}\left(p\delta_{jk} + 4pu_j u_k\right)=0,
    \quad j = 1,\dots,d.
\end{equation}
We want to rewrite the conserved variables $\mathbf{w}:= (w_1,\ldots,w_d, w_{d+1})=(\overline{\mathbf{w}},w_{d+1})\in\R^{d+1}$, $\overline{\mathbf{w}}:=(w_1,\ldots,w_d)\in\R^d$
in terms of the primitive variables  determined by the velocity $\mathbf{u}:=(u_1,\ldots,u_d)\in\R^d$ and the pressure $p$, i.e.,
\begin{subequations}
    \label{eq:cons-var}
    \begin{align}
        \label{eq:cons-var-a}
        &w_j  := 4 p u_j \sqrt{1 + |\mathbf{u}|^2},\ j=1,\ldots,d, \\
        \label{eq:cons-var-b}
        & w_{d+1}  := p( 3 + 4 |\mathbf{u}|^2).
    \end{align}
\end{subequations}
Reversely, the primitive variables can be rewritten in terms of the conserved variables as
\begin{subequations}
    \label{eq:prim-var}
    \begin{align}
        & u_j = \frac{w_j}{\sqrt{4p(w_{d+1}+p)}},\ j=1,\ldots, d,\label{eq:velo}\\
        & p = \frac{1}{3} \left(-w_{d+1} + \sqrt{-3 |\overline{\mathbf{w}}|^2 + 4 w_{d+1}^2} \right) \ge 0.\label{eq:press}
    \end{align}
\end{subequations}
Concerning the entropy for this system we want to remark the following. First, it should be noted that the equation for the particle number density $n$ decouples from the system in the ultra--relativistic regime, cf. \cite{Ruggeri2021b,Weinberg1972,Kunik2005}. This provides the opportunity to find an entropy, additionally to the one from the full system, particularly for the system under consideration.
According to \cite{Kunik2005} the system is equipped with a concave physical entropy given by
\begin{align}
	\eta = p^{\frac{3}{4}}\sqrt{1 + |\mathbf{u}|^2}\label{def:entropy}
\end{align}
and the corresponding flux reads
\begin{align}
	\mathbf{q} = p^{\frac{3}{4}}\mathbf{u}.\label{def:entropy_flux}
\end{align}
To justify this entropy one can use arguments from \cite{Dreyer1998} as well as from Weinberg \cite[Chap. 10, p. 51]{Weinberg1972}.
In order to be self--contained we provide the proof that the presented entropy is a mathematical entropy for this system in the Appendix~\ref{app:entropy}. In particular we show concavity for the physical entropy and the relation of the entropy and the entropy flux, i.e.,
\begin{align}
	\nabla_\mathbf{w}q_k^T = \nabla_\mathbf{w}\eta^T\cdot\mathbf{D}_\mathbf{w}\mathbf{F}_k\label{entropy_flux_rel}
\end{align}
where $\mathbf{D}_\mathbf{w}\mathbf{F}_k$ denotes the Jacobian of the flux with respect to the conservative variables and $\mathbf{F}_k$ denotes the flux of the system in $x_k$ direction for $k=1,\dots,d$.
The entropy variables $\mathbf{h} = (h_1,\dots,h_{d+1})^T$ which are also known as main field, see \cite{Ruggeri2021b}, are given by
\begin{align}
	\label{deriv_eta_cons}
      h_j := \partial_{w_j} \eta =
      \begin{cases}
       -\dfrac{1}{4}\dfrac{\eta}{p}\dfrac{u_j}{\sqrt{1 + |\mathbf{u}|^2}},\, &j=1,\ldots,d, \\[1.5em]
           \hphantom{-}\dfrac{1}{4}\dfrac{\eta}{p},\, &j=d+1.
       \end{cases}
\end{align}
The main field of the present system may be obtained from the main field of the relativistic Euler equations since the ultra--relativistic Euler equations are a subsystem of the latter. We refer the interested reader to \cite{Ruggeri2021,Ruggeri2021b} and to the recent survey \cite{Warnecke2024}.
To construct the entropy-stable numerical flux, cf.\ \cite{tadmor1987numerical,chen2017entropy}, we need to calculate the entropy - flux potential.
The corresponding flux is given by
\begin{align}
    \psi_k = \nabla_\mathbf{W}\eta\cdot\mathbf{F}_k - q_k = -\frac{1}{4}\eta\frac{u_k}{\sqrt{1 + |\mathbf{u}|^2}}\quad\Leftrightarrow\quad\bb{\psi} = -\frac{1}{4}\eta\frac{\mathbf{u}}{\sqrt{1 + |\mathbf{u}|^2}}.\label{eq:flux_potential}
\end{align}
%
%
\section{Derivation of an entropy-conservative numerical flux}\label{sec:deriv_consflux}
In this section we present the detailed derivation of the entropy-conservative numerical flux. Following the seminal works of Tadmor \cite{tadmor1987numerical,tadmor2003entropy}, the condition
for an entropy-conservative numerical flux reads
\begin{equation}
	\dbl\mathbf{h}\dbr \cdot \tilde{\mathbf{F}}_k - \dbl\psi_k\dbr = 0,
\end{equation}
where $\mathbf{h}$ are the entropy variables given by \eqref{deriv_eta_cons}, $\tilde{\mathbf{F}}_k$ is the numerical flux in
space direction $k$, $\psi_k$ is the flux potential in direction $k$ defined in \eqref{eq:flux_potential},
and $\dbl a\dbr = a^+ - a^-$ denotes the jump of a quantity $a$.
The numerical fluxes $\tilde{\mathbf{F}}_k$ should be consistent approximations of the
fluxes
\begin{equation}
	\mathbf{F}_k = \begin{pmatrix}
		p\delta_{1k} + 4pu_1 u_k \\
		\vdots \\
		p\delta_{dk} + 4pu_d u_k \\[1em]
    	4 p u_k \sqrt{1 + |\mathbf{u}|^2}
  \end{pmatrix}
  =
  \begin{pmatrix}
		p\delta_{1k} + \dfrac{w_1w_k}{w_{d+1} + p} \\
		\vdots \\
		p\delta_{dk} + \dfrac{w_dw_k}{w_{d+1} + p} \\[1em]
    	w_k
  \end{pmatrix}\label{cons_sys_flux}
\end{equation}
of the ultra-relativistic Euler equations \eqref{energy_general} and \eqref{momentum_general} which form a conservative system
\begin{align*}
	\partial_t \mathbf{w} + \sum_{k=1}^d\partial_{x_k} \mathbf{F}_k(\mathbf{w}) = 0.
\end{align*}
We follow the algorithmic approach of \cite{ranocha2018comparison} to
derive an entropy-conservative numerical flux for the ultra-relativistic
Euler equations. Thus, we choose a set of variables and express all
jumps in terms of these variables. We then solve the resulting system
of equations in order to obtain entropy-conservative numerical fluxes which are ``affordable'' in the sense of \cite{ismail2009affordable}.

To express jumps of general expressions, we use mean values such as the
arithmetic mean
\begin{equation}
	\mean{a} = \frac{a^+ + a^-}{2}.
\end{equation}
Applying straight-forward algebraic manipulations we obtain the analogues of well-known rules for differentiation, i.e.,
\begin{align}
	\dbl a b\dbr &= \mean{a} \dbl b\dbr + \dbl a\dbr \mean{b},\\
	\dbl a^2\dbr &= 2 \mean{a} \dbl a\dbr,
\end{align}
and for positive $a$ we also have
\begin{equation}
	\dbl\sqrt{a}\dbr = \frac{\dbl a\dbr}{2 \mean{\sqrt{a}}}.\label{eq:jump-sqrt}
\end{equation}
Moreover, we will use
\begin{equation}
	-(a_- \sqrt{a_+} + a_+ \sqrt{a_-}) \dbl a^{-1/2}\dbr
  	=
  	-(a_- + \sqrt{a_+} \sqrt{a_-} - \sqrt{a_+} \sqrt{a_-} - a_+)
  	=
  	\dbl a\dbr,
\end{equation}
i.e.,
\begin{equation}
	\label{eq:jump-a-1/2-over-jump-a}
  	\frac{\dbl a^{-1/2}\dbr}{\dbl a\dbr}
  	=
  	-\frac{1}{a_- \sqrt{a_+} + a_+ \sqrt{a_-}}.
\end{equation}
Because of the structure of the entropy variables $\mathbf{h}$, we choose the
variables $u_i p^{-1/4}$ and $p$ to express all jumps. Then, we get
\begin{equation}
	\dbl h_i\dbr = -\frac{1}{4} \jump{u_i p^{-1/4}}, \quad i \in \{1, \dots, d\}.
\end{equation}
The jump of the entropy variable $h_{d+1}$ is given by
\begin{equation}
	\begin{aligned}
		\dbl h_{d+1}\dbr
		&=
  		\frac{1}{4} \dbl p^{-1/4} \sqrt{1 + |\mathbf{u}|^2}\dbr
  		=
  		\frac{1}{4} \dbl\sqrt{p^{-1/2} + \left(|\mathbf{u}| p^{-1/4}\right)^2}\dbr
  		=
  		\frac{1}{4} \frac{\dbl p^{-1/2} + \left(|\mathbf{u}| p^{-1/4}\right)^2\dbr}{2 \mean{\sqrt{p^{-1/2} + \left(|\mathbf{u}| p^{-1/4}\right)^2}}}
  		\\
  		&=
  		\frac{1}{8 \mean{\sqrt{p^{-1/2} + \left(|\mathbf{u}| p^{-1/4}\right)^2}}} \left(
    	\frac{\dbl p^{-1/2}\dbr}{\dbl p\dbr} \dbl p\dbr + 2 \mean{u_i p^{-1/4}} \dbl u_i p^{-1/4}\dbr
  		\right),
  	\end{aligned}
\end{equation}
where \eqref{eq:jump-sqrt} has been used. Moreover, the ratio
of the jumps of the pressure terms can be expressed using \eqref{eq:jump-a-1/2-over-jump-a}.
For the jump of the flux potential we get
\begin{equation}
	\dbl\psi_k\dbr
	=
	-\frac{1}{4} \dbl p^{3/4} u_k\dbr
	=
	-\frac{1}{4} \dbl p \, u_k p^{-1/4}\dbr
	=
	-\frac{1}{4} \mean{p} \dbl u_k p^{-1/4}\dbr
	-\frac{1}{4} \dbl p\dbr \mean{u_k p^{-1/4}}.
\end{equation}
Thus, the condition for the entropy-conservative numerical flux in the present situation reads
\begin{equation}
	\begin{aligned}
		0
		=
		\dbl\mathbf{h}\dbr \cdot \tilde{\mathbf{F}}_k - \dbl\psi_k\dbr
		&=
		\left(
		-\frac{1}{4} \tilde{\mathbf{F}}_k^{(i)}
		+ \frac{\mean{u_i p^{-1/4}}}{4 \mean{\sqrt{p^{-1/2} + \left(|\mathbf{u}| p^{-1/4}\right)^2}}} \tilde{\mathbf{F}}_k^{(d+1)}
		+ \frac{1}{4} \mean{p} \delta_{ik}
		\right) \jump{u_i p^{-1/4}}
		\\
		&\quad
		+ \left(
		\frac{1}{8 \mean{\sqrt{p^{-1/2} + \left(|\mathbf{u}| p^{-1/4}\right)^2}}} \frac{\dbl p^{-1/2}\dbr}{\dbl p\dbr} \tilde{\mathbf{F}}_k^{(d+1)}
		+ \frac{1}{4} \mean{u_k p^{-1/4}}
		\right) \dbl p\dbr.
	\end{aligned}
\end{equation}
This equation is satisfied if we choose the numerical flux $\tilde{\mathbf{F}}_k$ so that
each term in parentheses vanishes. First, we obtain the numerical energy flux
\begin{equation}
	\begin{aligned}
		\tilde{\mathbf{F}}_k^{(d+1)}
		&=
		-2 \frac{\dbl p\dbr}{\dbl p^{-1/2}\dbr} \mean{p^{-1/4} \sqrt{1 + |\mathbf{u}|^2}} \mean{u_k p^{-1/4}}
		\\
		&=
		2 \left( p_- \sqrt{p_+} + p_+ \sqrt{p_-} \right)
		\mean{p^{-1/4} \sqrt{1 + |\mathbf{u}|^2}} \mean{u_k p^{-1/4}}.
	\end{aligned}
\end{equation}
This is clearly a consistent approximation of the energy flux
$\mathbf{F}_k^{(d+1)} = 4 p u_k \sqrt{1 + |\mathbf{u}|^2}$. Inserting this numerical energy flux
into the other term in parentheses, we get
\begin{equation}
	\begin{aligned}
		&\quad
		-\frac{1}{4} \tilde{\mathbf{F}}_k^{(i)}
		+ \frac{\mean{u_i p^{-1/4}}}{4 \mean{\sqrt{p^{-1/2} + \left(|\mathbf{u}| p^{-1/4}\right)^2}}} \tilde{\mathbf{F}}_k^{(d+1)}
		+ \frac{1}{4} \mean{p} \delta_{ik}
		\\
		&=
		-\frac{1}{4} \tilde{\mathbf{F}}_k^{(i)}
		- \frac{1}{2} \frac{\dbl p\dbr}{\dbl p^{-1/2}\dbr} \mean{u_k p^{-1/4}} \mean{u_i p^{-1/4}} + \frac{1}{4} \mean{p} \delta_{ik}.
	\end{aligned}
\end{equation}
This term vanishes if we choose the numerical momentum fluxes
\begin{equation}
	\begin{aligned}
		\tilde{\mathbf{F}}_k^{(i)}
		&=
		-2 \frac{\dbl p\dbr}{\dbl p^{-1/2}\dbr} \mean{u_i p^{-1/4}} \mean{u_k p^{-1/4}}
		+ \mean{p} \delta_{ik}
		\\
		&=
		2 \left( p_- \sqrt{p_+} + p_+ \sqrt{p_-} \right) \mean{u_i p^{-1/4}} \mean{u_k p^{-1/4}}
		+ \mean{p} \delta_{ik}.
	\end{aligned}
\end{equation}
This is a consistent approximation of the momentum flux
$\mathbf{F}_k^{(i)} = 4 p u_i u_k + \delta_{ik} p$ of the ultra-relativistic Euler equations.
Thus, we have proven
\begin{theorem}\label{thm:entropy-conservative-flux}
	The numerical fluxes
	\begin{equation}
		\label{eq:entropy-conservative-flux}
		\begin{aligned}
			\tilde{\mathbf{F}}_k^{(i)}
			&=
			2 \left( p_- \sqrt{p_+} + p_+ \sqrt{p_-} \right) \mean{u_i p^{-1/4}} \mean{u_k p^{-1/4}}
			+ \mean{p} \delta_{ik},
			\\
			\tilde{\mathbf{F}}_k^{(d+1)}
			&=
			2 \left( p_- \sqrt{p_+} + p_+ \sqrt{p_-} \right)
			\mean{p^{-1/4} \sqrt{1 + |\mathbf{u}|^2}} \mean{u_k p^{-1/4}}
		\end{aligned}
	\end{equation}
  	are symmetric, consistent, and entropy-conservative two-point fluxes
  	for the ultra-relativistic Euler equations \eqref{momentum_general} and \eqref{energy_general}.
\end{theorem}

\begin{remark}
	The pressure in the momentum flux of \eqref{eq:entropy-conservative-flux}
	is approximated using an arithmetic mean. For the classical compressible
	Euler equations, this is required for the stronger version of
	kinetic energy preservation \cite{jameson2008formulation,chandrashekar2013kinetic}
	discussed in \cite{ranocha2020entropy,ranocha2018thesis},
	see also \cite{ranocha2021preventing,ranocha2022note}.
    Even without provable entropic properties, kinetic energy-preserving split forms have been widely used for compressible computational fluid dynamics and turbulence simulations because of the increased robustness \cite{gassner2016split,sjogreen2017skew,klose2020assessing}.
\end{remark}
\section{Numerical Methods}\label{sec:num_meth}
We want to compare numerical results for radially symmetric solutions of the original multi--dimensional initial-value problem (IVP) for the ultra-relativistic Euler equations \eqref{energy_general}--\eqref{momentum_general}. To avoid confusion of the pressure in the multi--dimensional case and the pressure in the radially symmetric case we denote the pressure in the multi--dimensional case by $\tilde{p}$ from now on.
For $j=1,\dots,d$ the IVP reads
\begin{equation}\label{orisys}
    \left\{
    \begin{aligned}
        \frac{\partial}{\partial t}\left(3\tilde{p} + 4\tilde{p}|\mathbf{u}|^2 \right)
        + \sum\limits_{k=1}^d\frac{\partial}{\partial x_k}\left(4\tilde{p}u_k\sqrt{1 + |\mathbf{u}|^2} \right) &= 0,\\
        \frac{\partial}{\partial t}\left(4\tilde{p}u_j \sqrt{1 + |\mathbf{u}|^2} \right)
        + \sum\limits_{k=1}^d\frac{\partial}{\partial x_k}\left(\tilde{p}\delta_{jk} + 4\tilde{p}u_j u_k\right) &= 0,\\
        \tilde{p}(0,\mathbf{x}) = \tilde{p}_0(|\mathbf{x}|), \quad \mathbf{u}(0,\mathbf{x}) = \mathbf{u}_0(\mathbf{x}) &= u_0(|\mathbf{x}|) \frac{\mathbf{x}}{|\mathbf{x}|}.
    \end{aligned}\right.
\end{equation}
We require radially symmetric solutions of this system, i.e., for $t \geq 0$ and $x=|\mathbf{x}|>0$ the restrictions for pressure and velocity are given by
\begin{align}
    \tilde{p}(t,\mathbf{x})=p(t,x)\,,\quad \mathbf{u}(t,\mathbf{x})=\frac{v(t,x)}{\sqrt{1 - v(t,x)^2}} \cdot \frac{\mathbf{x}}{x}\,.
\end{align}
For $t=0$ we obtain given initial data $\tilde{p}(0,\mathbf{x}) =p(0,x)=p_0(x)$ and
\begin{align}
    \mathbf{u}(0,\mathbf{x}) = \mathbf{u}_0(\mathbf{x}) = \frac{v_0(x)}{\sqrt{1 - v_0(x)^2}} \cdot \frac{\mathbf{x}}{x}\,,
\end{align}
and for $\mathbf{x}=\mathbf{0} \in \R^d$ we may also put $\mathbf{u}(0,\mathbf{0})=\mathbf{0}$ and $v_0(0)=0$.

\subsection{Multi--dimensional discontinuous Galerkin methods}
\label{sec:dg}

We use the numerical simulation framework Trixi.jl \cite{ranocha2022adaptive,schlottkelakemper2021purely} implemented in Julia \cite{bezanson2017julia} to solve the multi--dimensional IVP \eqref{orisys} using discontinuous Galerkin methods.
We follow the standard approach to first partition the spatial domain $\Omega$ into quadrilateral/hexahedral elements with disjoint interiors.
In each element, the numerical solution is represented by its values at tensor product Gauss-Lobatto-Legendre nodes able to represent polynomials of degree three in each coordinate direction.
The associated quadrature rule and differentiation matrices are used to compute integrals and derivatives, respectively, i.e., we use a collocation approach.
Thus, the resulting scheme is a high-order accurate nodal discontinuous Galerkin spectral element method (DGSEM) \cite{kopriva2009implementing}.

DG schemes are often derived from a weak form of the equations.
First, multiply the conservation law $\partial_t \mathbf{w} + \sum_{k=1}^d\partial_{x_k} \mathbf{F}_k(\mathbf{w}) = 0$ by a test function $\phi$ and integrate over an element $E$ to get
\begin{align*}
    \int_E \phi \partial_t \mathbf{w} \, d\mathbf{x} + \sum_{k=1}^d \int_E \phi \partial_{x_k} \mathbf{F}_k(\mathbf{w}) \, d\mathbf{x} = 0.
\end{align*}
Next, integrate the volume term by parts to get
\begin{align*}
    \int_E \phi \partial_t \mathbf{w} \, d\mathbf{x} - \sum_{k=1}^d \int_E (\partial_{x_k} \phi) \mathbf{F}_k(\mathbf{w}) \, d\mathbf{x} + \sum_{k=1}^d \int_{\partial E} \phi \mathbf{F}_k(\mathbf{w}) n_k \, dS = 0,
\end{align*}
where $n_k$ is the $k$-th component of the outward normal vector on
the boundary $\partial E$ of the element $E$.
We replace the flux in the surface integral by a numerical flux $\hat{\mathbf{F}}$ to get
\begin{align*}
    \int_E \phi \partial_t \mathbf{w} \, d\mathbf{x} - \sum_{k=1}^d \int_E (\partial_{x_k} \phi) \mathbf{F}_k(\mathbf{w}) \, d\mathbf{x} + \int_{\partial E} \phi \hat{\mathbf{F}}(\mathbf{w}_-, \mathbf{w}_+, \mathbf{n}) \, dS = 0,
\end{align*}
where $\mathbf{w}_-$ and $\mathbf{w}_+$ are the values of the solution on the interior and exterior of the element, respectively, and $\mathbf{n}$ is the outward normal vector on the boundary $\partial E$ of the element $E$.
Concretely, we use a local Lax-Friedrichs/Rusanov flux $\hat{\mathbf{F}}$.

Flux differencing \cite{lefloch2002fully,fisher2013high,ranocha2018comparison,chen2017entropy} is more similar to the strong form DG scheme.
Thus, we integrate the volume term once more by parts to get
\begin{align*}
    \int_E \phi \partial_t \mathbf{w} \, d\mathbf{x} + \sum_{k=1}^d \int_E \phi \partial_{x_k} \mathbf{F}_k(\mathbf{w}) \, d\mathbf{x} + \int_{\partial E} \phi \left( \hat{\mathbf{F}}(\mathbf{w}_-, \mathbf{w}_+, \mathbf{n}) - \mathbf{F}(\mathbf{w}_-) \cdot \mathbf{n} \right) \, dS = 0.
\end{align*}
A classical strong-form DGSEM is obtained by replacing the volume flux $\mathbf{F}_k(\mathbf{w})$ by the polynomial approximation obtained by collocation at the Gauss-Lobatto-Legendre nodes, replacing integrals by the associated quadrature rule \cite{kopriva2009implementing}, and requiring the resulting equations to hold for all test functions $\phi$ in the same polynomial space.
The resulting volume term can be written as
\begin{align*}
    \mathbf{VOL}_i = \sum_{k=1}^d \sum_j (D_k)_{i,j} \mathbf{F}_k(\mathbf{w}_j),
\end{align*}
where the indices $i$ and $j$ run over the Gauss-Lobatto-Legendre nodes in the element and $D_k$ is the discrete (polynomial) differentiation matrix in the $k$-th coordinate direction.

Flux differencing methods \cite{lefloch2002fully,fisher2013high,ranocha2018comparison,chen2017entropy} replace the volume term above by
\begin{align*}
    \mathbf{VOL}_i = \sum_{k=1}^d \sum_j 2 (D_k)_{i,j} \tilde{\mathbf{F}}_k(\mathbf{w}_i, \mathbf{w}_j),
\end{align*}
where $\tilde{\mathbf{F}}_k$ is a two-point numerical flux.
An efficient implementation of such flux differencing volume terms is described in \cite{ranocha2023efficient}.
We use the newly developed entropy-conservative two-point flux given in Theorem~\ref{thm:entropy-conservative-flux} for the volume term.
Thus, the flux differencing theory \cite{lefloch2002fully,fisher2013high,ranocha2018comparison,chen2017entropy} guarantees that we obtain an entropy-stable high-order method.

To ensure robustness also in the presence of discontinuities, we use the finite volume subcell shock capturing approach of \cite{hennemann2021provably}.
There, the high-order flux differencing volume terms are interpreted as high-order subcell finite volume methods.
For robustness, these high-order methods are blended with a first-order finite volume method using a local Lax-Friedrichs/Rusanov flux.
The blending factor determining the convex combination of high- and low-order methods is based on a shock indicator which is applied to the variable $\tilde{p} \sqrt{1 + |\mathbf{u}|^2}$.
For efficiency, we use adaptive mesh refinement (AMR) based on \texttt{p4est} \cite{burstedde2011p4est} and an adaptation of the indicator of \cite{lohner1987adaptive} applied to the same indicator variable used for shock capturing.

The resulting system of ordinary differential equations is integrated in time using the third-order, four-stage explicit strong stability preserving Runge-Kutta method of \cite{kraaijevanger1991contractivity} with embedded method of \cite{conde2022embedded} and adaptive step size controller of \cite{ranocha2021optimized} implemented in OrdinaryDiffEq.jl \cite{rackauckas2017differentialequations}.
Thus, the step size is adapted automatically to the maximal resolution of the mesh and the local wave speeds, see also \cite{ranocha2023error} for further details on the error control and step size adaptation.
For the self-similar expansion (Example~2 below), we also apply the positivity-preserving limiter of Zhang and Shu for the pressure \cite{zhang2011maximum}.

All data and the Julia source code required to reproduce the numerical results presented in this paper are available online in our reproducibility repository \cite{thein2025computingRepro}.

\subsection{One-dimensional radially symmetric solver RadSymS}
For a detailed presentation of the scheme \emph{RadSymS} used to compute the solution to the quasi one-dimensional radially symmetric form of the ultra-relativistic Euler equations  we refer to \cite{kunik2024radially} and \cite{Kunik2021}.
The solution of the radially symmetric quasi one-dimensional problem satisfies
\begin{equation}\label{radsys}
    \left\{
    \begin{aligned}
        \quad\quad\quad\frac{\partial}{\partial t}\left(x^{d-1}a \right) + \frac{\partial}{\partial x}\left(x^{d-1}b \right) &= 0,\\
        \quad\quad\quad\frac{\partial}{\partial t}\left(x^{d-1}b \right) + \frac{\partial}{\partial x}\left(x^{d-1}c \right) &= \frac{d-1}{2} x^{d-2} (a - c)
        %
    \end{aligned}\right.
\end{equation}
for $t>0$ and $x > 0$ representing the radial direction and with the initial conditions
\begin{align}
	\lim \limits_{t \searrow 0}a(t,x) = a_0(x),\quad \lim \limits_{t \searrow 0}b(t,x) &= b_0(x)\quad\text{with}\quad|b_0(x)| < a_0(x).\label{radsys_init_cond}
\end{align}
%
Here $d=2$ or $d=3$ denotes the space dimension and the variables $(a,b)$ are linked to $(p,u)$ through a one-to-one transformation on the sets
\[
  \mathcal{S} = \{(p,u)\in\R^2\,|\, p > 0\}\quad\text{and}\quad\tilde{\mathcal{S}} = \{(a,b)\in\R^2\,|\, |b| < a\}.
\]
The mapping $\Theta:\mathcal{S}\to\tilde{\mathcal{S}}$ and its inverse are given by
\begin{align}
	\Theta(p,u) &= \begin{pmatrix} p\left(3 + 4u^2\right)\\ 4pu\sqrt{1 + u^2}\end{pmatrix} = \begin{pmatrix} a\\ b\end{pmatrix},\label{statetrans}\\
	\Theta^{-1}(a,b) &= \begin{pmatrix} \dfrac{1}{3}\left(\sqrt{4a^2 - 3b^2} - a\right)\\ \dfrac{b}{\sqrt{4p(p + a)}}\end{pmatrix} = \begin{pmatrix} p\\ u\end{pmatrix}.\label{inverse_statetrans}
\end{align}
Using the transformation \eqref{statetrans}, its inverse \eqref{inverse_statetrans} and the velocity
\begin{align}
\label{eq:velocity_transformation}
    v = \frac{u}{\sqrt{1+u^2}}\quad\text{with}\;|v| < 1,
\end{align}
we replace the state variables $a$ and $b$ by $p$ and $v$, respectively.
The variable $c = c(a,b)$ is given by
\begin{align}
	c = c(a,b) = \frac{5}{3}a - \frac{2}{3}\sqrt{4a^2 - 3b^2}\label{cdefinition}.
\end{align}
Observe, that the variables $a$ and $b$ correspond to the conservative quantities in radial direction, while $c(a,b) = p(1 + 4u^2)$.

We are looking for weak solutions $a = a(t,x)$, $b = b(t,x)$ with $|b| < a$, i.e.
\begin{align}
	\int_{\partial\Omega}x^{d-1}a\,\dd x - x^{d-1}b\,\dd t = 0,\,\int_{\partial\Omega}x^{d-1}b\,\dd x - x^{d-1}c\,\dd t = \frac{d-1}{2}\iint_{\Omega}x^{d-2}(a - c)\,\dd t\dd x.\label{weak_contour}
\end{align}
Here $\Omega \subset Q$ is a convex domain with piecewise smooth boundary $\partial\Omega$ with $Q = \{(t,x)\,|\,t > 0,\;x > 0\}$ and the integral formulation is based on \cite{Oleinik1963,Kunik2005}. Note that the integrals along $\partial\Omega$ are line integrals over a 1-form written in their local coordinates, cf. \cite{Lee2013,Dineen2014}.

Based on this reformulation of the problem the one-dimensional numerical scheme (RadSymS) for the initial value problem of the radially symmetric ultra-relativistic Euler equations in two and three space dimensions is developed. 
It further can be shown that this scheme preserves positive pressure, see \cite{kunik2024radially,Kunik2021}.
The method of contour-integration for the formulation of the balance laws \eqref{weak_contour} is used to construct a function called ``Euler''. 
This function enables us to obtain the time evolution of the numerical solution on a staggered grid.
More precisely it allows us to construct the solution $(a',b')$ at the next time step
from the solution $(a_{\pm},b_{\pm})$ in two neighboring grid points at the previous time step according to Figure \ref{figeul}.
First we determine the computational domain and define some quantities which are needed for its discretization.
\begin{enumerate}[1)]
    \item Given are $t_*,x_*>0$ in order to calculate a numerical solution of the initial value problem \eqref{weak_contour}, with the initial conditions given in \eqref{radsys},
          in the time range $[0,t_*]$ and the spatial range $[0,x_*]$.
    \item We want to use a staggered grid scheme. Any given number $N \in \N$ with $N \cdot x_* \geq t_*$ determines the time step size $\Delta t = t_*/2N$ and the time steps are $t_{n} = (n-1) \Delta t, n=1,\ldots,2N+1$.
    \item Put
          \[
            M = \left \lfloor \frac{x_*}{t_*}\,N \right \rfloor \geq 1,
          \]
          then the spatial mesh size is $\Delta x= x_*/M$ with the spatial grid points $x_{j} = (j-1) \Delta x\,,  j=1,\ldots,N+M+1$.
          Note that our scheme uses a trapezoidal computational domain $\mathcal{D}$ defined below that includes the target domain $[0,t_*] \times [0,x_*]$.
          Thereby, we can use all initial data that influence the solution on the target domain. In this way we avoid using a numerical boundary condition at $x_*$.
    \item The number $\lambda= \Delta x /(2 \Delta t) \geq 1$ is used to satisfy the CFL-condition and to define the computational domain
          \[
          \mathcal{D} = \left\{(t,x)\in \R^2\,:\,0 \leq t \leq t_*,\quad 0 \leq x \leq x_*+\lambda(t_*-t)\right \}.
          \]
\end{enumerate}
For the numerical discretization of the integral balance laws \eqref{weak_contour} we choose the triangular balance domain $\Omega$ depicted in Figure \ref{figeul}.
We assume that the midpoints $P_-=(\overline{t}, \overline{x}-\Delta x/2)$, $P_+=(\overline{t}, \overline{x}+\Delta x/2)$
and $P'=(\overline{t}+\Delta t, \overline{x})$ of the cords of $\partial \Omega$ are numerical grid points for the computational domain $\mathcal{D}$.
Let the numerical solution $(a_\pm,b_\pm)$ be given at the grid points $P_{\pm}$.
We have to require $|b_{\pm}|< a_{\pm}$ for the numerical solution in the actual time step $\overline{t}=t_n$ with $n=1,\ldots,2N$.
The major task is to calculate the numerical solution $(a',b')$ for the next time step $\overline{t}+\Delta t= t_{n+1}$ at its grid point $P'$, see Figure \ref{figeul}.
\begin{figure}[h!]
    \centering
    \begin{tikzpicture}[scale=0.9]
        \draw[thick] (0,2) -- (0.6,2);
        \fill[black] (0.8,2) circle (2pt);
        \fill[black] (1,2) circle (2pt);
        \fill[black] (1.2,2) circle (2pt);
        \draw[thick] (1.4,2) -- (2,2);
        \draw[thick] (0,4) -- (0.6,4);
        \fill[black] (0.8,4) circle (2pt);
        \fill[black] (1,4) circle (2pt);
        \fill[black] (1.2,4) circle (2pt);
        \draw[thick] (1.4,4)-- (2,4);
        \draw[thick] (0,6) -- (0.6,6);
        \fill[black] (0.8,6) circle (2pt);
        \fill[black] (1,6)circle (2pt);
        \fill[black] (1.2,6) circle (2pt);
        \draw[thick] (1.4,6) -- (2,6);
        \node at (-0.2,-2.3) {$0$};
        \node at (-1,2) {$\overline{x}-\Delta x/2$};
        \fill[black] (4,6) circle (2pt);
        \node at (-1,6) {$\overline{x}+\Delta x/2$};
        \node at (-0.3,4) {$\overline{x}$};
        %
        \draw[thick] (0,-2) -- (0.6,-2);
        \fill[black] (0.8,-2) circle (2pt);
        \fill[black] (1,-2) circle (2pt);
        \fill[black] (1.2,-2) circle (2pt);
        \draw[very thick,-stealth] (1.4,-2) -- (7,-2);
        \node at (7,-2.3) {$t$};
        \draw[thick] (0,-2) -- (0,-1.4);
        \fill[black] (0,-1.2) circle (2pt);
        \fill[black] (0,-1) circle (2pt);
        \fill[black] (0,-0.8) circle (2pt);
        \draw[very thick,-stealth] (0,-0.6) -- (0,8.5);
        \node at (-0.3,8.4) {$x$};
        \draw[thick] (2,-2) -- (2,-1.4);
        \fill[black] (2,-1.2) circle (2pt);
        \fill[black] (2,-1) circle (2pt);
        \fill[black] (2,-0.8) circle (2pt);
        \draw[thick] (2,-0.6) -- (2,4);
        \draw[thick] (4,-2) -- (4,-1.4);
        \fill[black] (4,-1.2) circle (2pt);
        \fill[black] (4,-1) circle (2pt);
        \fill[black] (4,-0.8) circle (2pt);
        \draw[thick] (4,-0.6) -- (4,2);
        \draw[thick] (6,-2) -- (6,-1.4);
        \fill[black] (6,-1.2) circle (2pt);
        \fill[black] (6,-1) circle (2pt);
        \fill[black] (6,-0.8) circle (2pt);
        \draw[thick] (6,-0.6) -- (6,0);
        \draw[thick] (2,4) -- (6,4);
        \draw[thick] (2,2) -- (4,2);
        \draw[thick] (2,6) -- (4,6);
        \draw[thick] (4,2) -- (4,6);
        \node at (2,-2.3) {$\overline{t}-\Delta t$};
        \node at (4,-2.3) {$\overline{t}$};
        \node at (6,-2.3) {$\overline{t}+\Delta t$};
        \draw[thick] (2,4) -- (6,8);
        \draw[thick] (2,4) -- (6,0);
        \draw[thick] (6,0) -- (6,8);
        \draw[very thick,-stealth] (4,2) -- (5,1);
        \node at (4.8,0.8) {$\gamma_2$};
        \draw[thick,-stealth] (4,6) -- (3,5);
        \node at (2.72,5.2) {$\gamma_1$};
        \draw[very thick,-stealth] (6,4) -- (6,6);
        \node at (6.4,6) {$\gamma_0$};
        \node at (4.5,4.3) {$\Omega$};
        \node at (4.8,5.9) {$(a_+,b_+)$};
        \node at (4.8,2.1) {$(a_-,b_-)$};
        \node at (6.6,4) {$(a',b')$};
    \end{tikzpicture}
    \caption{The balance region $\Omega$}
    \label{figeul}
\end{figure}
The spatial value $\overline{x} \geq 0$ is given. We have to determine a function
\begin{equation}\label{euler_solve}
    \text{Euler}(a_-,b_-,a_+,b_+,\overline{x},\Delta x, \lambda)=(a',b')
\end{equation}
for the calculation of $(a',b')$. This leads to the structure of a staggered grid scheme.
Note that at the boundary the balance region $\Omega$ may have parts outside $\mathcal{D}$, e.g.~points below the half-space $x\geq 0$.
In the latter case we employ a simple reflection principle for the numerical solution in order to use the function $\text{Euler}$ as well for the boundary points with $\overline{x}=0$.\\
Next we make use of the fact that the points $P_{\pm}$ with the numerical values $(a_{\pm},b_{\pm})$ and $P'$ with the unknown value $(a',b')$ are the 
\textit{midpoints of the three boundary cords} of the balance region $\Omega$.
\begin{theorem}[Numerical solution $(a',b')$ for the balance region $\Omega$]\label{apbp} 
    Given are real quantities $\overline{x}\geq 0$ and $a_{\pm}$, $b_{\pm}$. Assume that $|b_{\pm}|< a_{\pm}$.
    We recall $\lambda\geq 1$ defined in terms of $\Delta t$ and $\Delta x$ and put $c_{\pm}=c(a_{\pm},b_{\pm})$ in  \eqref{cdefinition}.
    \begin{align*}
    	a' &=
    	\begin{cases}
    		\dfrac12\left(a_- +\dfrac{b_-}{\lambda}\right)\left(1-\dfrac{\tilde{q}}{d-1}\right) + \dfrac12\left(a_+ - \dfrac{b_+}{\lambda}\right)\left(1+\dfrac{\tilde{q}}{d-1}\right),\,&\overline{x} > 0,\\
    		a_+-\dfrac{b_+}{\lambda},\,&\overline{x} = 0
    	\end{cases}\\
    	b' &=
    	\begin{cases}
    		\dfrac{\xi + \nu\sqrt{4a'^2(1+3\nu^2)-3\xi^2}}{1+3\nu^2},\,&\overline{x} > 0,\\
    		0,\,&\overline{x} = 0
    	\end{cases}\\
    	\text{with}&\\
    	\tilde{q} &= \dfrac{(d-1)\overline{x}^{d-2}\Delta x}{2\left(\overline{x}^{d-1} + \frac{d-2}{3}(\Delta x)^2\right)} < 1,\quad
    	\nu = \frac{\tilde{q}}{3(d-1)\lambda}\\
    	\text{and}\;\xi &= \frac12\left(b_- +\frac{c_-}{\lambda}\right)\left(1-\frac{\tilde{q}}{d-1}\right) + \frac12\left(b_+ - \frac{c_+}{\lambda}\right)\left(1+\frac{\tilde{q}}{d-1}\right)-a' \nu
    \end{align*}
    Then we have $|b'|<a'$ in both cases.
\end{theorem}
\begin{definition}[The function Euler]\label{eulerfunction}
    The state $(a',b')$ from Theorem \ref{apbp} defines the function $\text{Euler}$ in \eqref{euler_solve}.
\end{definition}
Now we formulate the numerical scheme for the solution of the initial-boundary value problem \eqref{weak_contour}.
We construct staggered grid points in the computational domain $\mathcal{D}$ and compute the numerical solution at these grid points.
Using the function $\emph{Euler}$ we obtain the evolution of the numerical solution in time,
i.e., it allows us to construct the solution at time $t=t_{n+1}$ from the solution which is already calculated in the grid points   at the former time step $t=t_n$.
\begin{enumerate}[(I)]
    \item The staggered grid points are $(t_n,x_{n,j}) \in \mathcal{D}$ for $t_n=(n-1)\Delta t$, $n=1,\ldots,2N+1$ and $j=1,\ldots,M+N-\lfloor (n-1)/2 \rfloor$ with
          \begin{align*}
            x_{n,j} = 
            \begin{cases}
                (x_{j}+x_{j+1})/2\; &\text{if $n$ is odd}\\
                x_j\; &\text{if $n$ is even}.
            \end{cases}
          \end{align*}
          We want to calculate the numerical solution $(a_{n,j},b_{n,j})$ at $(t_n,x_{n,j})$.
    \item For $j = 1,\ldots,M+N$ we calculate the numerical solution $(a_{1,j}, b_{1,j})$ at the grid point $(t_1,x_{1,j})=(0,(x_j+x_{j+1})/2)$ from the given initial data by
          \[
            a_{1,j} = a_0\left(x_{1,j}\right),\quad b_{1,j} = b_0\left(x_{1,j}\right).
          \]
          This corresponds to taking the integral average of the initial data on $(x_j,x_{j+1})$ and using the midpoint rule as quadrature.
    \item Assume that for a fixed \textit{odd index} $n\in \{1,\ldots,2N\}$ we have already determined the numerical solution $(a_{n,j}, b_{n,j})$
          at the grid points $(t_n,x_{n,j})$, $j=1,\ldots,M + N - (n - 1)/2$.\\
          First we determine the solution $(a_{n+1,1},b_{n+1,1})$ at the boundary point $(t_{n+1},x_{n+1,1}) = (t_{n+1},0)$ according to Thm.\ \ref{apbp}.
          For this purpose we put $a_+ = a_{n,1}$, $b_+ = b_{n,1}$, $a_- = a_{n,1}$, $b_- = -b_{n,1}$ and have
          \[
            (a_{n+1,1}, b_{n+1,1}) = \text{Euler}(a_-,b_-,a_+,b_+,0,\Delta x,\lambda)\;\text{with}\; b_{n+1,1} = 0.
          \]
          Next we put $a_- = a_{n,j-1}$, $b_- = b_{n,j-1}$ and $a_+ = a_{n,j}$, $b_+ = b_{n,j}$ for $j = 2,\ldots,M + N  - (n - 1)/2$ and determine the values $a_{n+1,j}$, $b_{n+1,j}$ 
          at time $t_{n+1}$ and position $\overline{x} = x_{n+1,j} = x_{j}$ from
          \[
            (a_{n+1,j}, b_{n+1,j})=\text{Euler}(a_-,b_-,a_+,b_+,\overline{x},\Delta x,\lambda).
          \]
    \item Assume that for a fixed \textit{even index} $n \in \{1,\ldots,2N\}$ we have already determined the numerical solution $(a_{n,j}, b_{n,j})$
    at the grid points $(t_n,x_{n,j})$, $j=1,\ldots,M + N - n/2 + 1$.\\
    We put $a_- = a_{n,j}$, $b_- = b_{n,j}$ and $a_+ = a_{n,j+1}$, $b_+ = b_{n,j+1}$ for $j = 1,\ldots,M + N - n/2$ and determine the values $a_{n+1,j}$, $ b_{n+1,j}$ at time $t_{n+1}$ and position 
    $\overline{x} = x_{n+1,j} = (x_{j} + x_{j+1})/2$ from
    \[
      (a_{n+1,j}, b_{n+1,j}) = \text{Euler}(a_-,b_-,a_+,b_+,\overline{x},\Delta x,\lambda).
    \]
\end{enumerate}
The quasi one-dimensional problem \eqref{radsys} is approximately solved by the one-dimensional radially symmetric scheme with $N=5000$.
We prescribe the initial pressure $p_0=p(0,\cdot)$ as well as the initial velocity $v_0=v(0,\cdot)$.
The variable $v$ is chosen, because the restriction $|v| < 1$ leads to better color plots (the variable $u$ is unbounded).
Note that for constant pressure $p_0 > 0$  and $v_0 = 0$, we obtain a stationary solution.
Such a steady part is contained in the solutions to the examples presented below.
\section{Numerical Results}\label{sec:num_res}
\begin{figure}[h!]
    \hspace{-1cm}
    \subfigure[Results for the pressure $p$ in 2D]{
    \includegraphics[width=0.575\textwidth]{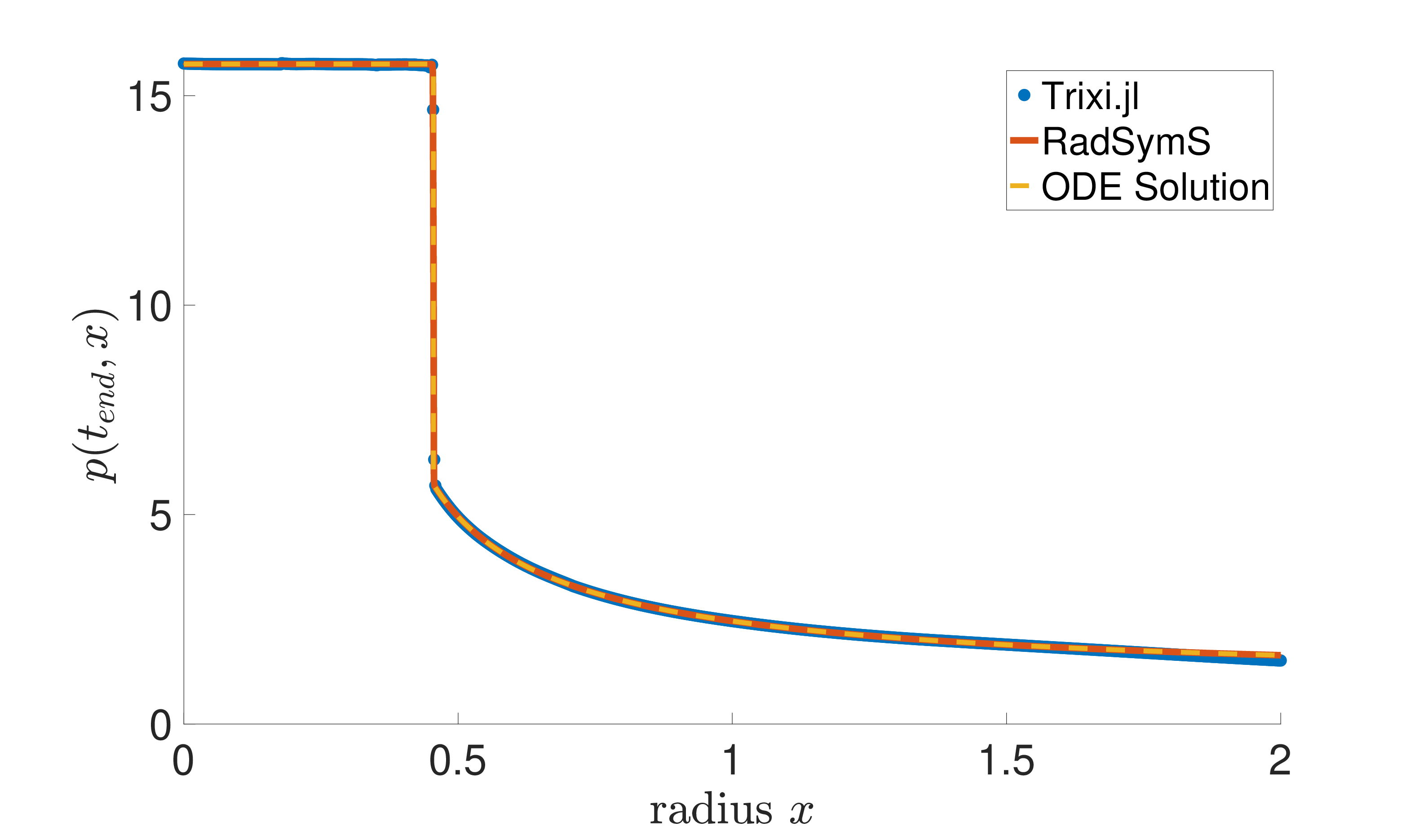}}
    \hspace{-1cm}
    \subfigure[Results for the velocity $v$ in 2D]{
    \includegraphics[width=0.575\textwidth]{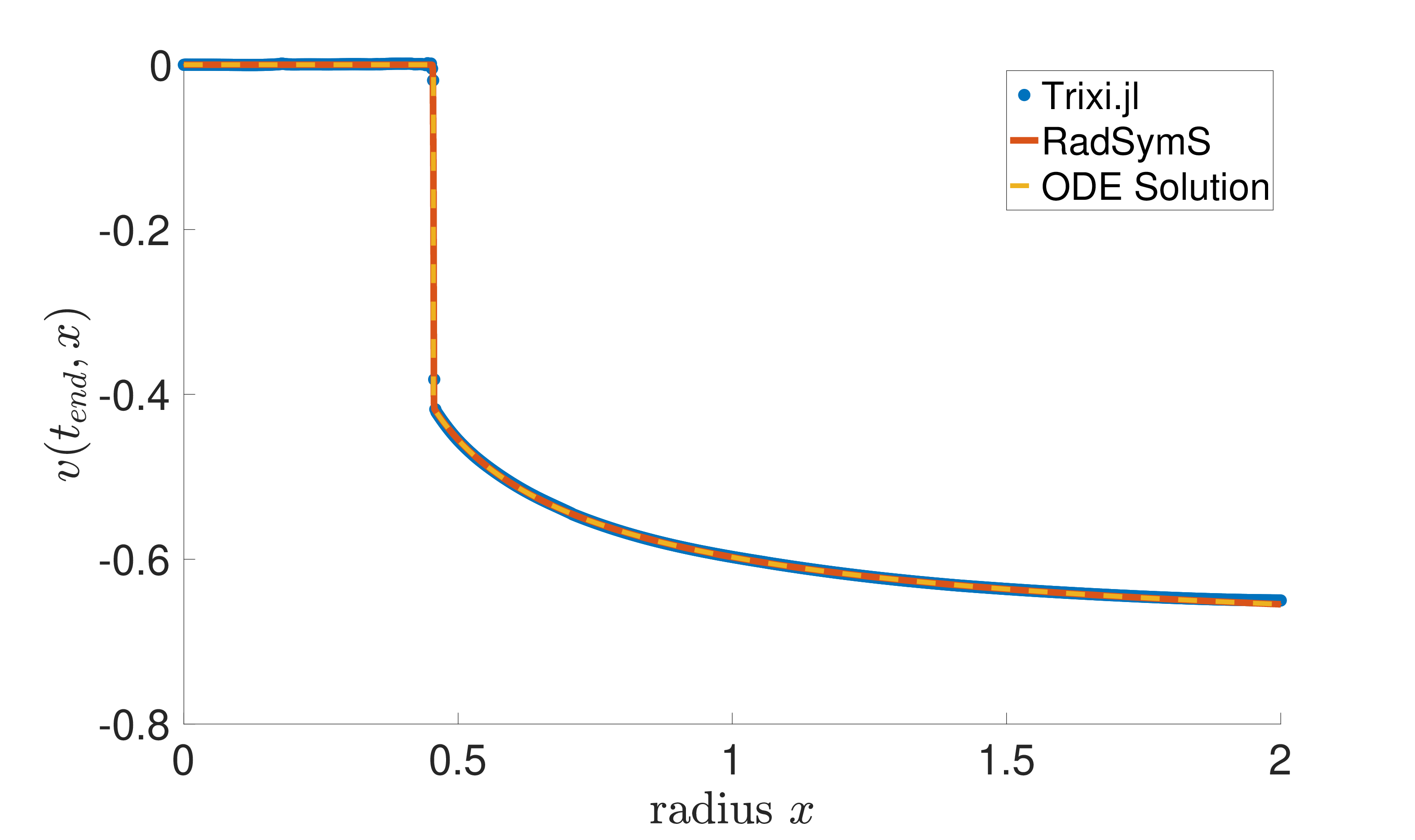}}
    \caption{   \textbf{Example 1:}
    Results for the pressure $p$ and velocity $v$ in 2D
     Comparison for Trixi.jl (blue) and RadSymS (red) in radial direction  with the solution of ODE system \eqref{GL_ODE_sys} (yellow) at final time $t_{end} = 1$ for the pressure $p$, see (a), and velocity $v$, see (b).
  }
    \label{fig:example-1-2D}
\end{figure}

\begin{figure}[h!]
    \hspace{-1cm}
    \subfigure[Results for the pressure $p$ in 3D]{
    \includegraphics[width=0.575\textwidth]{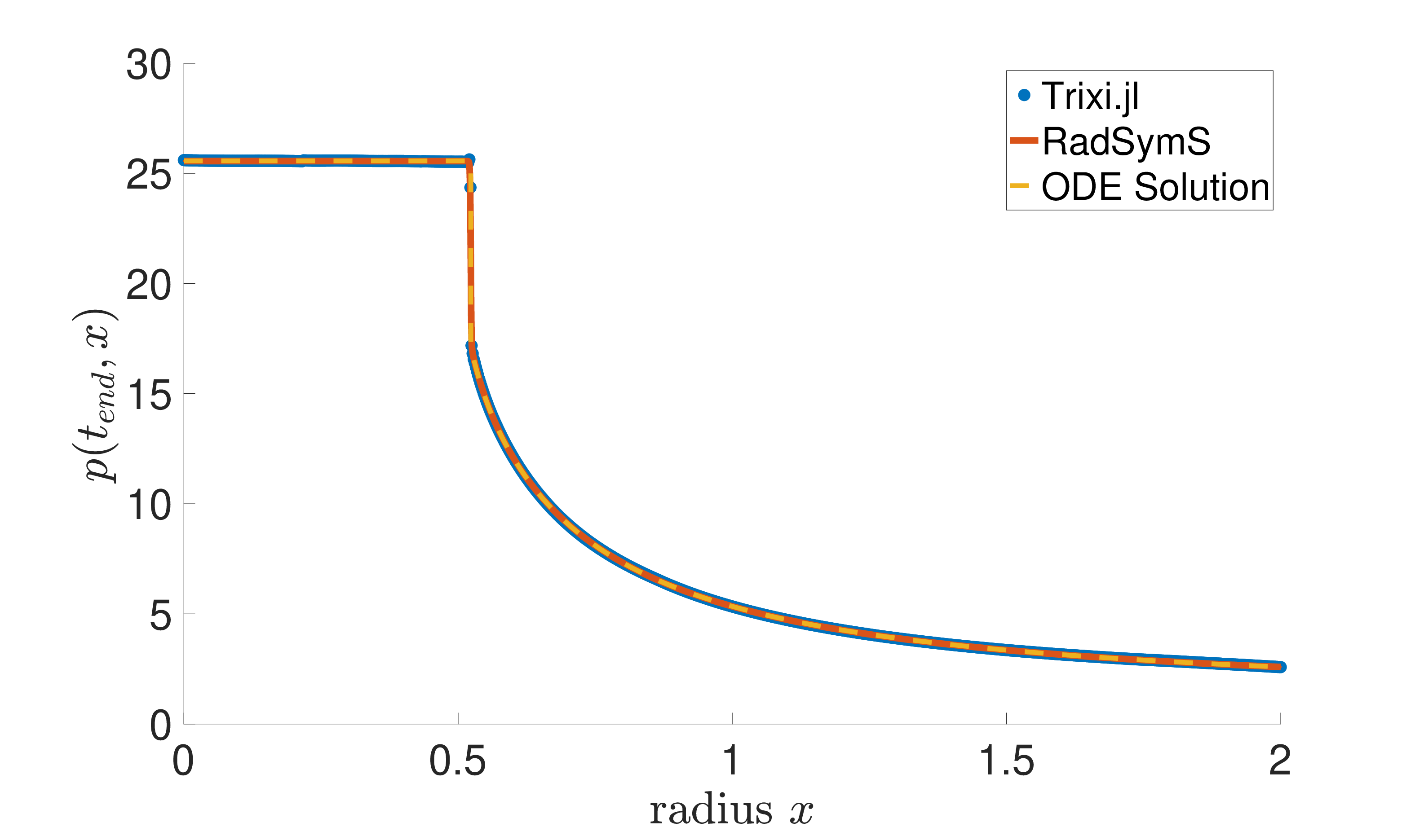}}
    \hspace{-1cm}
    \subfigure[Results for the velocity $v$ in 3D]{
    \includegraphics[width=0.575\textwidth]{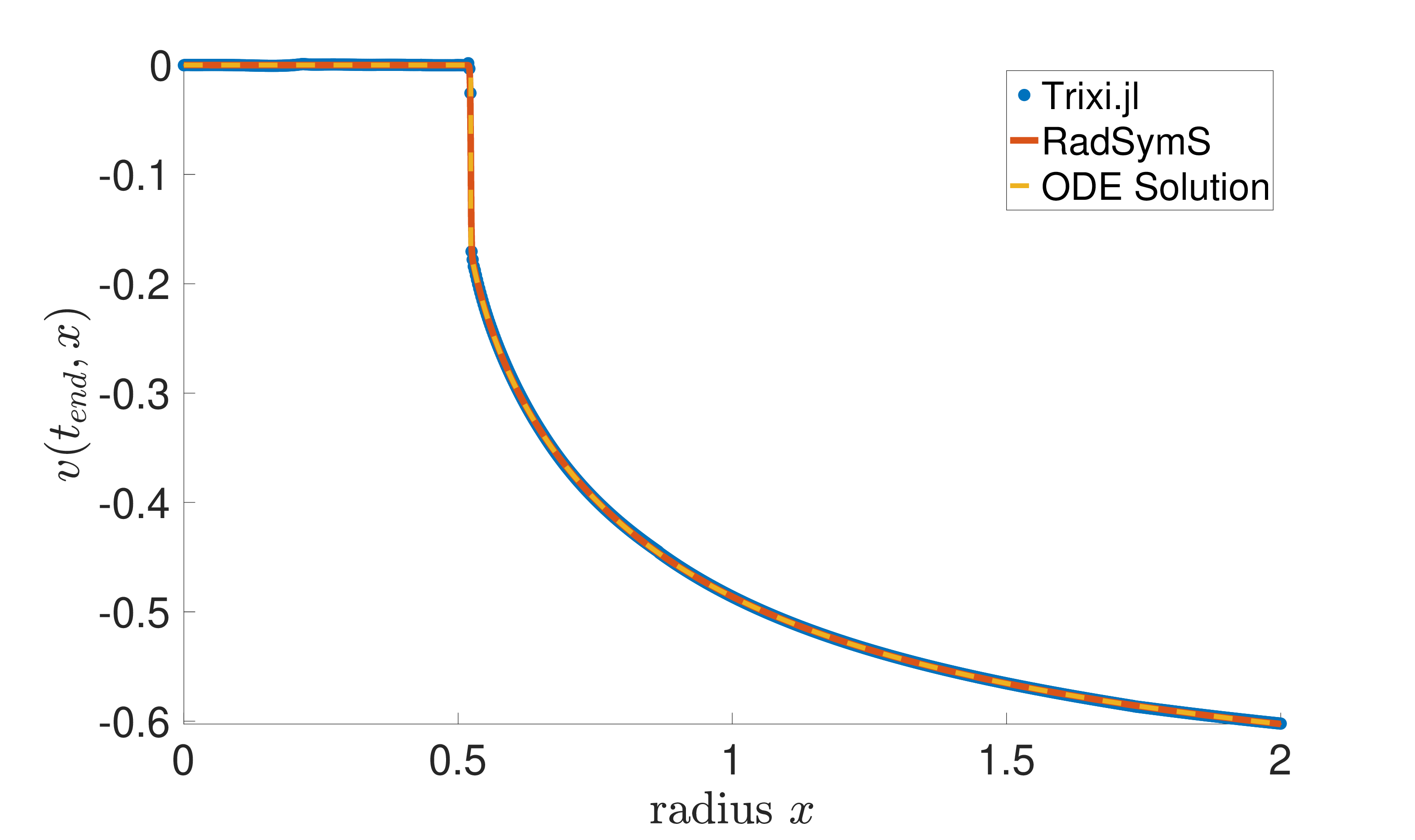}}
    \caption{   \textbf{Example 1:}
    Results for the pressure $p$ and velocity $v$ in 3D
     Comparison for Trixi.jl (blue) and RadSymS (red) in radial direction  with the solution of ODE system \eqref{GL_ODE_sys} (yellow) at final time $t_{end} = 1$ for the pressure $p$, see (a), and velocity $v$, see (b).
  }
    \label{fig:example-1-3D}
\end{figure}

\begin{figure}[h!]
    \hspace{-1cm}
    \subfigure{
    \includegraphics[width=0.575\textwidth]{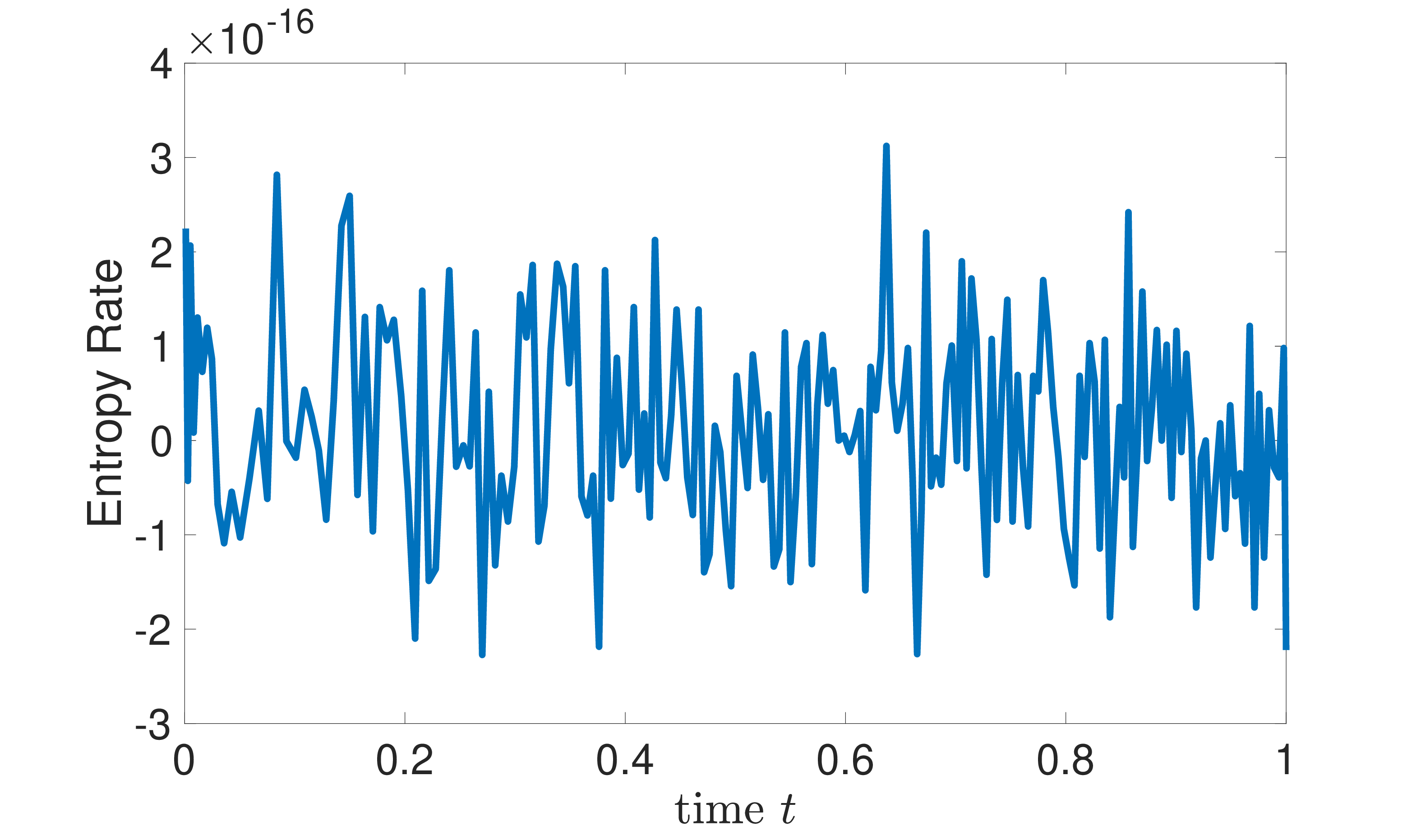}}
    \hspace{-1cm}
    \subfigure{
    \includegraphics[width=0.575\textwidth]{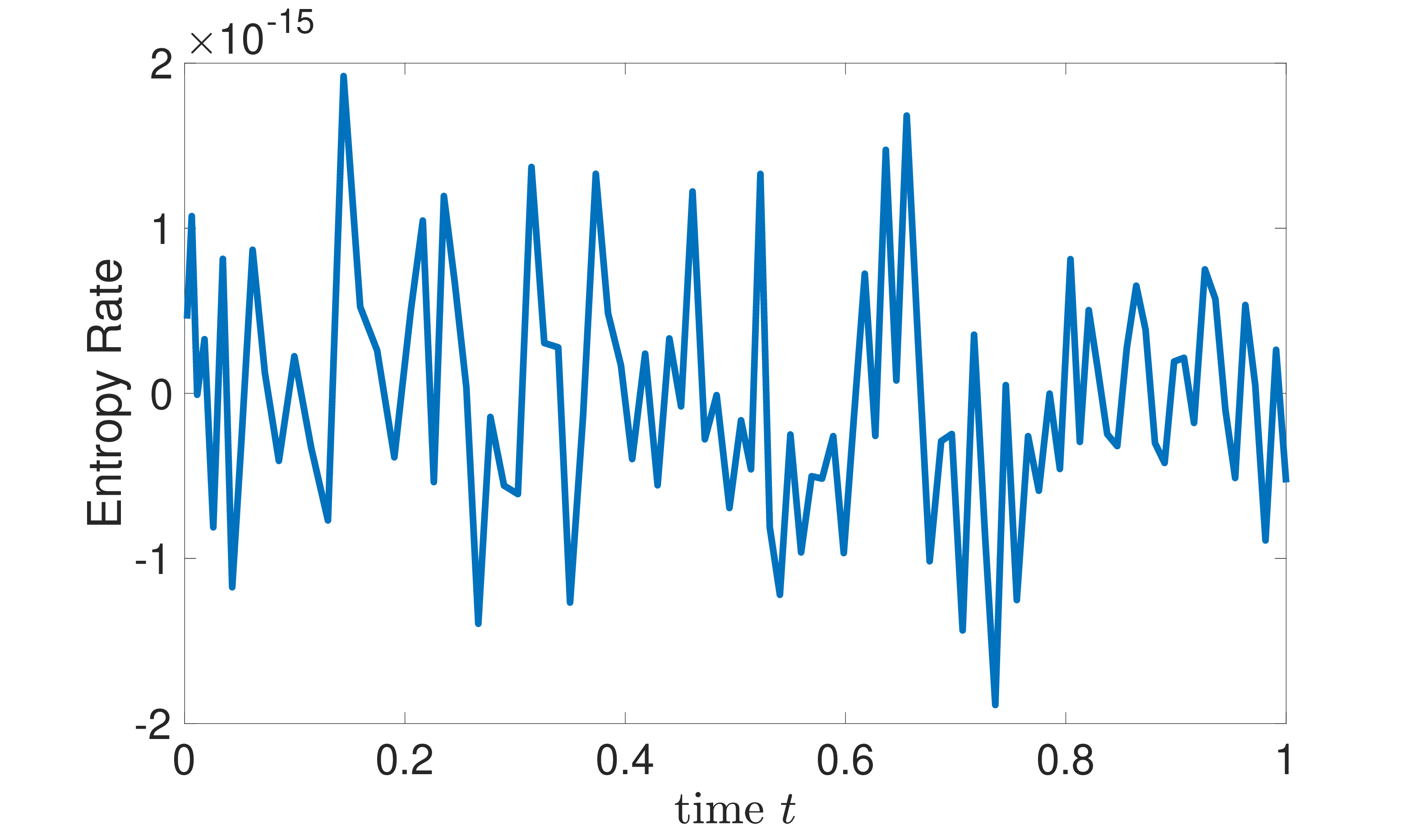}}
    \caption{Entropy rate in 2D (left) and 3D (right).}
    \label{fig:entropy_rate}
\end{figure}

\begin{figure}[h!]
    \hspace{-1cm}
    \subfigure[Results for the pressure $p$ in 2D]{
    \includegraphics[width=0.575\textwidth]{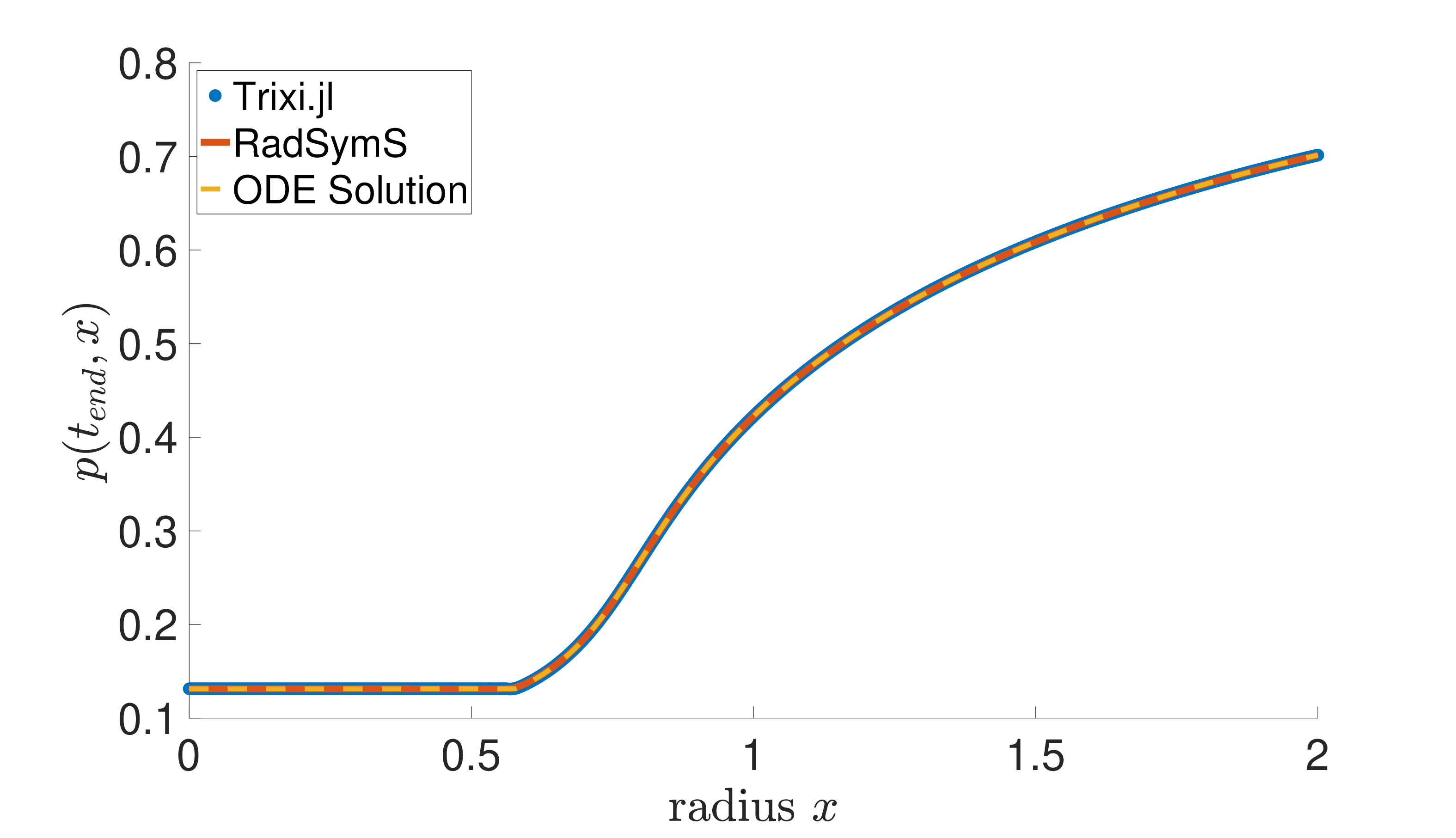}}
    \hspace{-1cm}
    \subfigure[Results for the velocity $v$ in 2D]{
    \includegraphics[width=0.575\textwidth]{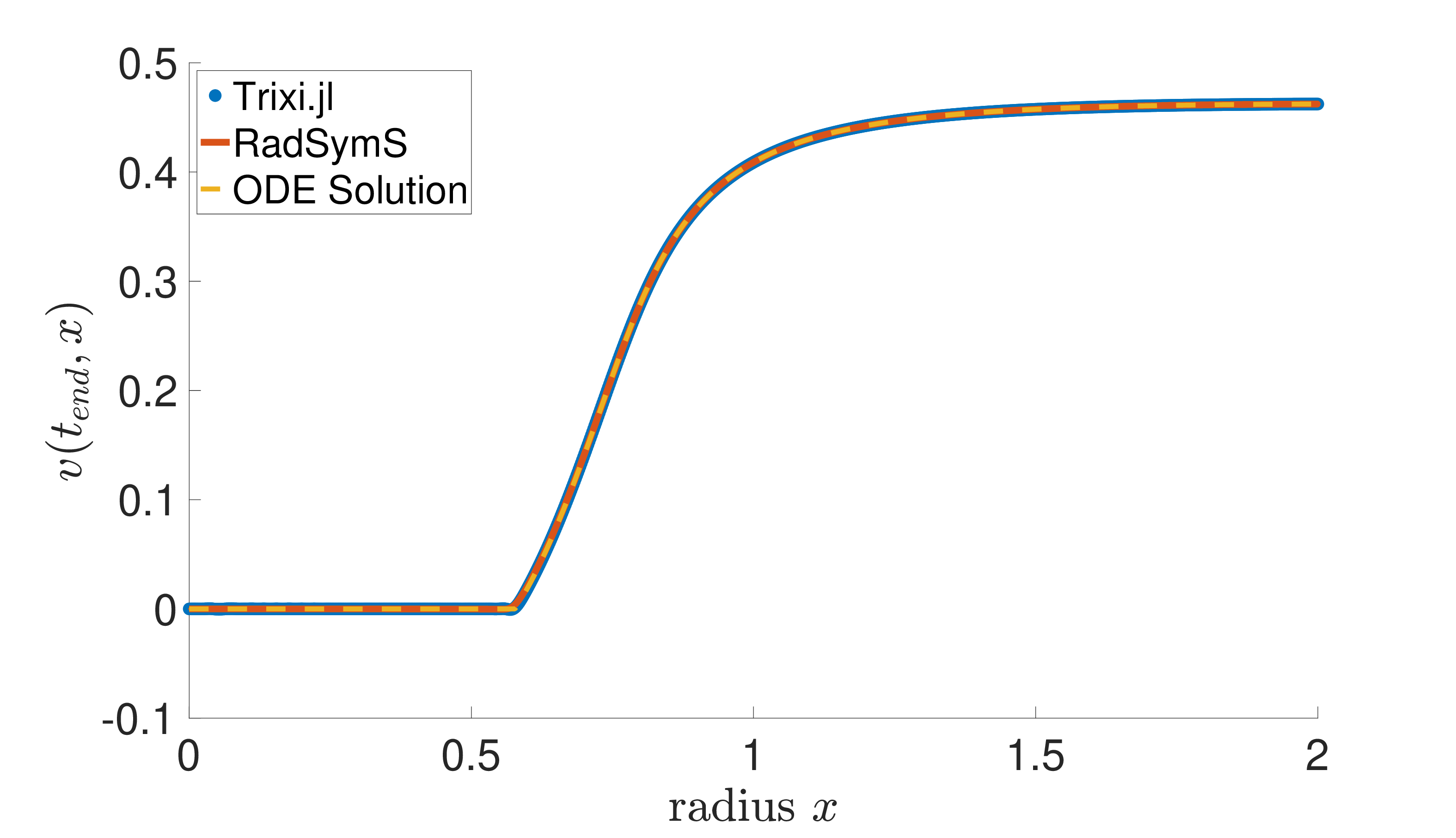}}
    \caption{\textbf{Example 2:}
    Results for the pressure $p$ and velocity $v$ in 2D.
     Comparison for Trixi.jl (blue) and RadSymS (red) in radial direction  with the solution of ODE system \eqref{GL_ODE_sys} (yellow) at final time $t_{end} = 1$ for the pressure $p$, see (a), and velocity $v$, see (b).
  }
    \label{fig:example-2-2D}
\end{figure}

\begin{figure}[h!]
    \hspace{-1cm}
    \subfigure[Results for the pressure $p$ in 3D]{
    \includegraphics[width=0.575\textwidth]{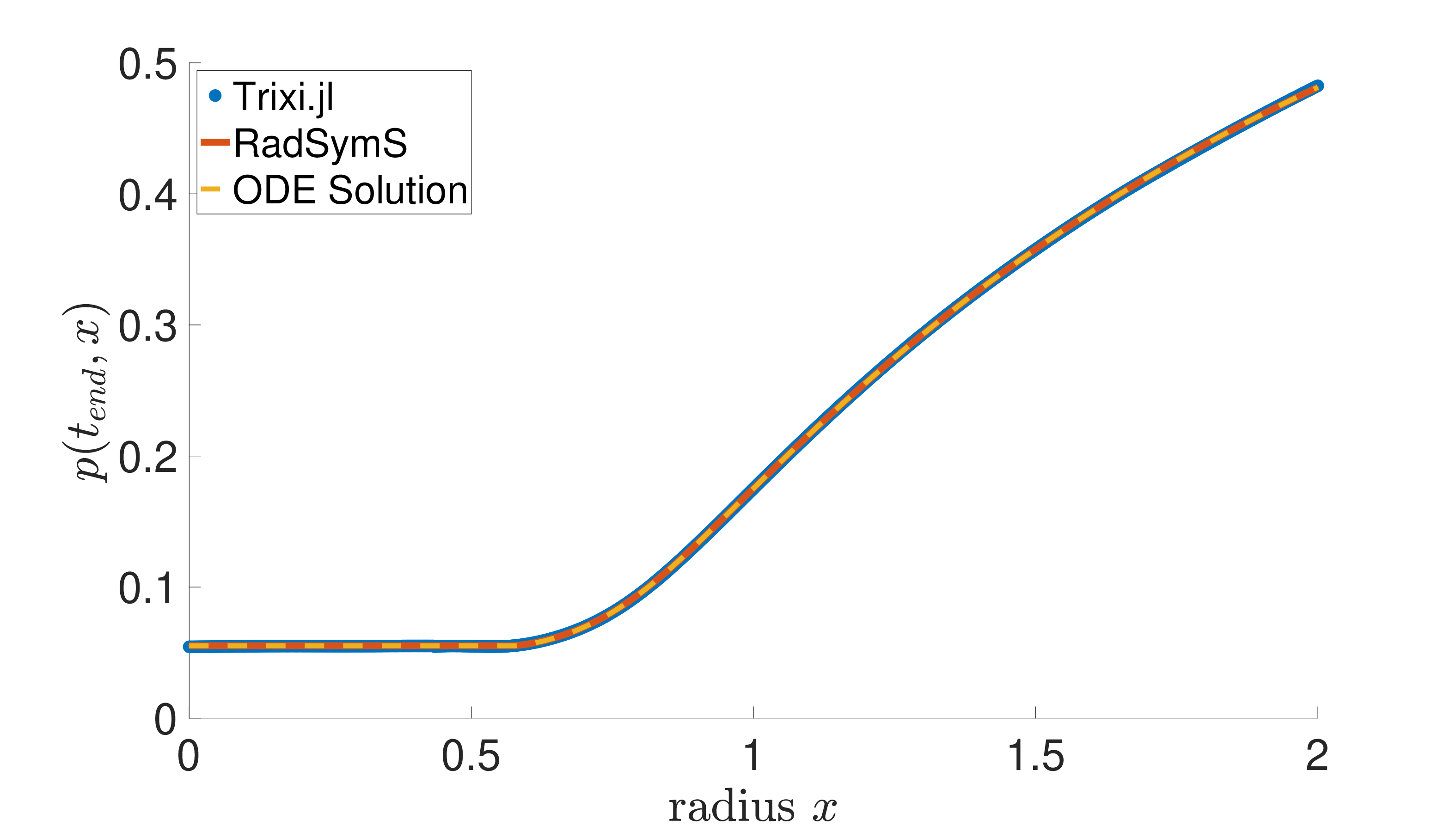}}
    \subfigure[Results for the velocity $v$ in 3D]{
    \includegraphics[width=0.575\textwidth]{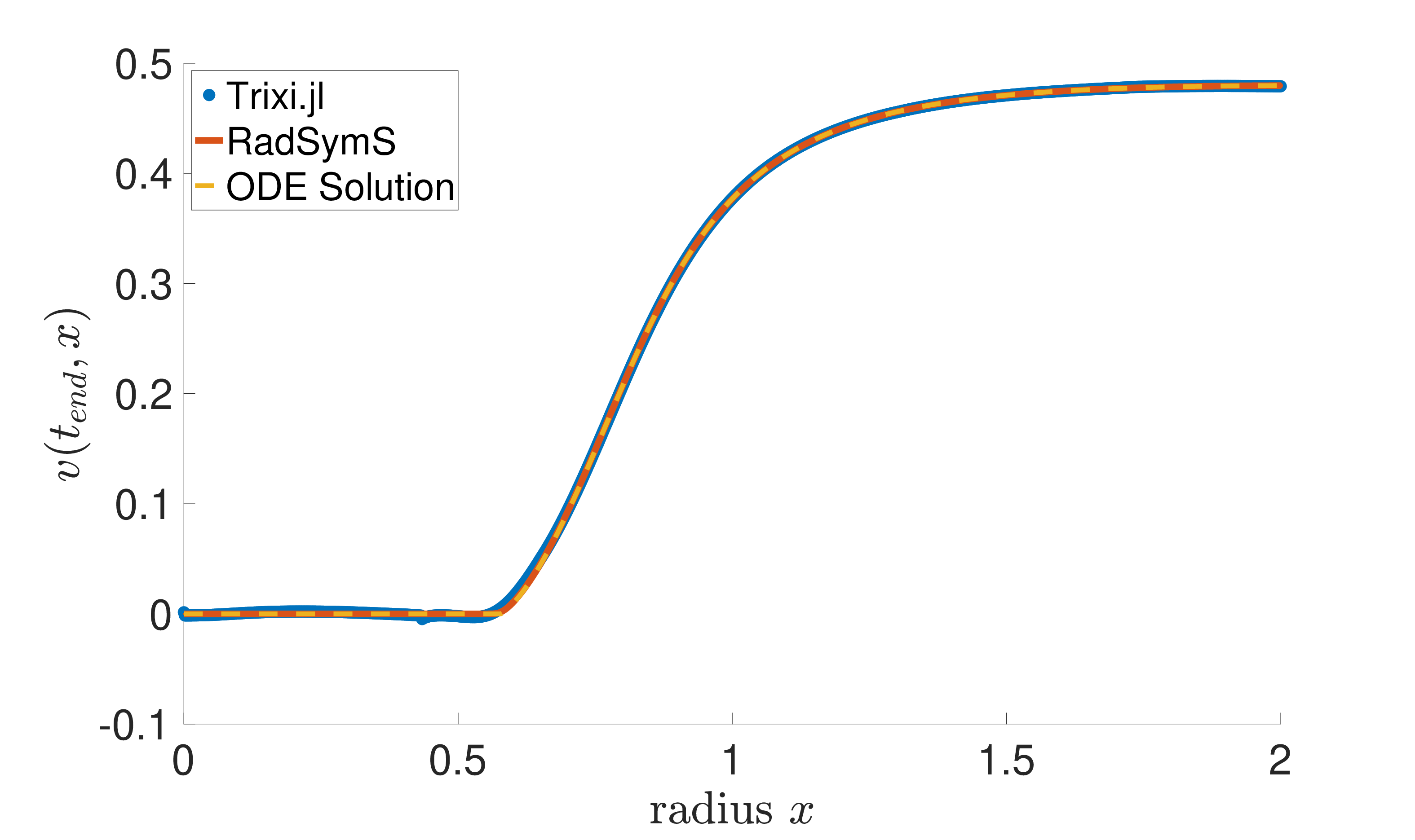}}
    \caption{\textbf{Example 2:}
    Results for the pressure $p$ and velocity $v$ in 3D.
     Comparison for Trixi.jl (blue) and RadSymS (red) in radial direction  with the solution of ODE system \eqref{GL_ODE_sys} (yellow) at final time $t_{end} = 1$ for the pressure $p$, see (a), and velocity $v$, see (b).
  }
    \label{fig:example-2-3D}
\end{figure}

\begin{figure}[htb!]
	%
    %
    \subfigure{
    \includegraphics[width=0.55\textwidth]{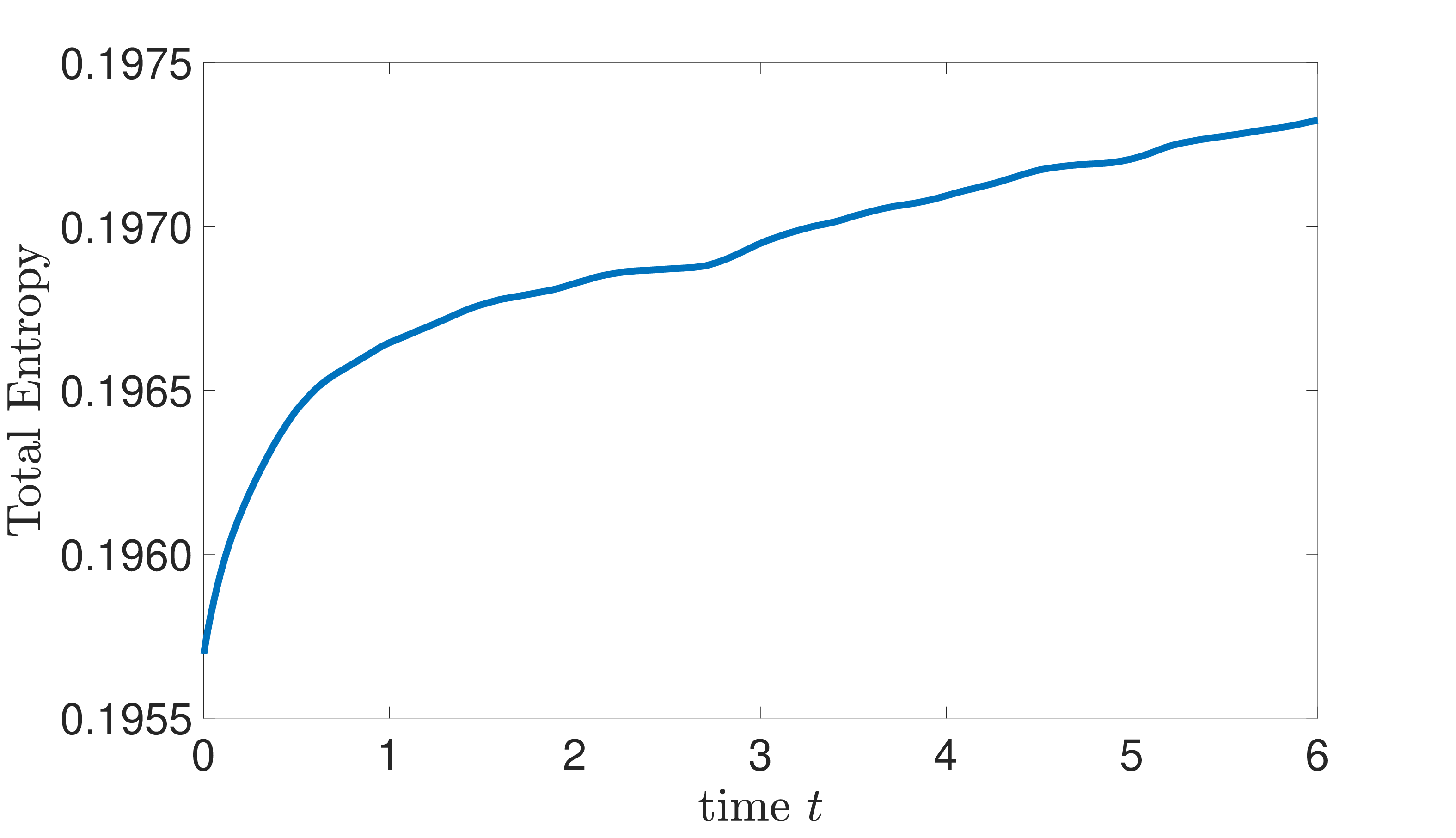}}
    \hspace{-1cm}
    \subfigure{
    \includegraphics[width=0.55\textwidth]{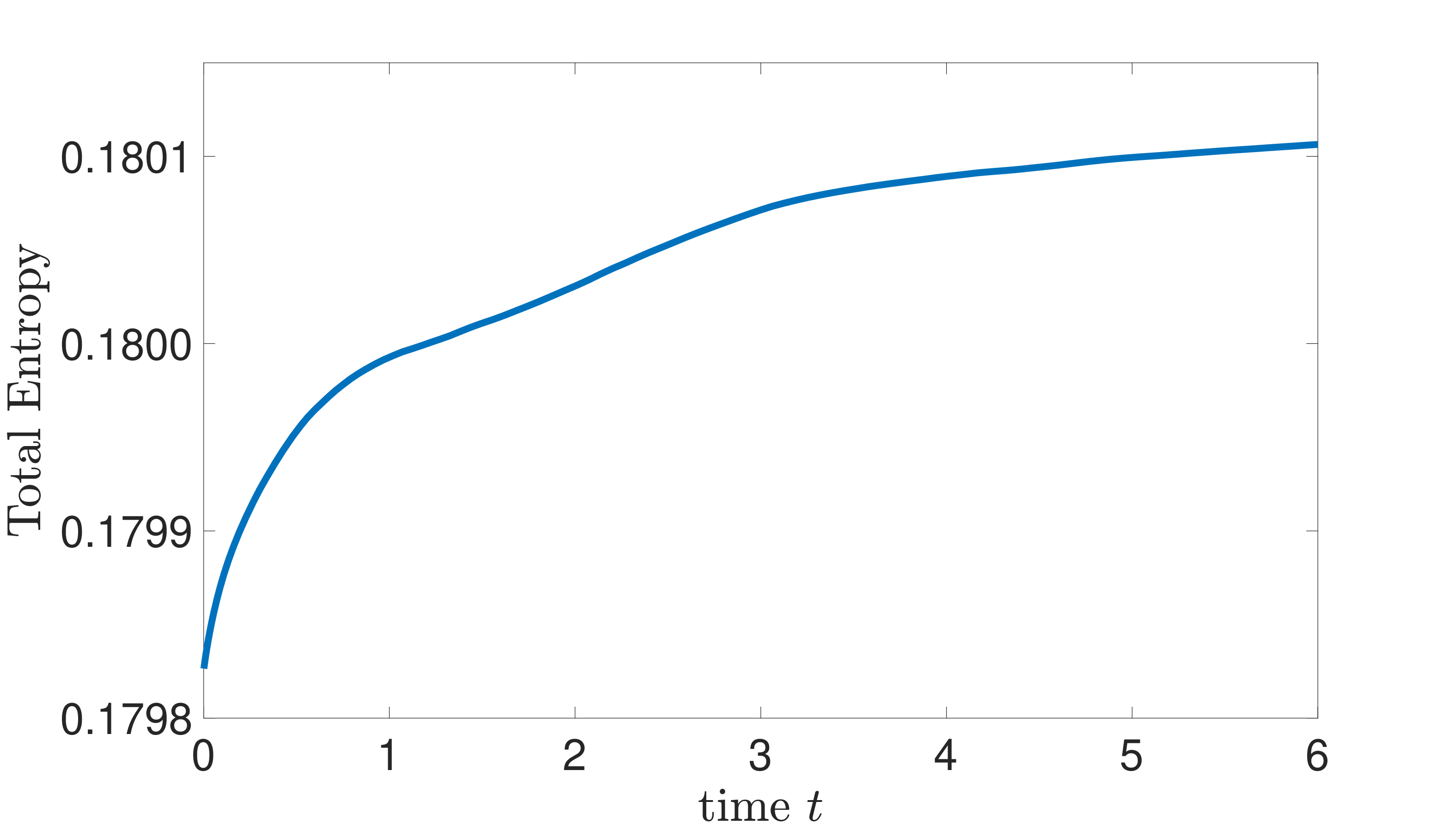}}
    \caption{Total entropy in 2D (left) and 3D (right). Since the (physical) entropy $\eta$ is concave, entropy stability means that the total entropy increases over time.}
    \label{fig:entropy}
\end{figure}

\begin{figure}[h!]
    \subfigure[Results for the pressure $p$ in 2D]{
    \includegraphics[width=\textwidth]{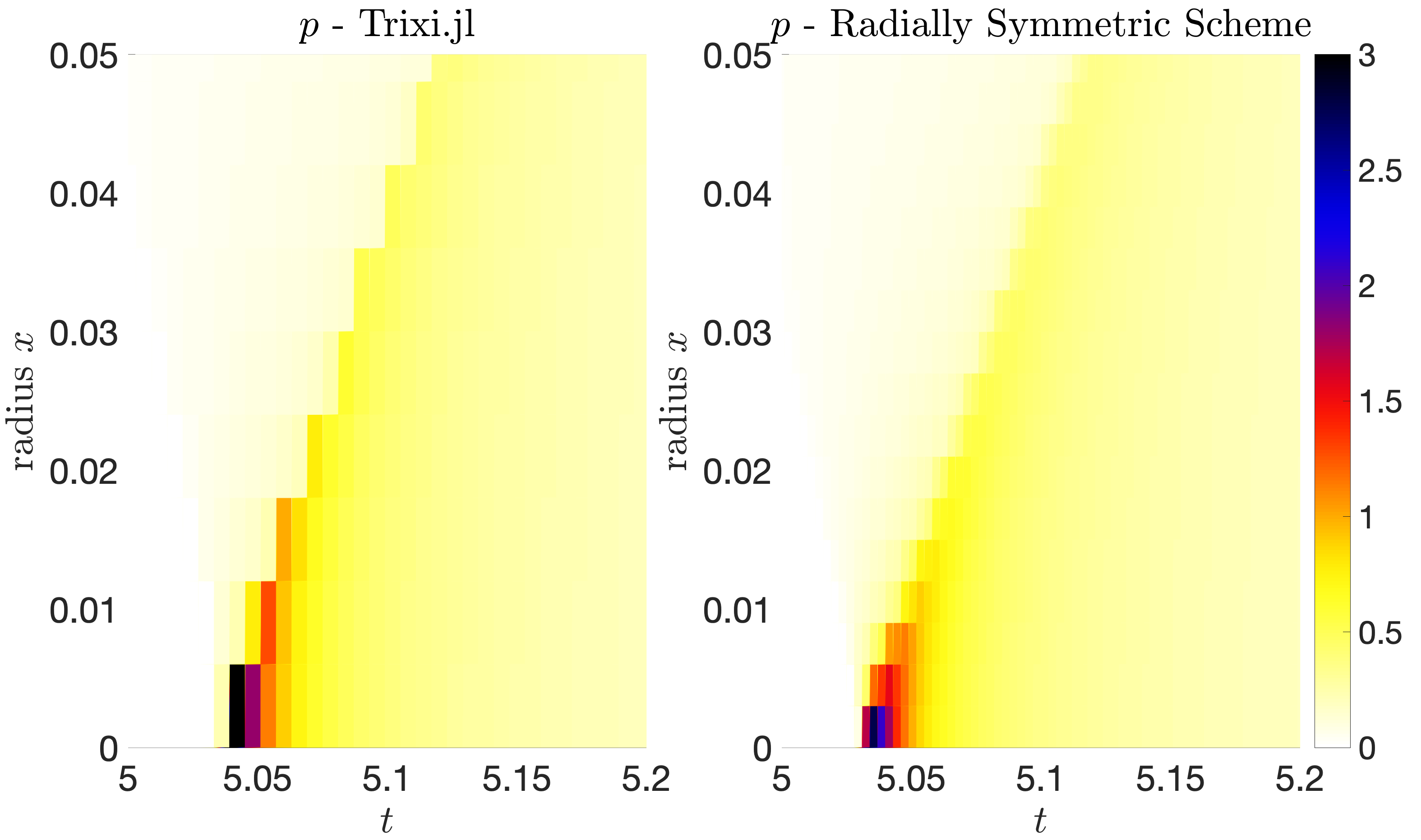}
    }\\
    \subfigure[Results for the velocity $v$ in 2D]{
    \includegraphics[width=\textwidth]{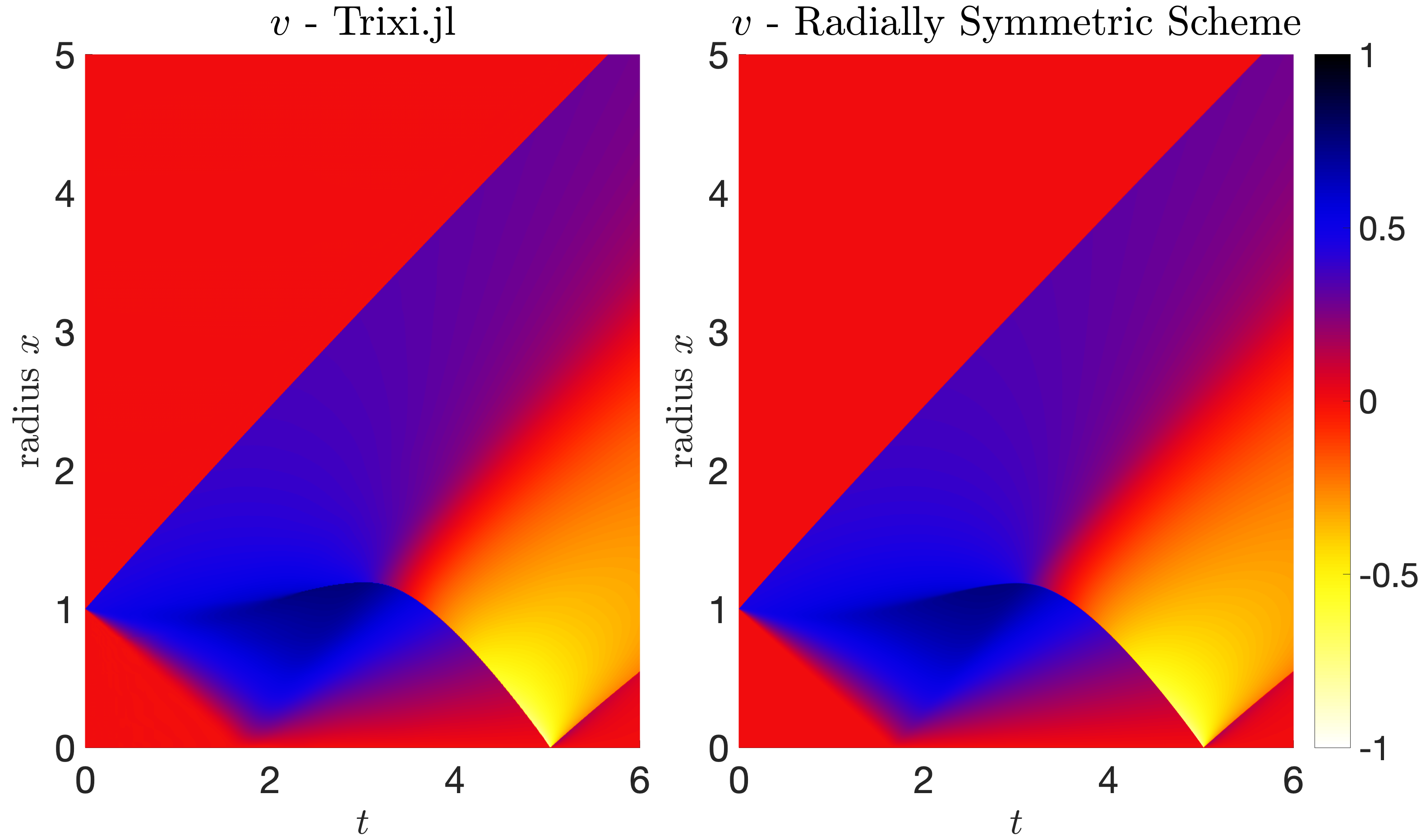}}
    \caption{\textbf{Example 3:}
    Comparison for Trixi.jl (left) and RadSymS (right) in the $t$--$x$ plane, see (a), (b). For the pressure we present a zoom plot to visualize the small vicinity of the origin, where the pressure increases localized in time.}
    \label{fig:example-3-2D-tr}
\end{figure}

\begin{figure}[h!]
    \subfigure[Results for the pressure $p$ in 3D]{
    \includegraphics[width=\textwidth]{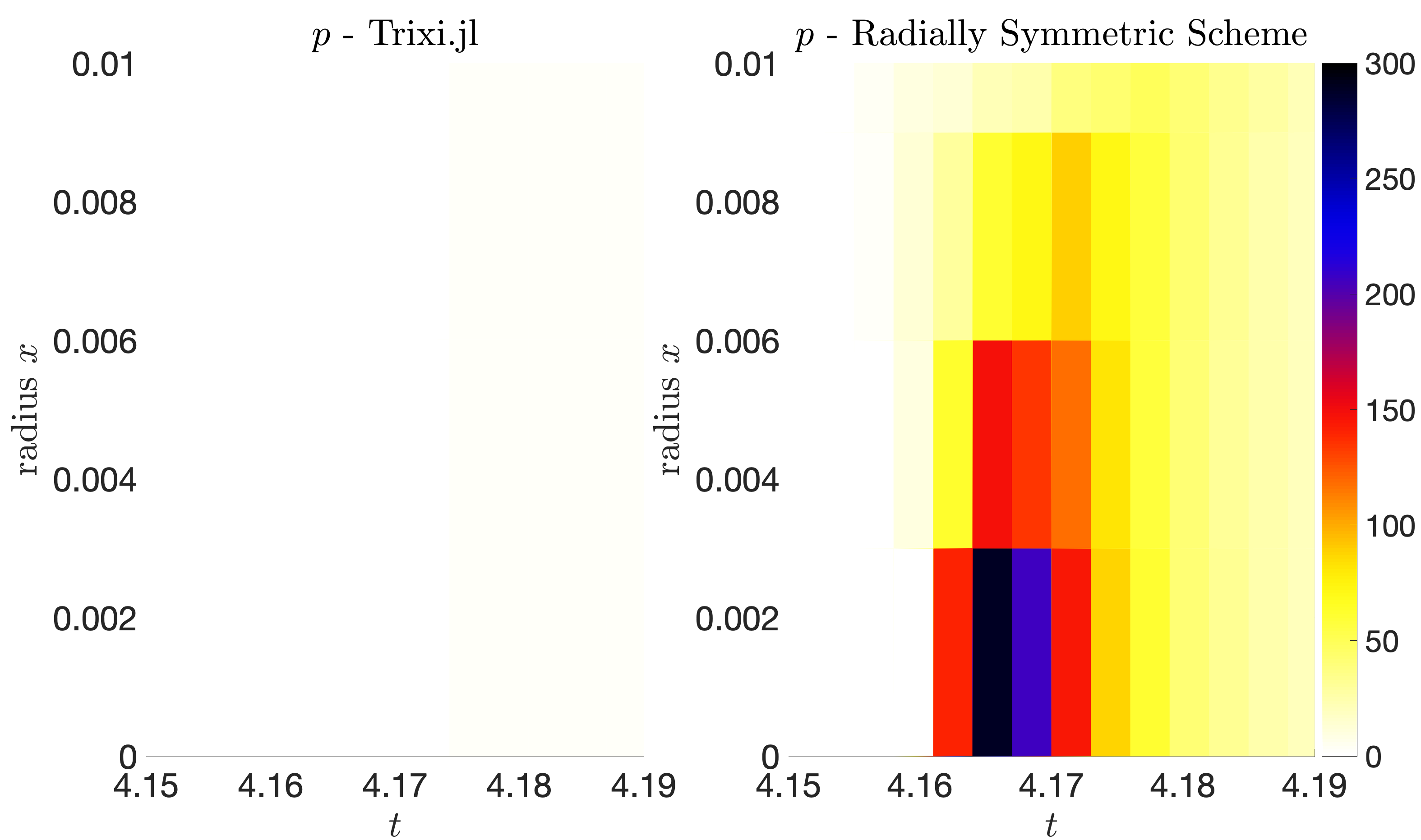}
    }\\
    \subfigure[Results for the velocity $v$ in 3D]{
    \includegraphics[width=\textwidth]{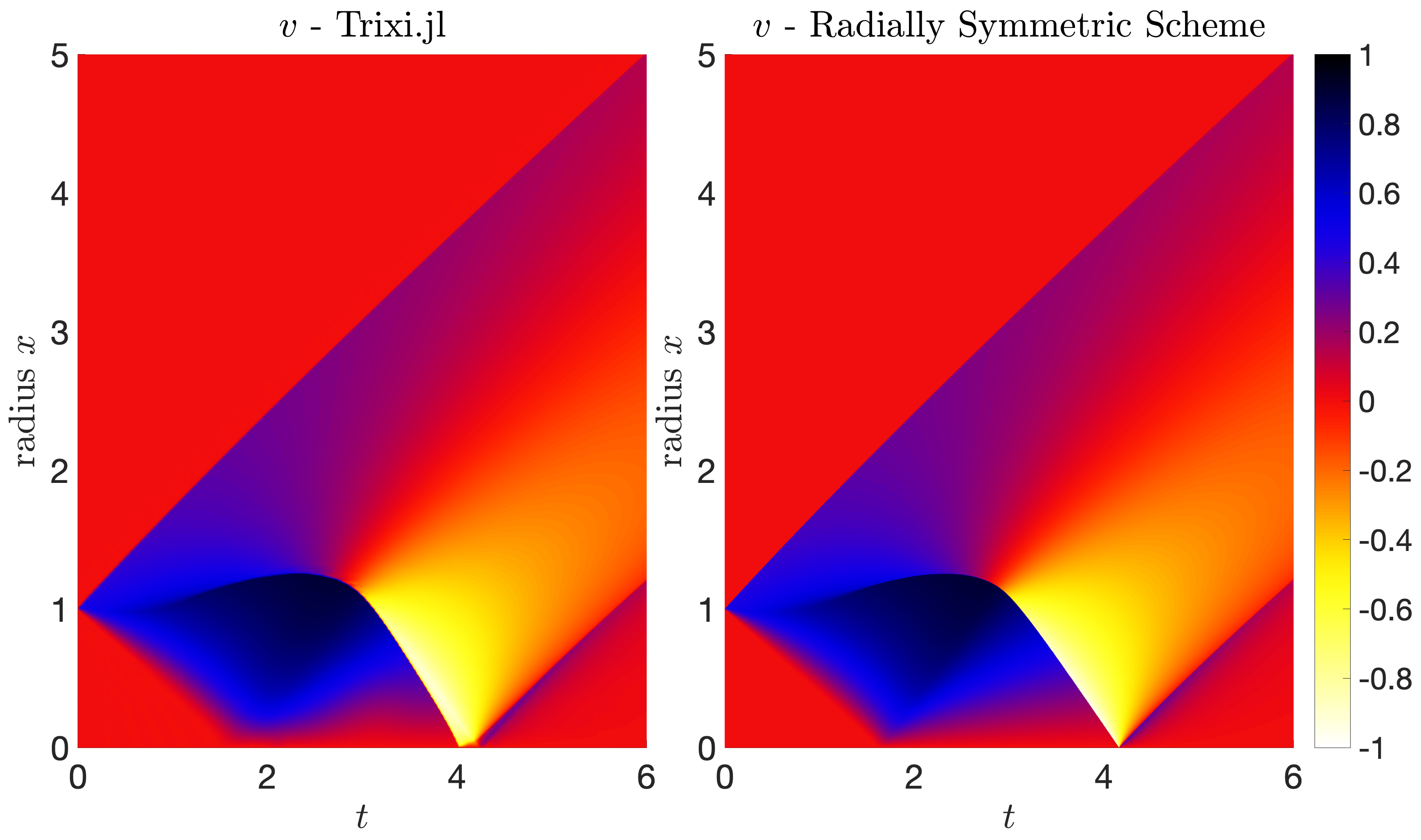}}
    \caption{\textbf{Example 3:}
    Comparison for Trixi.jl (left) and RadSymS (right) in the $t$--$x$ plane, see (a), (b). For the pressure we present a zoom plot to visualize the small vicinity of the origin, where the pressure increases localized in time.}
    \label{fig:example-3-3D-tr}
\end{figure}

\begin{figure}[h!]
    \hspace{-1cm}
    \subfigure[Results for the pressure $p$ in 2D]{
    \includegraphics[width=0.575\textwidth]{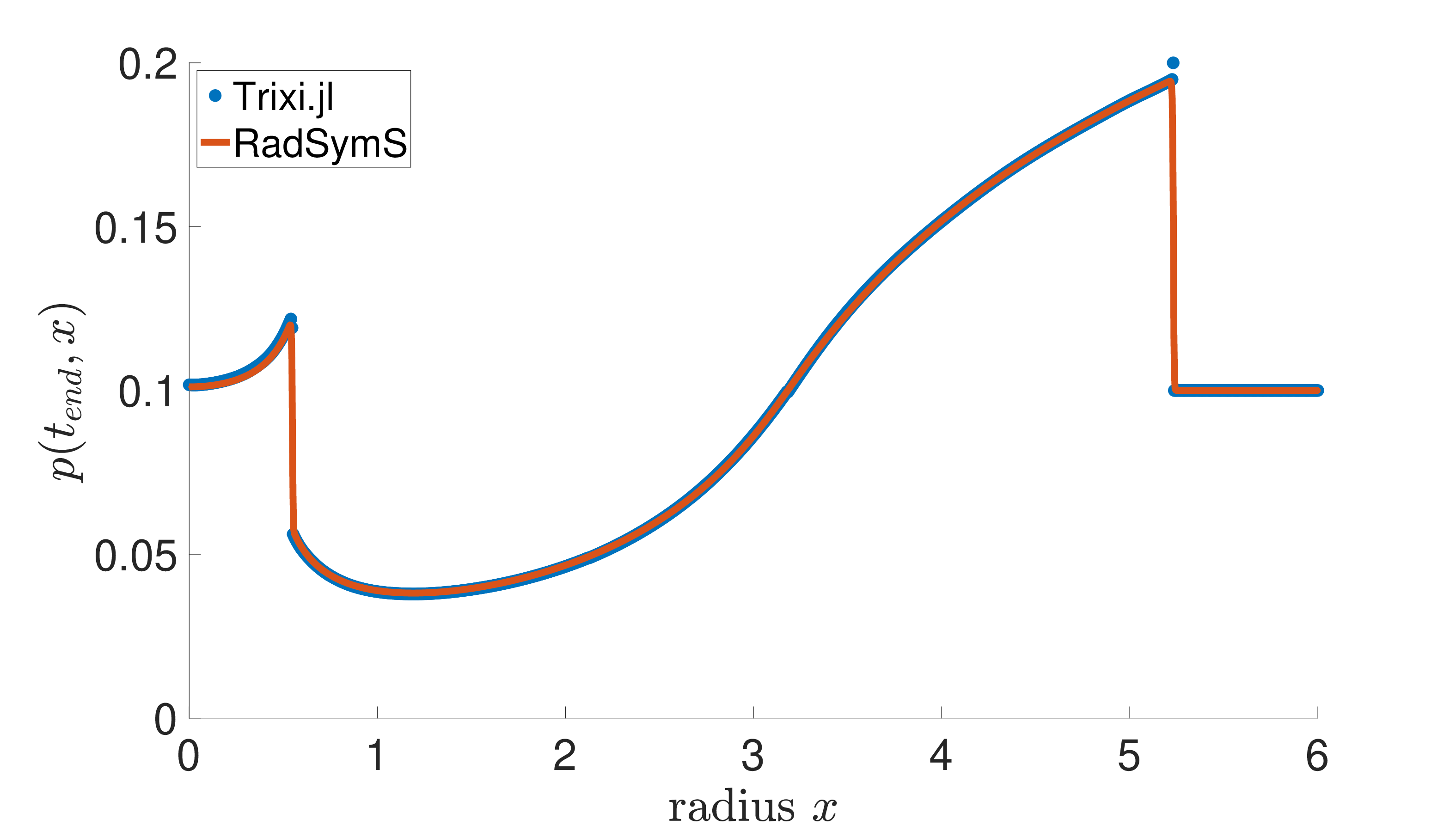}}
    \hspace{-1cm}
    \subfigure[Results for the velocity $v$ in 2D]{
    \includegraphics[width=0.575\textwidth]{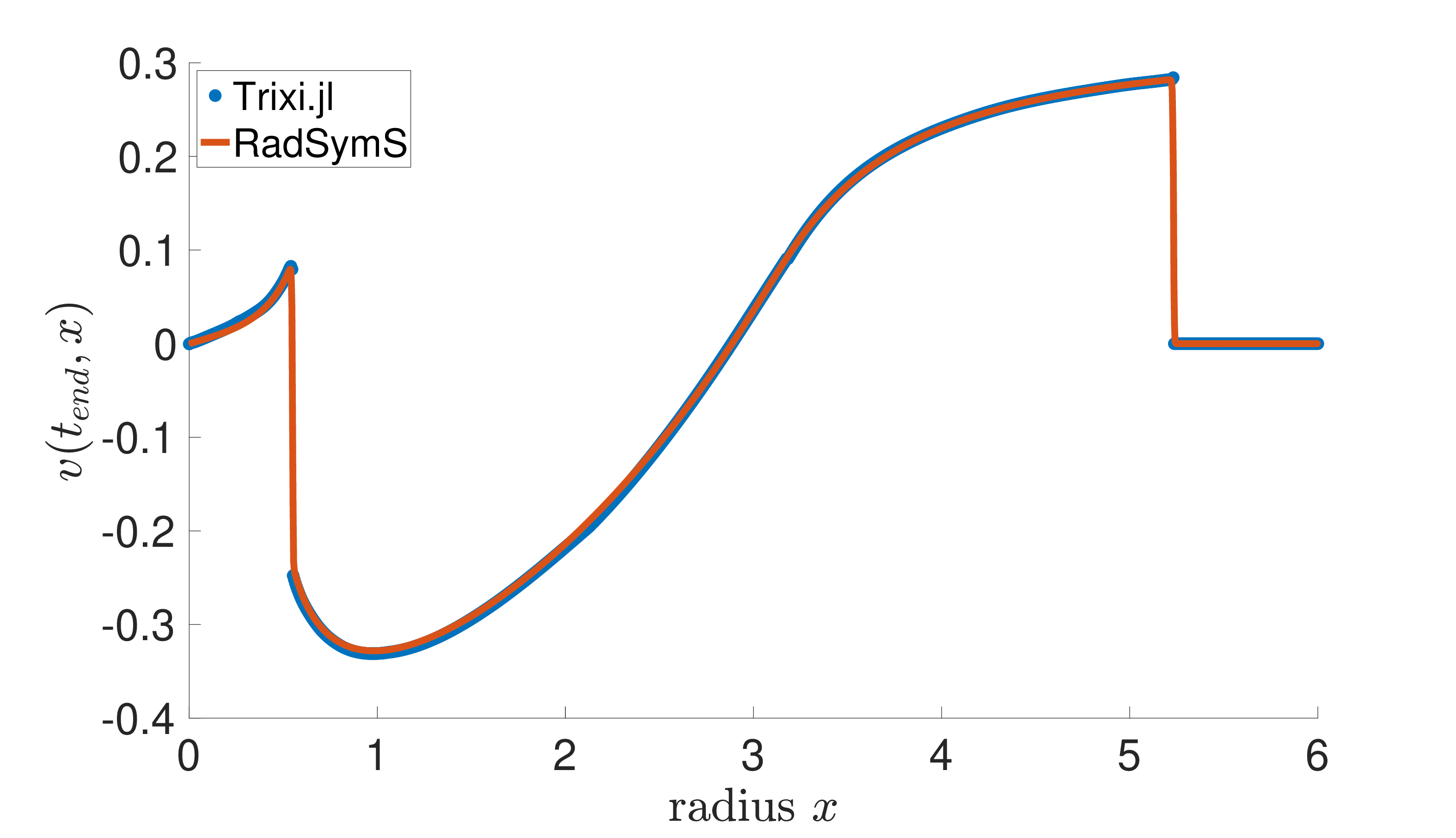}}
    \caption{\textbf{Example 3:}
    Comparison for RadSymS (red) and Trixi.jl (blue) in radial direction at final time $t_{end} = 6$, see (a), (b), for pressure $p$ and velocity $v$.}
    \label{fig:example-3-2D-tend}
\end{figure}

\begin{figure}[h!]
    \hspace{-1cm}
    \subfigure[Results for the pressure $p$ in 3D]{
    \includegraphics[width=0.575\textwidth]{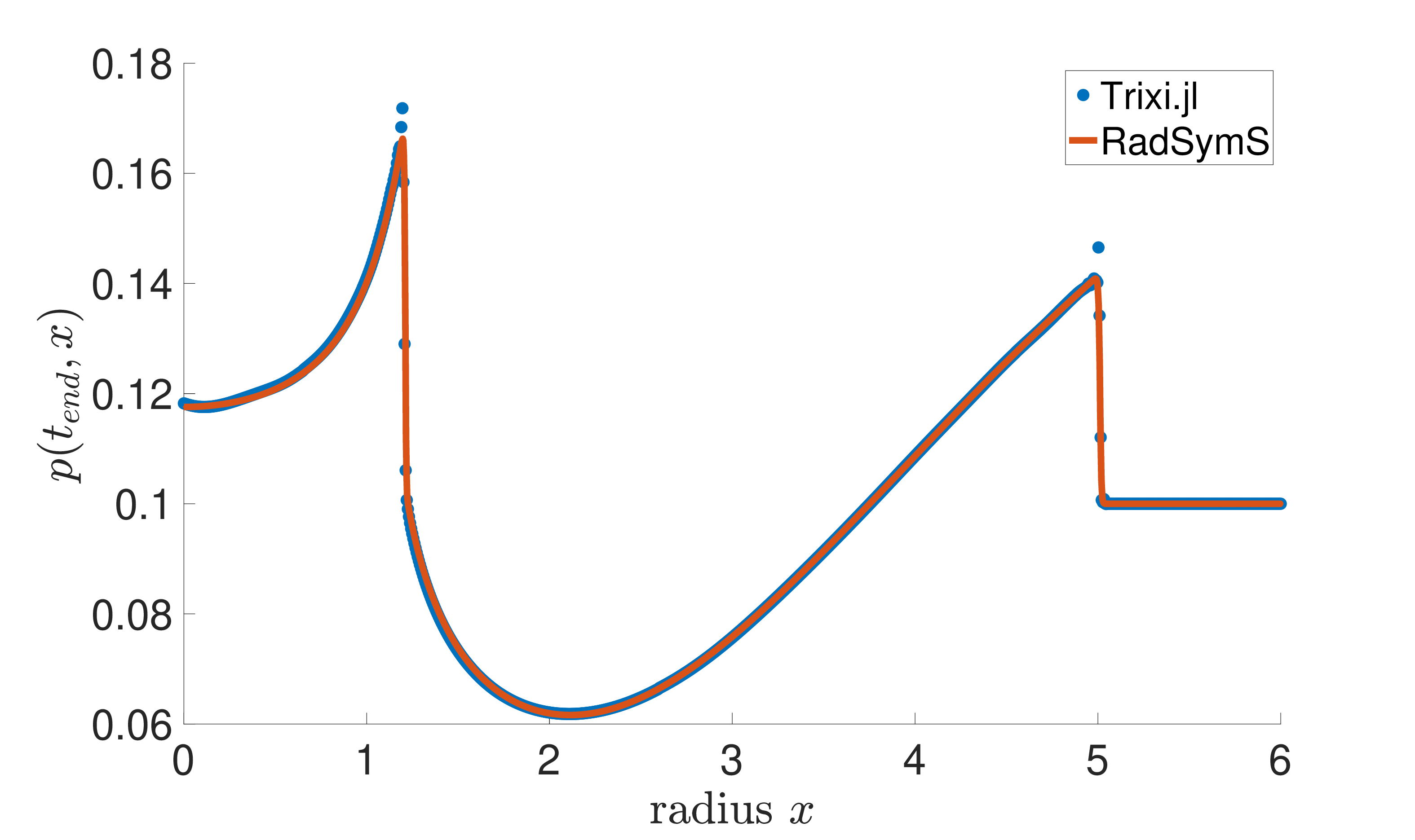}}
    \hspace{-1cm}
    \subfigure[Results for the velocity $v$ in 3D]{
    \includegraphics[width=0.575\textwidth]{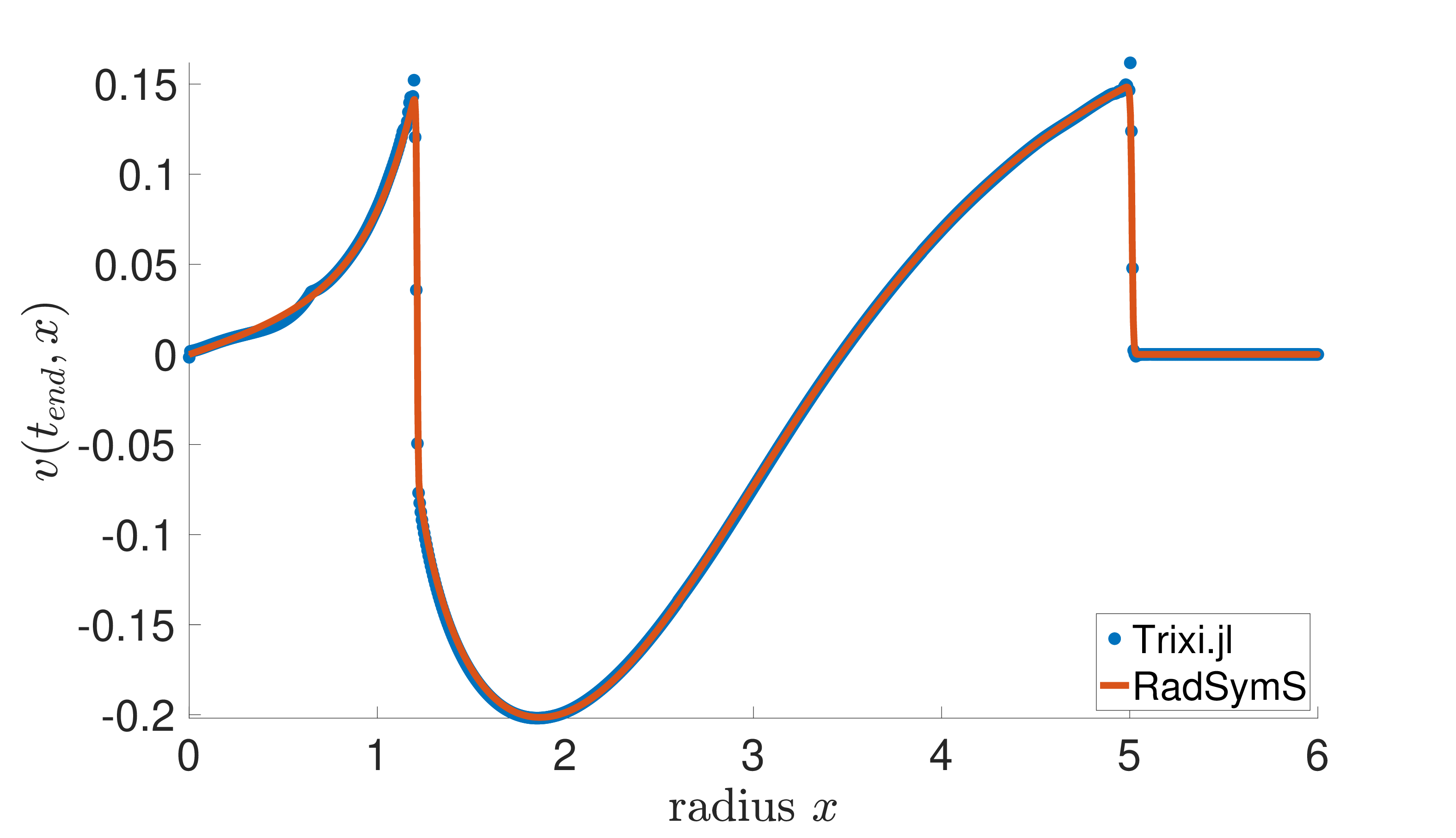}}
    \caption{\textbf{Example 3:}
    Comparison for RadSymS (red) and Trixi.jl (blue) in radial direction at final time $t_{end} = 6$, see (a), (b), for pressure $p$ and velocity $v$.}
    \label{fig:example-3-3D-tend}
\end{figure}

\begin{figure}[h!]
    \hspace{-1cm}
    \subfigure[Results for the pressure $p$ in 2D]{
    \includegraphics[width=0.575\textwidth]{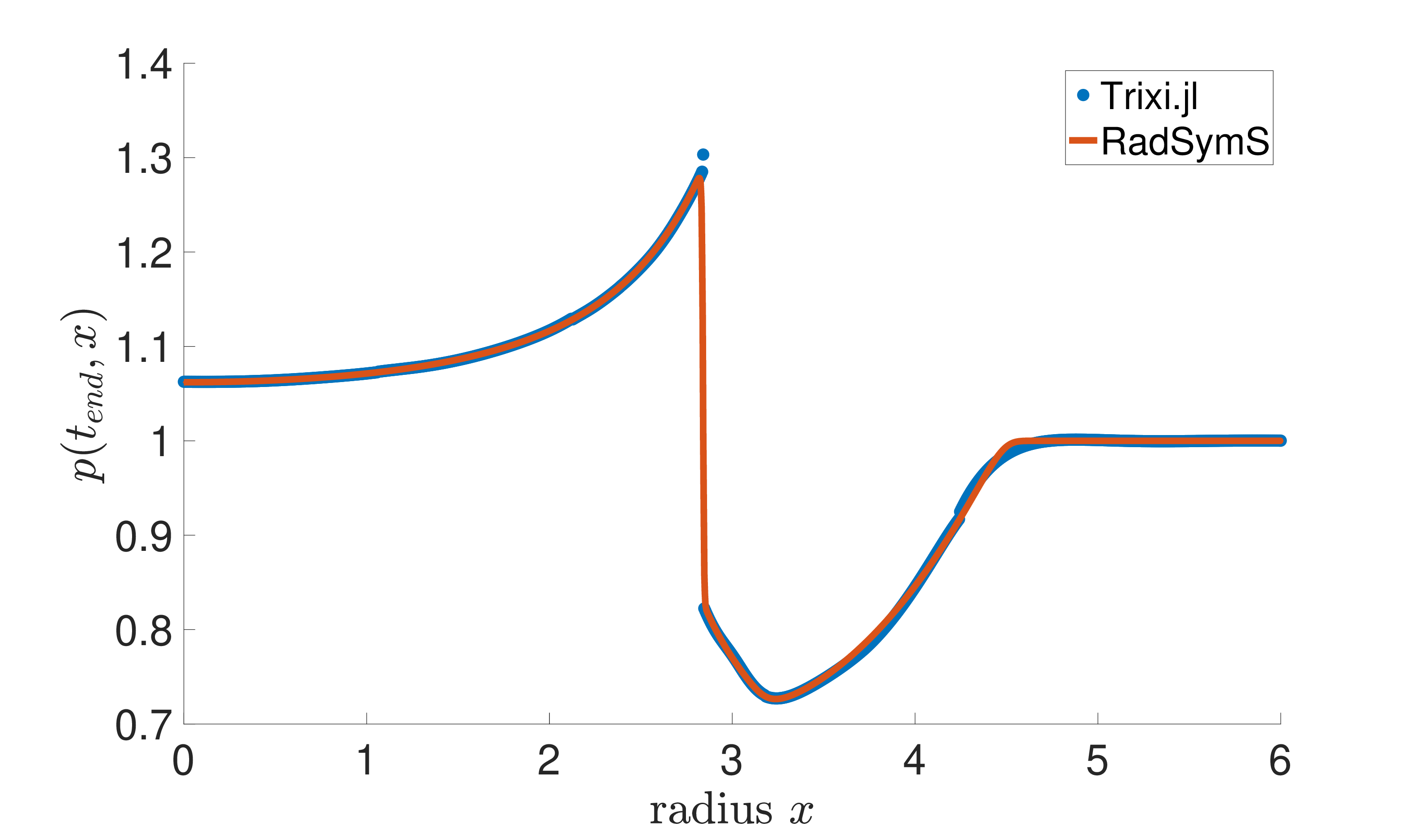}}
    \hspace{-1cm}
    \subfigure[Results for the velocity $v$ in 2D]{
    \includegraphics[width=0.575\textwidth]{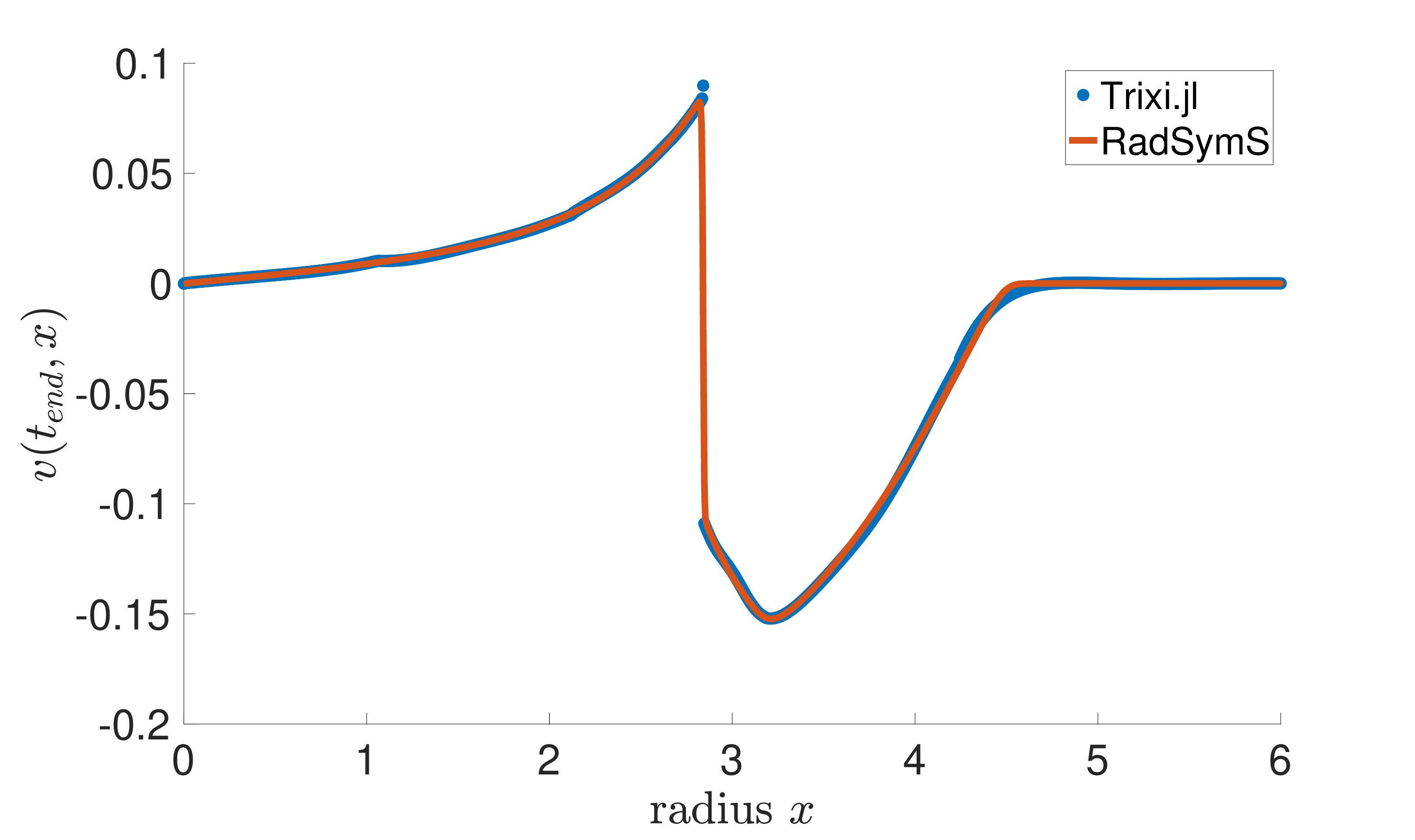}}
    \caption{\textbf{Example 4:}
    Results for the velocity $v$ and the pressure $p$ in 2D.
    Comparison for RadSymS (red) and Trixi.jl (blue) in radial direction at final time $t_{end} = 6$, see (a), (b), for pressure $p$ and velocity $v$.}
     \label{fig:example-4-2D-tend}
\end{figure}

\begin{figure}[h!]
    \subfigure[Results for the pressure $p$ in 2D]{
    \includegraphics[width=\textwidth]{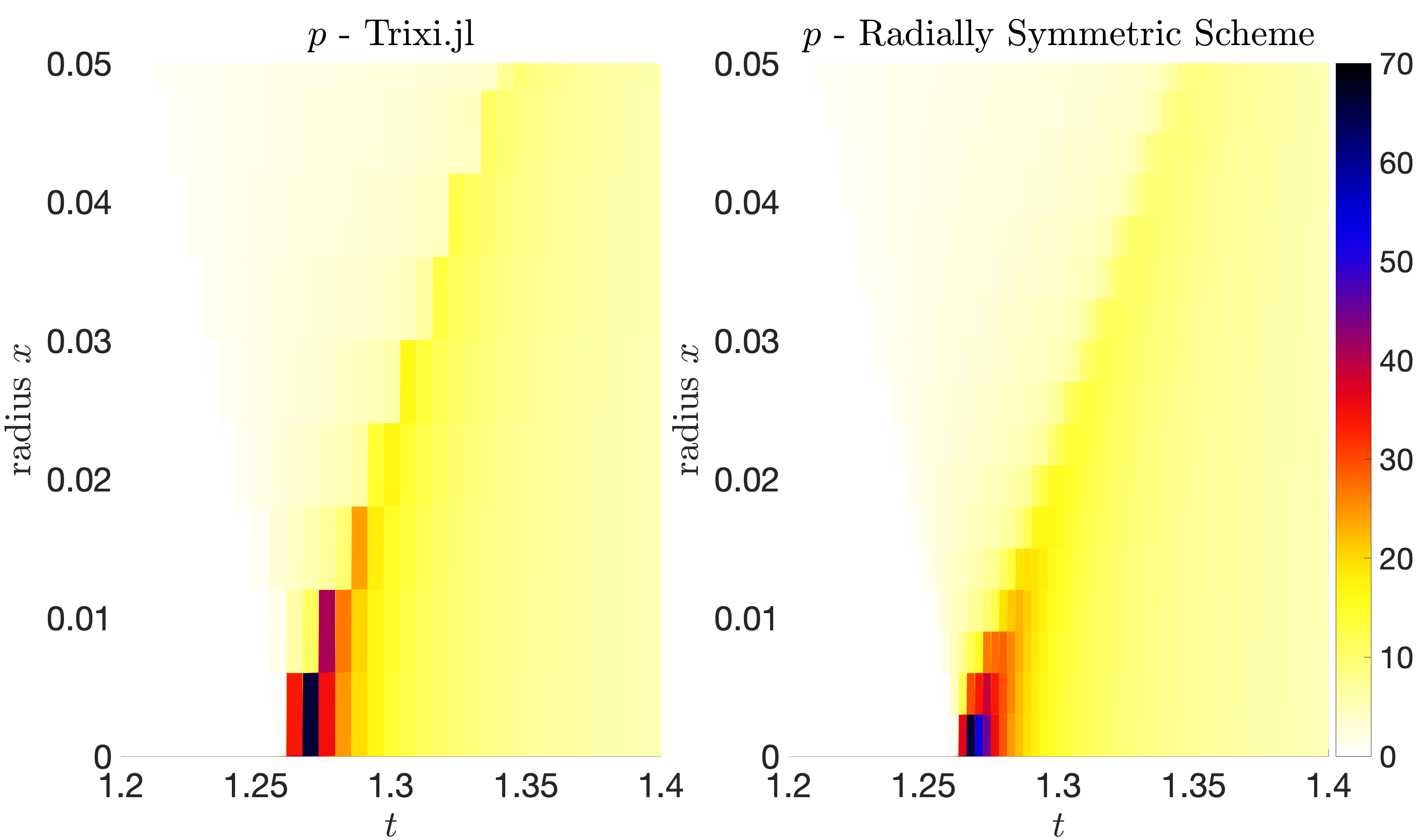}
    }\\
    \subfigure[Results for the velocity $v$ in 2D]{
    \includegraphics[width=\textwidth]{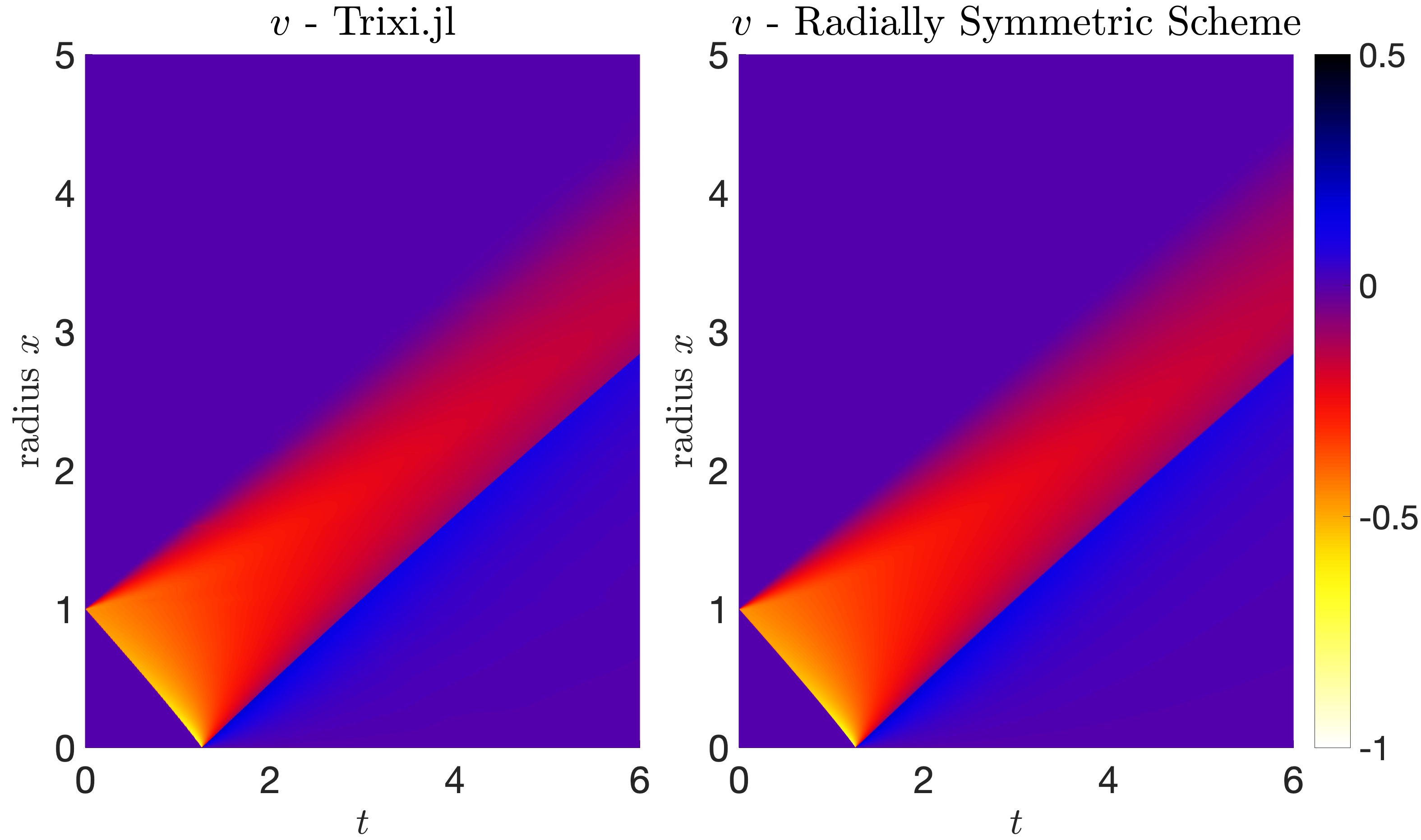}}
    \caption{\textbf{Example 4:}
    Results for the velocity $v$ and the pressure $p$ in 2D.
     Comparison for Trixi.jl (left) and RadSymS (right) in the $t$--$x$ plane, see (a), (b). For the pressure we present a zoom plot to visualize the small vicinity of the origin, where the pressure increases localized in time.}
     \label{fig:example-4-2D-tr}
\end{figure}

\begin{figure}[h!]
   	\subfigure[Results for the pressure $p$ in 3D]{
    \includegraphics[width=\textwidth]{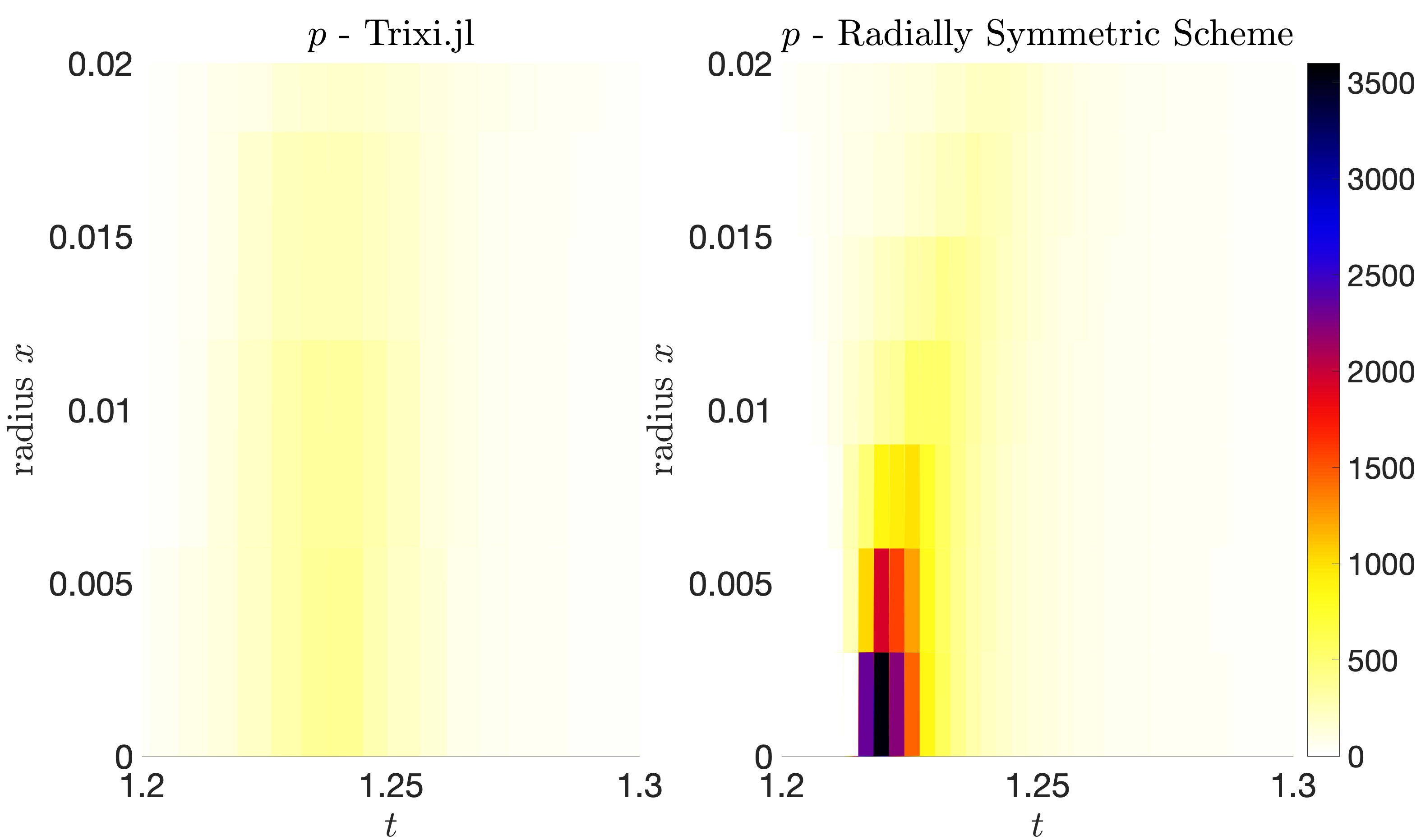}
    }\\
    \subfigure[Results for the velocity $v$ in 3D]{
    \includegraphics[width=\textwidth]{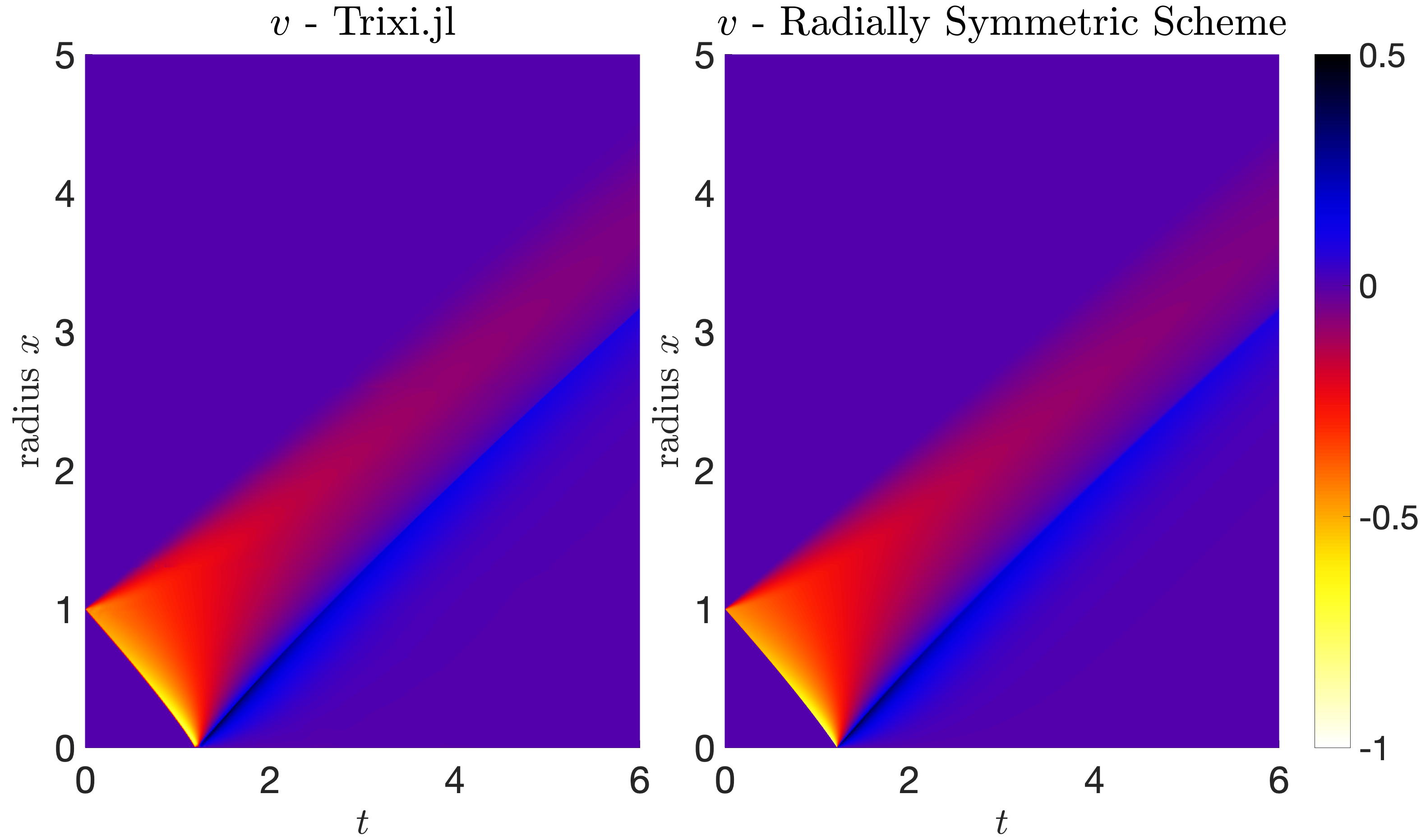}}
    \caption{\textbf{Example 4:}
    Results for the velocity $v$ and the pressure $p$ in 3D.
     Comparison for Trixi.jl (left) and RadSymS (right) in the $t$--$x$ plane, see (a), (b). For the pressure we present a zoom plot to visualize the small vicinity of the origin, where the pressure increases localized in time.}
     \label{fig:example-4-3D-tr}
\end{figure}

\begin{figure}[h!]
    \hspace{-1cm}
    \subfigure[Results for the pressure $p$ 3D]{
    \includegraphics[width=0.575\textwidth]{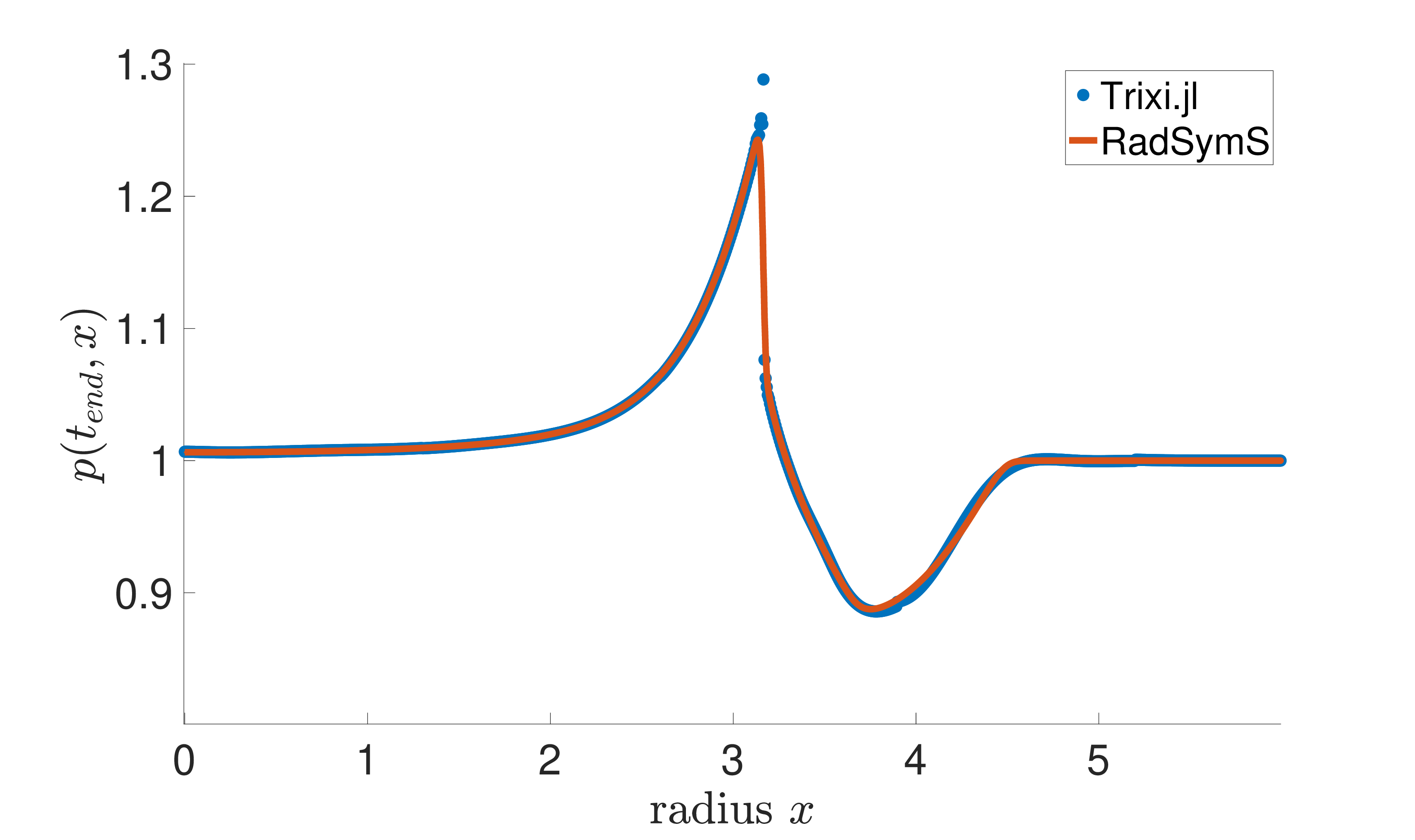}}
    \hspace{-1cm}
    \subfigure[Results for the velocity $v$ 3D]{
    \includegraphics[width=0.575\textwidth]{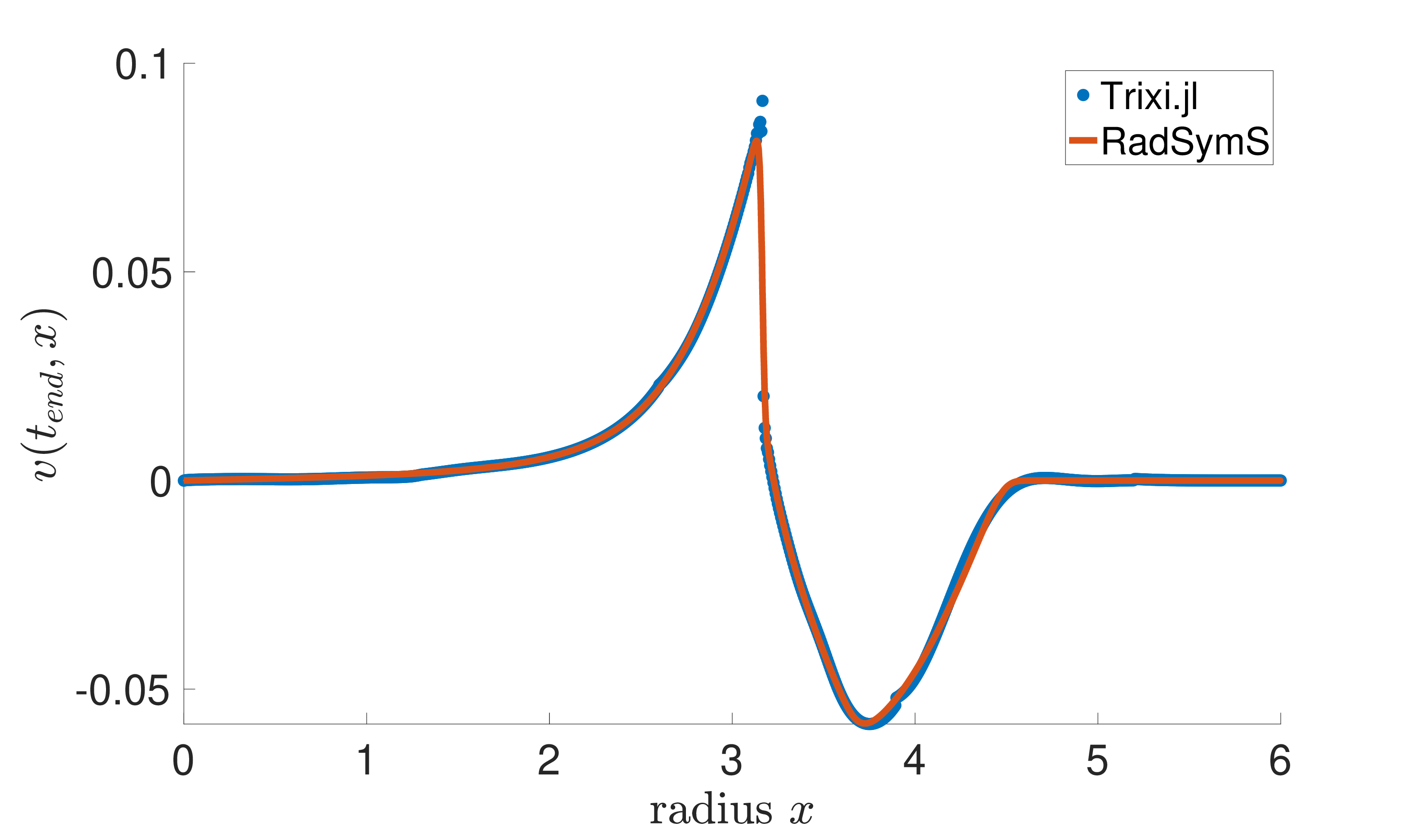}}

    \caption{\textbf{Example 4:}
    Results for the velocity $v$ and the pressure $p$ in 3D.
    Comparison for RadSymS (red) and Trixi.jl (blue) in radial direction at final time $t_{end} = 6$, see (a), (b), for pressure $p$ and velocity $v$.}
     \label{fig:example-4-3D-tend}
\end{figure}

\begin{figure}[h!]
    \subfigure[Results for the pressure $p$ in 2D]{
    \includegraphics[width=\textwidth]{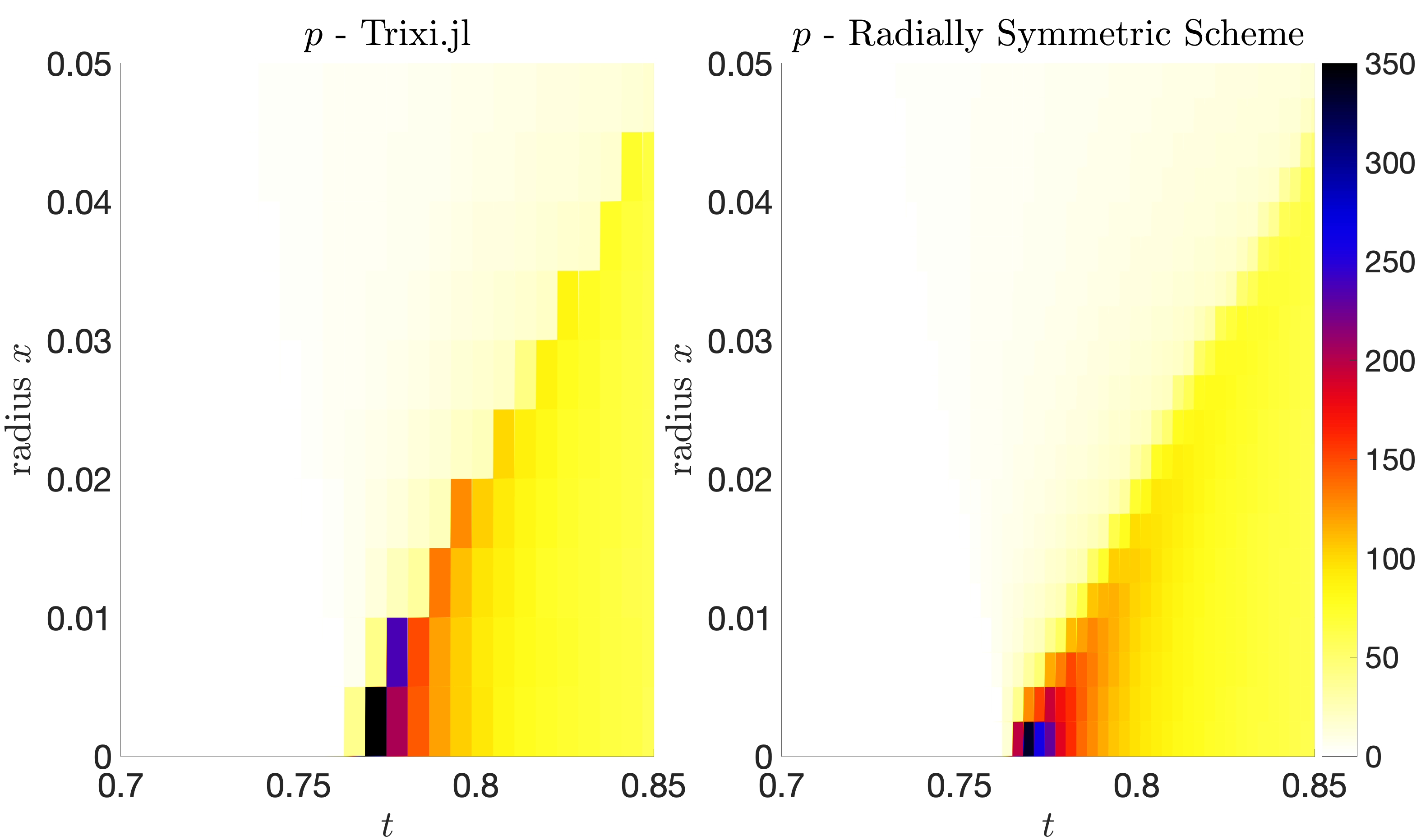}}\\
    \subfigure[Results for the velocity $v$ in 2D]{
    \includegraphics[width=\textwidth]{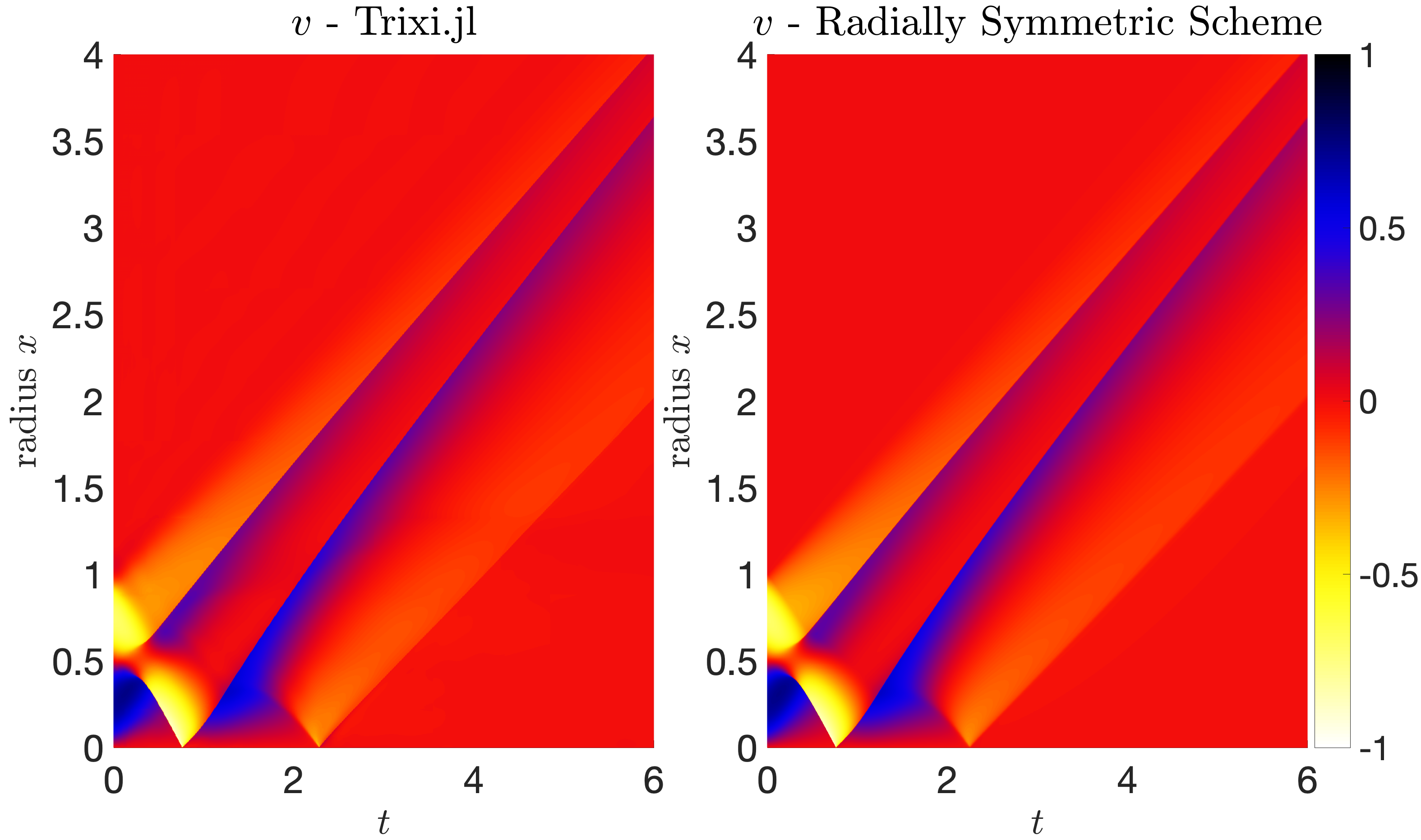}}\\
    \caption{\textbf{Example 5:}
    Results for the velocity $v$ and the pressure $p$ in 2D.
     Comparison for Trixi.jl (left) and RadSymS (right) in the $t$--$x$ plane, see (a), (b). For the pressure we present a zoom plot to visualize the small vicinity of the origin, where the pressure increases localized in time.}
    \label{fig:example-5-2D-tr}
\end{figure}

\begin{figure}[h!]
    \subfigure[Results for the pressure $p$ in 3D]{
    \includegraphics[width=\textwidth]{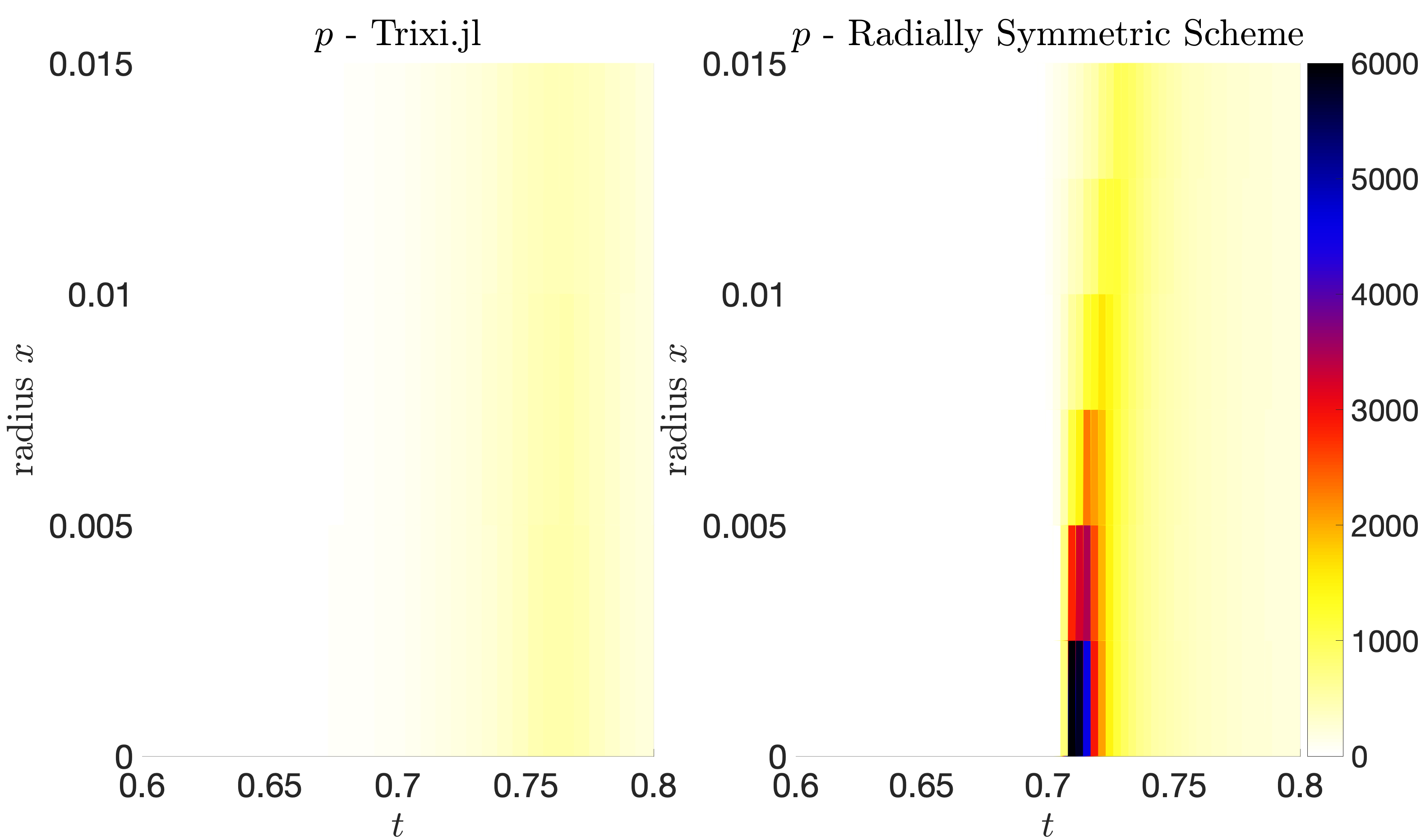}}\\
    \subfigure[Results for the velocity $v$ in 3D]{
    \includegraphics[width=\textwidth]{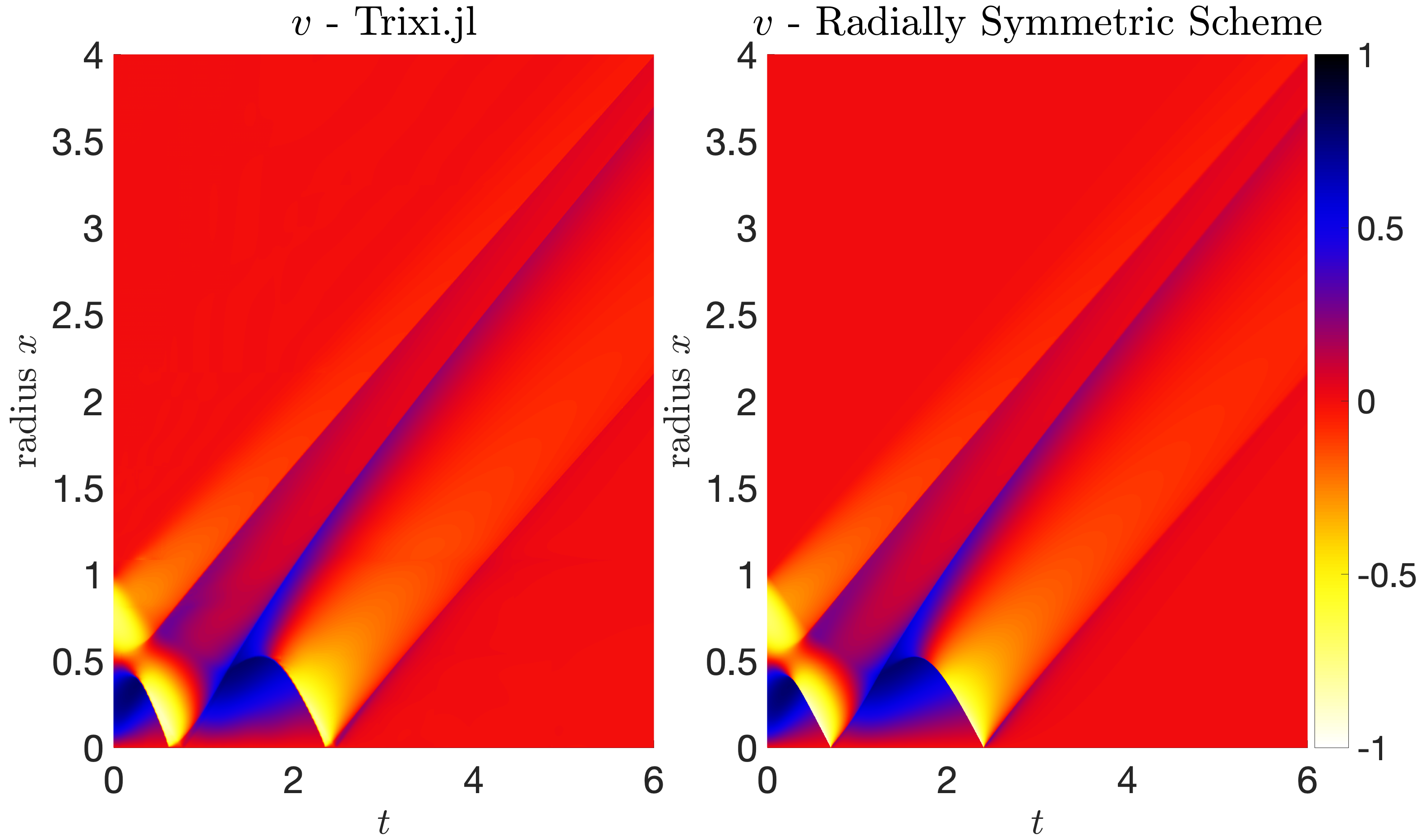}}\\
    \caption{\textbf{Example 5:}
    Results for the velocity $v$ and the pressure $p$ in 3D.
     Comparison for Trixi.jl (left) and RadSymS (right) in the $t$--$x$ plane, see (a), (b). For the pressure we present a zoom plot to visualize the small vicinity of the origin, where the pressure increases localized in time.}
    \label{fig:example-5-3D-tr}
\end{figure}

\begin{figure}[h!]
    \hspace{-1cm}
    \subfigure[Results for the pressure $p$ in 2D]{
    \includegraphics[width=0.575\textwidth]{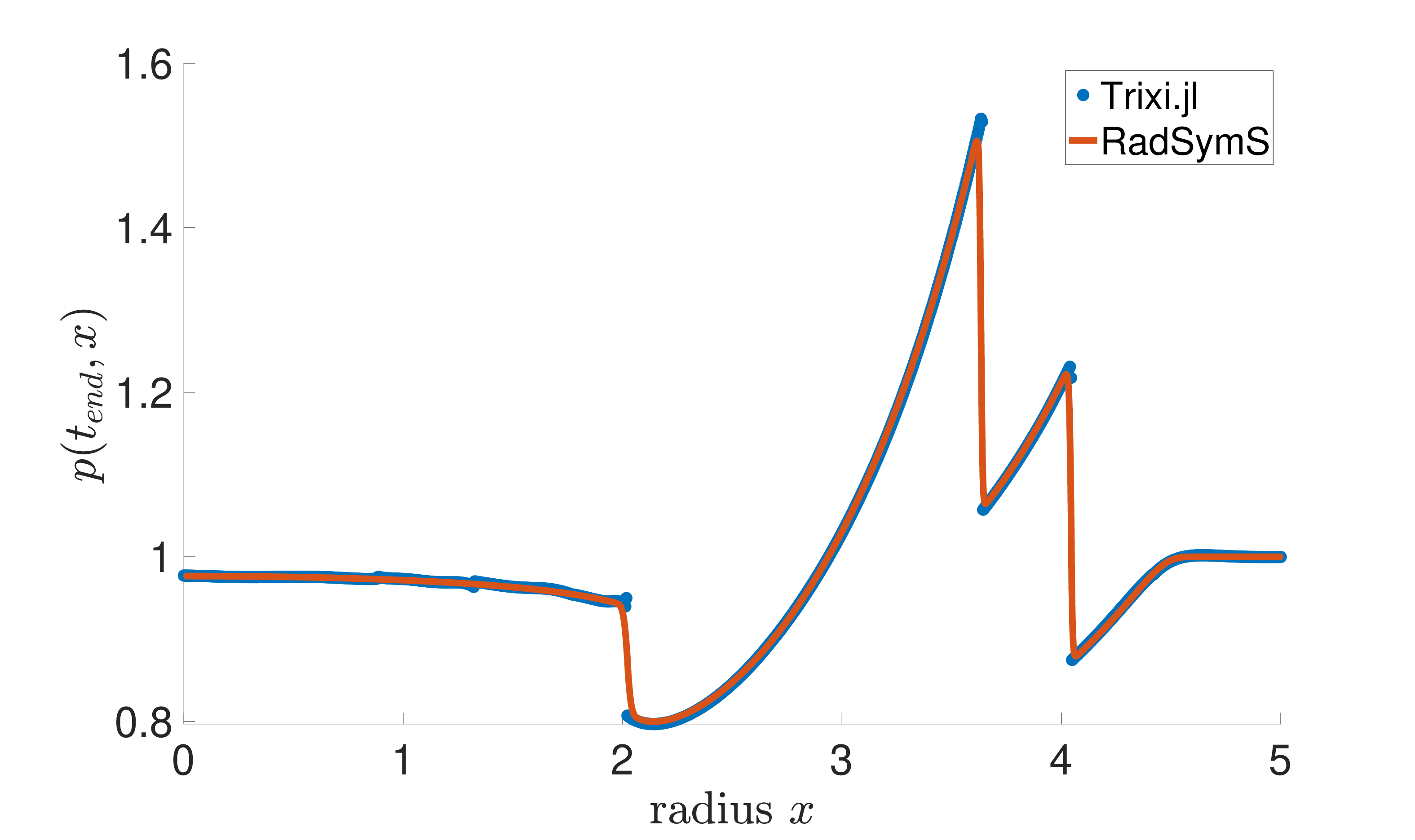}}
    \hspace{-1cm}
    \subfigure[Results for the velocity $v$ in 2D]{
    \includegraphics[width=0.575\textwidth]{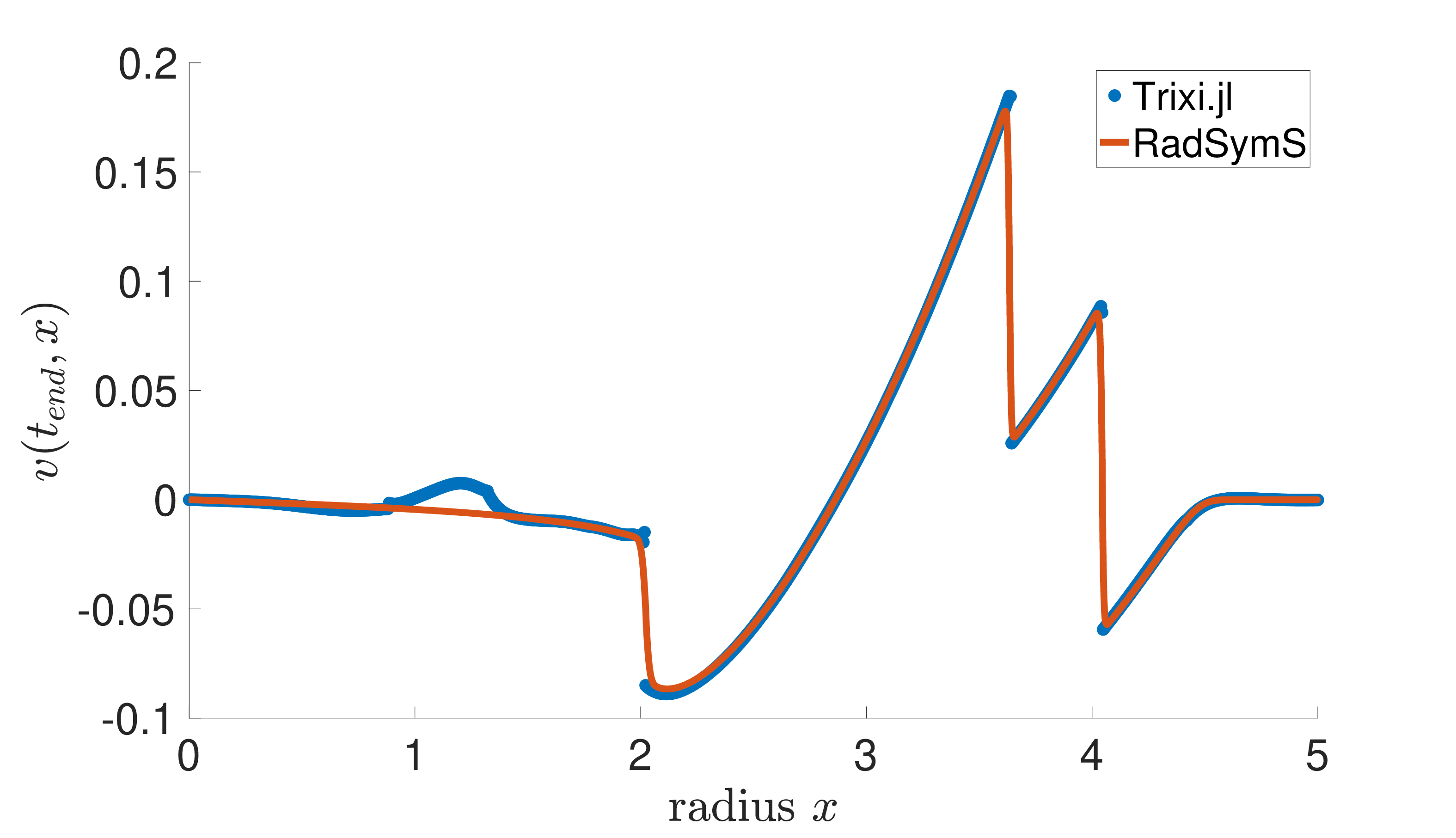}}
    \caption{\textbf{Example 5:}
    Results for the velocity $v$ and the pressure $p$ in 2D.
     Comparison for RadSymS (red) and Trixi.jl (blue) in radial direction at final time $t_{end} = 6$, see (a), (b), for pressure $p$ and velocity $v$. }
    \label{fig:example-5-2D-tend}
\end{figure}

\begin{figure}[h!]
    \hspace{-1cm}
    \subfigure[Results for the pressure $p$ in 3D]{
    \includegraphics[width=0.575\textwidth]{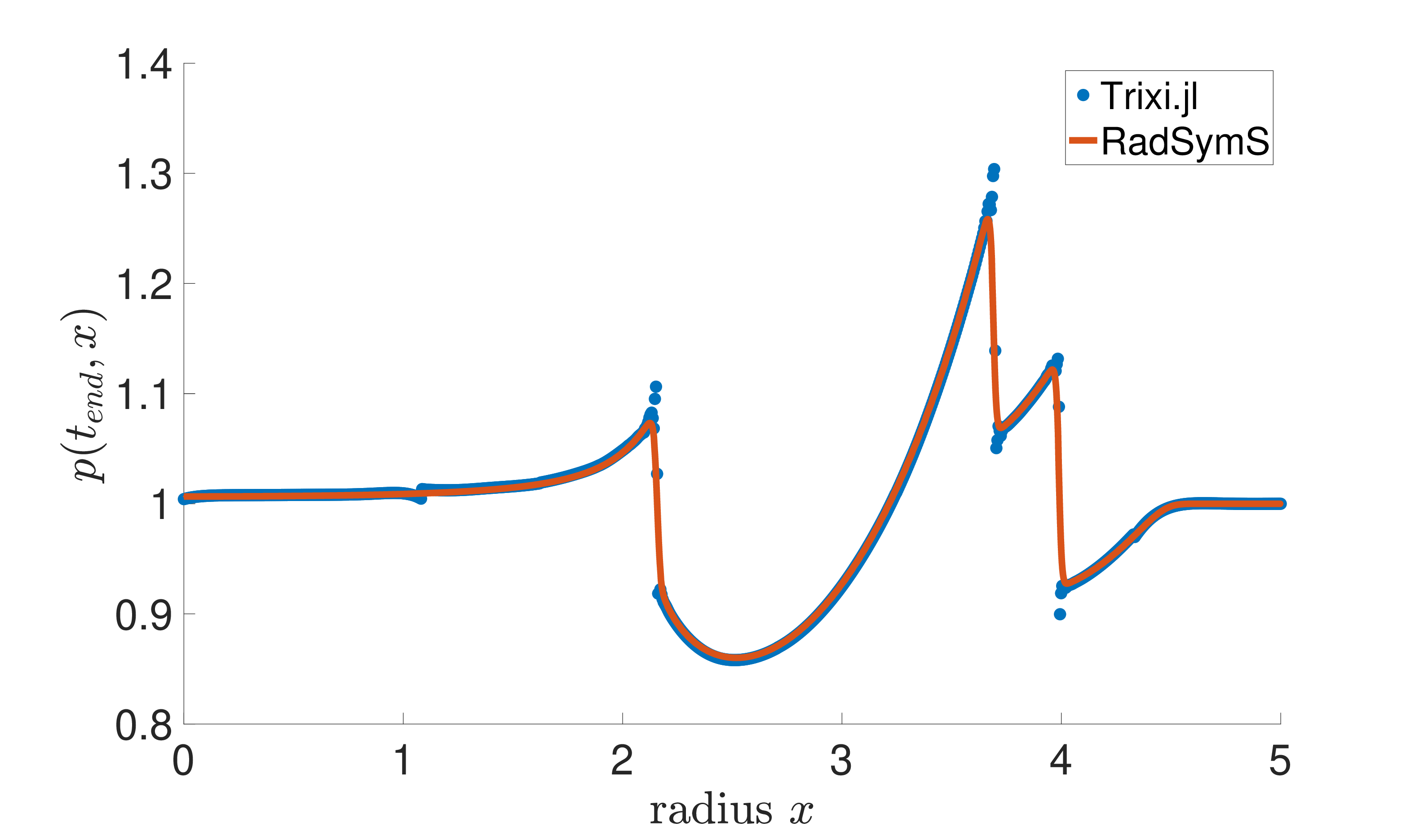}}
    \hspace{-1cm}
    \subfigure[Results for the velocity $v$ in 3D]{
    \includegraphics[width=0.575\textwidth]{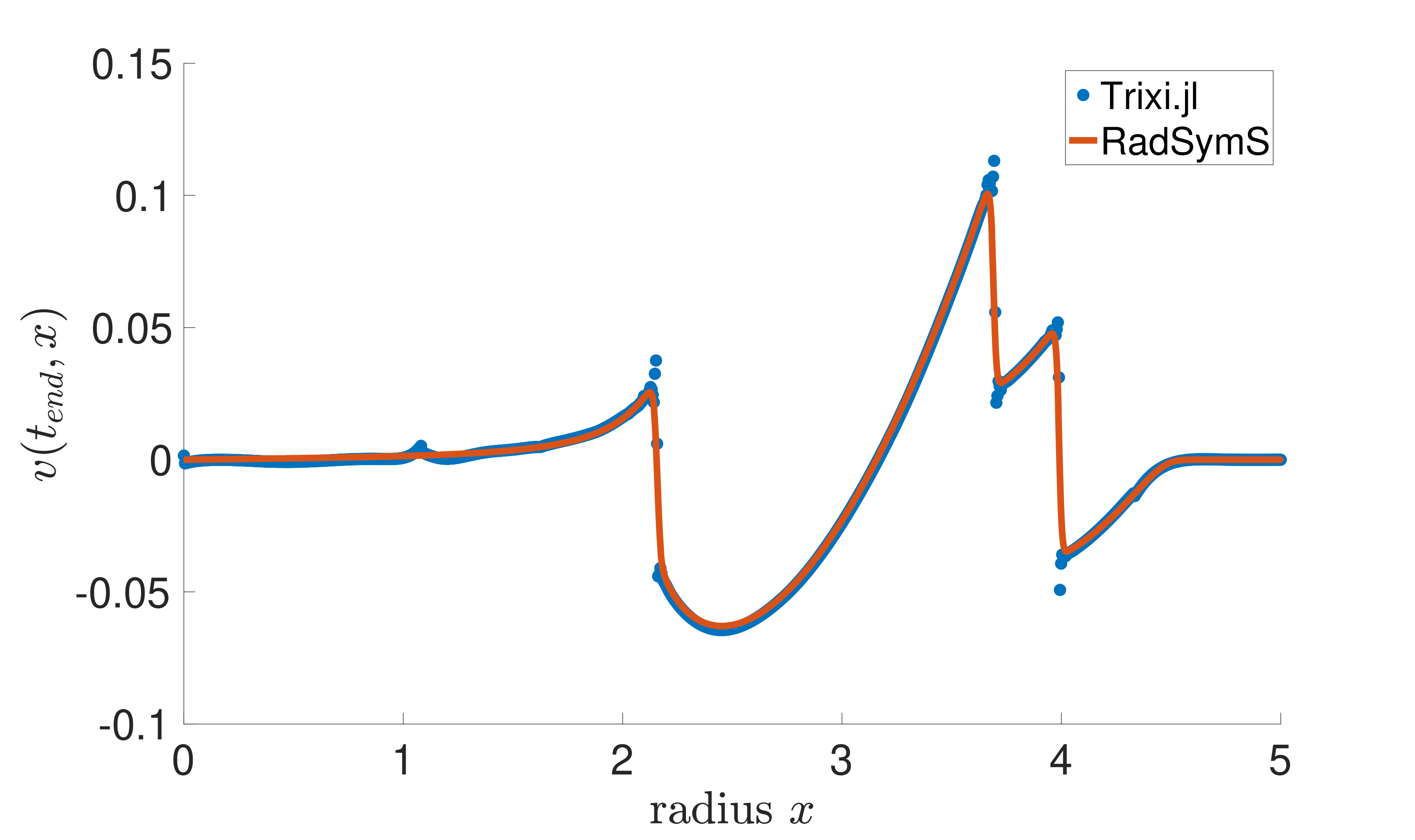}}
    \caption{\textbf{Example 5:}
    Results for the velocity $v$ and the pressure $p$ in 3D.
     Comparison for RadSymS (red) and Trixi.jl (blue) in radial direction at final time $t_{end} = 6$, see (a), (b), for pressure $p$ and velocity $v$. }
    \label{fig:example-5-3D-tend}
\end{figure}

In the following we provide the results for the benchmarks proposed and studied in \cite{kunik2024radially}.
We computed numerical solutions in two and three space dimensions using the genuine multi--dimensional solver Trixi.jl and compare the obtained results with the one-dimensional radially symmetric scheme RadSymS presented in \cite{Kunik2021,kunik2024radially}.
Note that all of the test cases are radially symmetric problems where two of them provide self-similar solutions, see Example 1 and 2.\\
In Figs.~\ref{fig:example-1-2D}--\ref{fig:example-5-3D-tend} for Examples 1--5, respectively,
we present 2D and 3D results in the $t$--$x$ plane and a line plot of the solution in radial direction at the final time $t_{end}$. We compare the results by means of the pressure $p$ and the velocity $v$ \eqref{eq:velocity_transformation}.
Additionally, for Examples 1 and 2 we compare the solutions also with the reference solution that can be determined for self-similar problems by solving the ODE \eqref{GL_ODE_sys} in the smooth part and using the jump conditions \eqref{shockvalues} to determine the stationary part in Example 1.\\
For all examples presented subsequently, we summarize the computational domain $\Omega\subset \R^d$ and
 the maximum refinement level $L\in\N$ in Table \ref{tab:parameter-trixi}. The initial cell size $\Delta\mathbf{x}_0 \in \R^d$ at the coarsest refinement level is always the whole computational domain $\Omega$. Due to the dyadic grid hierarchy in Trixi.jl the cell size at refinement level $l$ is $\Delta\mathbf{x}_l = 2^{-l}\Delta\mathbf{x}_0$, for $l=0,\dots,L$.
\begin{table}[h!]
    \centering
    \caption{Discretization parameters (domain $\Omega$ and maximum refinement level $L$) for Trixi.jl}
    \label{tab:parameter-trixi}
    \begin{tabular}{c|cccccccccc}
                    			& \multicolumn{2}{c}{Ex.~1}   & \multicolumn{2}{c}{Ex.~2}  & \multicolumn{2}{c}{Ex.~3} & \multicolumn{2}{c}{Ex.~4} & \multicolumn{2}{c}{Ex.~5}     \\
                    			& $d=2$ 		 & $d=3$ 	  & $d=2$	   & $d=3$		& $d=2$ 	 & $d=3$ 	  & $d=2$ 	   & $d=3$ 		& $d=2$ 	 & $d=3$\\
        \hline
        $\Omega$     			& \multicolumn{2}{c}{$[-2,2]^d$} & \multicolumn{2}{c}{$[-2,2]^d$} 	& \multicolumn{2}{c}{$[-6,6]^d$} & \multicolumn{2}{c}{$[-6,6]^d$} & \multicolumn{2}{c}{$[-5,5]^d$}\\
        $L$         			& 11 		 & 10 		  & 14		   & 13			& 13		 & 9		  & 13		   & 10 		& 14 		 & 9
    \end{tabular}
\end{table}
We want to remark that the refinement level in each 3D computation is less than the level used for the 2D computation of the same example, respectively.
Nevertheless, the computational time increases dramatically with the space dimension. Therefore, the final results for the pressure in examples 3, 4, and 5 in 3D cannot match the excellent agreement of the corresponding 2D results due to the computational costs and complexity. Here we decided to provide results that can be obtained within hours of computational time on a standard workstation rather than using HPC resources.\\
~\\
\textbf{\emph{Convergence test.}}
To verify the implementation and demonstrate the high-order accoracy of the multi--dimensional solver Trixi.jl, we perform a convergence test for a smooth manufactured solution of the ultra-relativistic Euler equations.
We choose the solution
\begin{align*}
    p = 2 + \sin\biggl(\sum_{j=1}^d x_j - \alpha_d t\biggr),
    \quad
    u_j = \frac{1}{2 \sqrt{d}},
\end{align*}
where $\alpha_2 = \sqrt{10}{4}$ and $\alpha_3 = \sqrt{15}{4}$.
The resulting source terms for the momentum equations are
\begin{align*}
    s_j = \frac{3}{4} \cos\biggl(\sum_{j=1}^d x_j - \alpha_d t\biggr),
\end{align*}
and no source term has to be added to the energy equation.
We use a uniformly refined mesh with periodic boundaries and compute the solution until the final time $t_{end} = 1$.
The time integration tolerance is set to $10^{-8}$ and the time step is adapted automatically by the step size controller of \cite{ranocha2021optimized} to ensure that the temporal error is negligible compared to the spatial error.
The resulting discrete $L^2$ errors and the corresponding experimental order of convergence (EOC) are presented in Tables~\ref{tab:convergence-trixi-2D} and \ref{tab:convergence-trixi-3D} for the 2D and 3D version of Trixi.jl, respectively.
As expected for polynomials of degree three, the EOC is close to four.

\begin{table}[htb]
\caption{Convergence results for the 2D version of Trixi.jl. The discrete $L^2$ errors are computed at the final time for the conserved variables.}
\label{tab:convergence-trixi-2D}
\centering
\begin{tabular}{c cc cc cc}
    \toprule
    Level & $L^2$ error $w_1$ & EOC $w_1$ & $L^2$ error $w_2$ & EOC $w_2$ & $L^2$ error $w_3$ & EOC $w_3$ \\
    \midrule
    2 &   \num{7.52e-03}  &           &   \num{7.52e-03}  &           &   \num{1.10e-02}  &           \\
    3 &   \num{4.87e-04}  &      3.95 &   \num{4.87e-04}  &      3.95 &   \num{8.51e-04}  &      3.69 \\
    4 &   \num{3.29e-05}  &      3.89 &   \num{3.29e-05}  &      3.89 &   \num{6.67e-05}  &      3.67 \\
    5 &   \num{1.93e-06}  &      4.09 &   \num{1.93e-06}  &      4.09 &   \num{3.43e-06}  &      4.28 \\
    6 &   \num{1.22e-07}  &      3.99 &   \num{1.22e-07}  &      3.99 &   \num{2.30e-07}  &      3.90 \\
    \bottomrule
\end{tabular}
\end{table}

\begin{table}[htb]
\caption{Convergence results for the 3D version of Trixi.jl. The discrete $L^2$ errors are computed at the final time for the conserved variables.}
\label{tab:convergence-trixi-3D}
\centering
\begin{tabular}{c cc cc cc cc}
    \toprule
Level & $L^2$ error $w_1$ & EOC $w_1$ & $L^2$ error $w_2$ & EOC $w_2$ & $L^2$ error $w_3$ & EOC $w_3$ & $L^2$ error $w_4$ & EOC $w_4$ \\
    \midrule
    2 &   \num{6.98e-03}  &           &   \num{6.98e-03}  &           &   \num{6.98e-03}  &           &   \num{1.18e-02}  &           \\
    3 &   \num{4.78e-04}  &      3.87 &   \num{4.78e-04}  &      3.87 &   \num{4.78e-04}  &      3.87 &   \num{1.00e-03}  &      3.56 \\
    4 &   \num{3.29e-05}  &      3.86 &   \num{3.29e-05}  &      3.86 &   \num{3.29e-05}  &      3.86 &   \num{7.48e-05}  &      3.74 \\
    5 &   \num{1.78e-06}  &      4.21 &   \num{1.78e-06}  &      4.21 &   \num{1.78e-06}  &      4.21 &   \num{4.08e-06}  &      4.20 \\
    6 &   \num{1.14e-07}  &      3.96 &   \num{1.14e-07}  &      3.96 &   \num{1.14e-07}  &      3.96 &   \num{2.48e-07}  &      4.04 \\
    \bottomrule
\end{tabular}
\end{table} 
\textbf{\emph{Example 1: Solutions including a shock and a stationary part.}}
Following \cite{Lai2019}, self-similar solutions can be constructed solving an ODE system in radial direction $x$. These solutions are in particular constant along rays $\xi=x/t$ or, equivalently, $\vartheta=t/x$. \\
We consider constant initial data with pressure $p_0=1$ and radial velocity $v_0 \in (-1,0)$.
Due to \cite[Section 2.3]{Lai2019} there is a solution $p(t,x)=P(\vartheta)$ and $v(t,x)=V(\vartheta)$ depending only on $\vartheta=t/x$ for $t,x >0$,
with a single straight line shock emanating from the zero point. Let $\tilde{s} \in (0,1)$ be the unknown constant shock speed.
Then we put $v_-=0$ and can find an unknown pressure $p_->0$ with
\[
  p(t,x)=p_-\,, \quad v(t,x)=v_-=0 \quad \mbox{for~} x < \tilde{s}\cdot t.
\]
Due to the Rankine-Hugoniot shock conditions introduced in \cite[Section 4.4]{Kunik2005} we obtain from \cite[page 82]{Kunik2005}
for a so called 3-shock after a lengthy calculation the algebraic shock conditions
\begin{equation}\label{pplus}
  \frac{p_+}{p_-}=\frac{9\tilde{s}^2-1}{3(1-\tilde{s}^2)},\quad v_+=\frac{3}{2}\tilde{s}-\frac{1}{2\tilde{s}}\,, \quad v_-=0,
\end{equation}
\begin{equation*}\label{ungl}
  0 < p_+ < p_-,\quad \frac13 < \tilde{s} < \frac{1}{\sqrt{3}}.
\end{equation*}
Here $p_-,v_-$ are the values of pressure and velocity left to the 3-shock, and $p_+,v_+$ are the values of pressure and velocity right to the 3-shock.
Due to Lai \cite[Section 2.3]{Lai2019} the solution $p=P(\vartheta)$, $v=V(\vartheta)$ satisfies the IVP of ordinary differential equations
\begin{equation}\label{GL_ODE_sys}
   \begin{split}
   \dot{V}(\vartheta) &= (d-1)\frac{V(V-\vartheta)(1-V^2)}{3(\vartheta V-1)^2-(V-\vartheta)^2},\quad V(0)=v_0 \in (-1,0),\\
   \dot{P}({\vartheta}) &= (d-1)\frac{4PV(\vartheta V-1)}{3(\vartheta V-1)^2-(V-\vartheta)^2},\quad P(0)=p_0=1.
   \end{split}
\end{equation}
Moreover, it is shown in the appendix of \cite{kunik2024radially} that Lai's results guarantee a unique solution for $0 \leq \vartheta < \vartheta_{max}$
with a certain value $\vartheta_{max} >\sqrt{3}$ and a unique value $\sqrt{3} < \tilde{\vartheta} < \sqrt{v_0^2 + 3} - v_0 < 3$, see \cite[(A.10)]{kunik2024radially} such that
\begin{equation}\label{vvalue}
   V(\tilde{\vartheta})=\frac{3}{2\tilde{\vartheta}}-\frac{\tilde{\vartheta}}{2}.
\end{equation}
After having determined $\tilde{\vartheta}$ we finally obtain
\begin{equation}\label{shockvalues}
   \tilde{s}=\frac{1}{\tilde{\vartheta}},\quad v_+=V(\tilde{\vartheta}),\quad p_+=P(\tilde{\vartheta}),\quad p_-=p_+\frac{3(1-\tilde{s}^2)}{9\tilde{s}^2-1}.
\end{equation}
For the numerical simulation we prescribe a constant initial pressure $p_0 = 1$ and a constant initial velocity $v_0 = -1/\sqrt{2}$.
This corresponds to the following initial data for the original IVP \eqref{orisys}
\[
  \tilde{p}(0,\mathbf{x}) = 1\quad\text{and}\quad\mathbf{u}(0,\mathbf{x}) = -\frac{\mathbf{x}}{|\mathbf{x}|}\,, \quad\mathbf{x}\in\R^d \setminus \{\mathbf{0} \}.
\]
The radially symmetric solution of \eqref{orisys} is determined by the solution of the ODE system \eqref{GL_ODE_sys} which gives a shock wave moving at constant speed $\tilde{s}$. We want to point out that these ODEs are written in terms of $\vartheta = t/x$. Thus, the initial data $P(0)$ and $V(0)$ for \eqref{GL_ODE_sys} prescribe the solution for \eqref{orisys} at infinity given by $p_0 = P(0)$ and $v_0 = V(0)$.
The ODE system \eqref{GL_ODE_sys} is solved by applying the classical fourth-order RK scheme with step size $h = 10^{-6}$.
For the computation of   $\tilde{\vartheta}$ the ODE solver is run until \eqref{vvalue}
 is satisfied with a tolerance of $\varepsilon = 10^{-9}$.
The solution values for the shock speed, the states ahead and behind of the shock, respectively, are presented in Table \ref{tab:shock_dat_ex1}.
\begin{table}[h!]
    %
    \caption{Shock states for $d=2,3$ in Example 1}
    \label{tab:shock_dat_ex1}
    \centering
    \begin{tabular}{c|ccccc}
                    & $\tilde{s}$   & $p_-$         & $v_-$ & $p_+$         & $v_+$     \\
        \hline
        $d = 2$     & $0.45503$     & $15.75505$    & $0$   & $\hphantom{1}5.71869$    & $-0.41629$ \\
        $d = 3$     & $0.52314$     & $25.56463$    & $0$   & $17.16524$    & $-0.17106$
    \end{tabular}
\end{table}
The results reported in Fig.~\ref{fig:example-1-2D} (2D) and Fig.~\ref{fig:example-1-3D} (3D) clearly show an excellent agreement of both numerical methods with the reference solution provided by an ODE solver for the final time $t_{end} = 1$ in radial direction. Due to the self-similar structure of the solution we omit the presentation of results in the $t$--$x$ plane.
%
%

In order to verify the entropy conservation property of the presented scheme we provide a numerical test for slightly modified initial data of the present example, i.e., we use a constant initial pressure $p_0 = 1$ and a constant initial velocity $v_0 = -1/\sqrt{26}$.
This corresponds to the following initial data for the original IVP \eqref{orisys}
\[
  \tilde{p}(0,\mathbf{x}) = 1\quad\text{and}\quad\mathbf{u}(0,\mathbf{x}) = -\frac{1}{5}\frac{\mathbf{x}}{|\mathbf{x}|}\,, \quad\mathbf{x}\in\R^d \setminus \{\mathbf{0} \}.
\]
To verify the entropy conservation property, we compute the entropy rate, i.e., the discrete version of
\[
    \int_\Omega \nabla_\mathbf{w} \eta(\mathbf{w}) \cdot \partial_t \mathbf{w} \,\dd\mathbf{x},
\]
where the integral is replaced by the quadrature associated to the DGSEM, and the time derivative is given by the entroy-conservative semidiscretization.
The results for the entropy rate over time in 2D and 3D are given in Fig.~\ref{fig:entropy_rate}. Clearly entropy is conserved up to machine precision.
Checking discrete entropy conservation like this is a demanding test for the derivation of the entropy-conservative flux in Theorem~\ref{thm:entropy-conservative-flux} as well as for its implementation.

\FloatBarrier
\textbf{\emph{Example 2: Self-Similar Expansion.}\label{selfsim_expansion}}
In this example we prescribe initial data leading to a smooth self-similar solution
by applying Lai's approach as in Example 1.
Note that for $d=3$ the solution will expand into vacuum at the speed of light when the initial value for $v_0$ is close enough to one.
We do not consider this case here and refer to \cite{Lai2019} for more details.
For the numerical simulation we consider a constant initial pressure $p_0 = 1$ and a constant initial velocity $v_0 = 1/\sqrt{5}$.
This corresponds to the following initial data for the original IVP \eqref{orisys}
\[
  \tilde{p}(0,\mathbf{x}) = 1\quad\text{and}\quad\mathbf{u}(0,\mathbf{x}) = \frac{1}{2}\frac{\mathbf{x}}{|\mathbf{x}|}\,, \quad\mathbf{x}\in\R^d \setminus \{\mathbf{0} \}.
\]
The solution consists of a rarefaction wave that is determined by the ODE \eqref{GL_ODE_sys}. In contrast to Example 1 there is no shock.\\
The results reported in Fig.~\ref{fig:example-2-2D} for 2D and Fig.~\ref{fig:example-2-3D} for 3D, respectively, clearly show an excellent agreement of both numerical methods with the reference solution provided by an ODE solver for the final time $t_{end} = 1$ in radial direction. Due to the self-similar structure of the solution we omit the presentation of results in the $t$--$x$ plane.
\FloatBarrier
\textbf{\emph{Example  3: Expansion of a Spherical Bubble.}}
Next we consider the expansion of a spherical bubble with the following initial data
\begin{equation*}
  p_0(x) = \begin{cases}
      1 \quad & \mbox{for~} 0 \leq x \leq 1\\
      0.1 \quad & \mbox{for~} x > 1,
  \end{cases}
  \qquad v_0(x) = 0.
\end{equation*}
These initial values correspond to the following initial values for the original IVP \eqref{orisys}
\begin{equation*}
  \tilde{p}(0,\mathbf{x}) = \begin{cases}
      1 \quad & \mbox{for~} 0 \leq |\mathbf{x}| \leq 1\\
      0.1 \quad & \mbox{for~} |\mathbf{x}| > 1,
  \end{cases}
  \qquad \mathbf{u}(0,\mathbf{x}) = \mathbf{0}.
\end{equation*}
Initially, the pressure inside the bubble is ten times larger than outside, which leads to a fast expansion of the bubble into the outer low pressure area.
This in turn gives rise to the formation of another low pressure area.
We observe the formation of a shock wave, running towards the new low pressure area and reaching the zero point around time $t=5.032$ for $d=2$ and $t = 4.165$ for $d=3$.
The formation of this new shock wave is a peculiar nonlinear phenomenon.
Shortly before the shock reaches the zero point the pressure takes very low values, but its reflection from the zero point depicted in
Fig.~\ref{fig:example-3-2D-tr} and Fig.~\ref{fig:example-3-3D-tr}
indicates a blow-up of the pressure in a very small time-space range near the boundary.
This is similar to \cite[Section 5, Example 4]{Kunik2021}, but here the blow-up is weaker than in three space dimensions, and its illustration requires a higher numerical resolution.
A similar behavior for the corresponding solution of the explosion problem with the classical Euler equations should be expected.

Before discussing the numerical results in details, we show the evolution of the total entropy $\int_\Omega \eta(\mathbf{w}(t,\mathbf{x}))\,\dd\mathbf{x}$ (computed numerically via the quadrature rule associated to the method) over time for the present example in Fig.~\ref{fig:entropy}.
In contrast to the entropy results shown earlier in Fig.~\ref{fig:entropy_rate}, we use the the local Lax-Friedrichs/Rusanov flux and do not show the discrete entropy rate, but the total entropy itself.
Since the entropy $\eta$ is concave, entropy stability means that the total entropy increases over time.
Clearly, the numerical results shown in Fig.~\ref{fig:entropy} confirm the entropy stability of the presented scheme, as expected based on the theory of flux differencing and our construction of the entropy-conservative flux in Theorem~\ref{thm:entropy-conservative-flux}.

From the depicted solution at time $t=6$ we observe the reflected shock curve at radius around $x=0.55$ for $d=2$ and near $x = 1.21$ for $d=3$, respectively.
The results in the $t$--$x$ plane given in Fig.~\ref{fig:example-3-2D-tr} and Fig.~\ref{fig:example-3-3D-tr} show the expansion of the bubble and the formation of the new shock which focuses at the origin and then is reflected again. At the focus point the pressure rises drastically in a small vicinity of the origin and away from that the pressure values do not vary a lot. Therefore, we present a zoom plot of the pressure. We observe an excellent agreement of both numerical methods for the velocity in Figs.~\ref{fig:example-3-2D-tr} and \ref{fig:example-3-3D-tr} (b). For the pressure it is visible in 2D that Trixi.jl is capable of resolving the solution close to the reflection point of the shock at the origin rather well, see Fig.~\ref{fig:example-3-2D-tr} (a).
However, in 3D the reflection point is not captured well, which is due to the computational efforts needed in 3D as stated before.
Indeed, the maximum pressure value computed by Trixi.jl is $p_{Trixi} \approx 9.6562$ versus $p_{RadSymS} \approx 289.2772$ for the case in 3D, see Fig.~\ref{fig:example-3-3D-tr} (a).
This is due to the lower resolution we for the 3D results compared to the 2D results, see Table~\ref{tab:parameter-trixi}.
Matching the maximum pressure value also in 3D would require at least the same resolution as in 2D, which was not feasible for the workstation we could use for this project. Further note that this pressure peak occurs in the very vicinity of the origin and is very much localized in time. Thus we present zoom plot for the pressure to visualize it. For the velocity this is not needed since we have proper scaling properties due to $|v| <1$.

Additionally, we compare both methods at final time $t_{end} = 6$ where the solutions coincide very well, see
Figs.~\ref{fig:example-3-2D-tend} and \ref{fig:example-3-3D-tend}.
\FloatBarrier
\textbf{\emph{Example 4: Collapse of a Spherical Bubble.}}
Next we study the collapse of a spherical bubble with the following initial data
\begin{equation*}
  p_0(x) = \begin{cases}
      0.1 \quad & \mbox{for~} 0 \leq x \leq 1\\
      1 \quad & \mbox{for~} x > 1,
  \end{cases}
  \qquad v_0(x) = 0.
\end{equation*}
These initial values correspond to the following initial values for the original IVP \eqref{orisys}
\begin{equation*}
  \tilde{p}(0,\mathbf{x}) = \begin{cases}
      0.1 \quad & \mbox{for~} 0 \leq |\mathbf{x}| \leq 1\\
      1 \quad & \mbox{for~} |\mathbf{x}| > 1,
  \end{cases}
  \qquad \mathbf{u}(0,\mathbf{x}) = \mathbf{0}.
\end{equation*}
The results in the $t$--$x$ plane, see Figs.~\ref{fig:example-4-2D-tr} and \ref{fig:example-4-3D-tr},
show the collapse of the bubble with a focus point at the origin which then is reflected again. At the focus point the pressure rises drastically in a small vicinity of the origin and away from that the pressure values do not vary a lot. Thus, we again present a zoom plot of the pressure for this example.
We observe for both, $d=2$ and $d=3$, an excellent agreement of both numerical methods for the velocity in Figs.~\ref{fig:example-4-2D-tr} and \ref{fig:example-4-3D-tr}, respectively. For the pressure it is again visible that in 2D the results agree well whereas due to the computational complexity in three dimensions the pressure peak is not very well resolved.
Indeed, the maximum pressure value computed by Trixi.jl is $p_{Trixi} \approx 390$ versus $p_{RadSymS} \approx 3550$ for the case in 3D, see Fig.~\ref{fig:example-4-3D-tr} (a).
This is again caused by the lower resolution we for the 3D results compared to the 2D results, see Table~\ref{tab:parameter-trixi}.
Matching the maximum pressure value also in 3D would require at least the same resolution as in 2D, which was not feasible for the workstation we could use for this project.

Additionally, we compare both methods at time $t_{end} = 6$
where the solutions coincide very well, see Figs.~\ref{fig:example-4-2D-tend} and \ref{fig:example-4-3D-tend}.
\FloatBarrier
%
%
\textbf{\emph{Example 5: Initially Sine--Shaped Radial Velocity.}}
Finally, we study initial data with a velocity in radial direction described by a sine function, i.e.,
\begin{align*}
    p_0(x) = 1,\quad u_0(x) = \begin{cases}
        \sin(2\pi x),\;&x < 1,\\
        0,\;&x \geq 1.
    \end{cases}
\end{align*}
These correspond to the following initial data for the multi--dimensional simulation
\begin{align*}
    \tilde{p}(0,\mathbf{x}) = 1,\quad \mathbf{u}(0,\mathbf{x}) = \begin{dcases}
        \sin(2\pi |\mathbf{x}|)\frac{\mathbf{x}}{|\mathbf{x}|},\;&|\mathbf{x}| < 1,\\
        0,\;&\mathbf{x} \geq 1.
    \end{dcases}
\end{align*}
The results in the $t$--$x$ plane, see Figs.~\ref{fig:example-5-2D-tr} and \ref{fig:example-5-3D-tr},
show a complex wave structure for the velocity and a high pressure focus at around time $t=0.77$ for $d=2$ and $t=0.71$ for $d=3$, respectively. Again, we present a zoom plot of the pressure for this example.
As for the previous examples it is again visible near the reflection point that in 2D the results agree well whereas due to the computational complexity in three dimensions the pressure peak is not very well resolved.
Indeed, the maximum pressure value computed by Trixi.jl is $p_{Trixi} \approx 503$ versus $p_{RadSymS} \approx 5957$ for the case in 3D, see Fig.~\ref{fig:example-5-3D-tr} (a).
This is again due to the lower resolution we for the 3D results compared to the 2D results, see Table~\ref{tab:parameter-trixi}.
Matching the maximum pressure value also in 3D would require at least the same resolution as in 2D, which was not feasible for the workstation we could use for this project.

We observe an excellent agreement of both numerical methods for the velocity in Figs.~\ref{fig:example-5-2D-tr} and \ref{fig:example-5-3D-tr} (b). Additionally, we compare both methods at final time $t_{end} = 6$ where the solutions coincide very well, see Figs.~\ref{fig:example-5-2D-tend} and \ref{fig:example-5-3D-tend}.
\FloatBarrier
\section{Conclusion}\label{sec:concl}
We have investigated the system of the ultra--relativistic Euler equations in two and three space dimensions. For the first time an entropy-conservative flux has been derived and applied for numerical computations. We studied five different benchmarks in two and three dimensions, where the full-dimensional system is treated using Trixi.jl and the results are compared to numerical reference solutions. It should be emphasized that this is the first time that genuine three-dimensional results for these benchmarks have been obtained and presented.
The obtained numerical results coincide well, although they emphasize the need to reduce the computational complexity for the three-dimensional case. Despite the interest in the particular system under consideration, these benchmarks may serve as challenging benchmarks for other multi--dimensional solvers.
It will be interesting to investigate the extension of the present work to different space--times.
\appendix
\section*{Acknowledgments}
Funded by the Deutsche Forschungsgemeinschaft (DFG, German Research Foundation)
through
SPP 2410 Hyperbolic Balance Laws in Fluid Mechanics: Complexity, Scales, Randomness (CoScaRa) within the Projects
525939417 (FT)
and 526031774 (HR).
HR additionally acknowledges the DFG individual research grants 513301895 and 528753982.
FT thanks Matthias Kunik for introducing him to the interesting topic of relativistic Euler equations and fruitful discussions.
HR thanks Joshua Lampert for reviewing HR's contributions to Trixi.jl that improved the efficiency of radially symmetric output data.
\section*{Data Availability Statement}
All Julia source code and data needed to reproduce the numerical results presented in this paper are available in the accompanying reproducibility repository \cite{thein2025computingRepro}.
\section{Derivatives and Supplementary Calculations for the Entropy}\label{app:entropy}
In this section we want to provide further details on the derivatives. Furthermore, we show the concavity of the entropy and prove relation \eqref{entropy_flux_rel}.
For simplicity, we will express the derivatives in terms of the primitive variables, which is more convenient to work with.\\
\newline
Using \eqref{eq:press} and \eqref{eq:cons-var} we can calculate the derivatives of the pressure
\begin{align}
    \label{eq:aux-pres-der-1}
      \partial_{w_j} p =
      \begin{cases}
       - \dfrac{w_j}{w_{d+1}+ 3p} = -2u_j\dfrac{\sqrt{1 + |\mathbf{u}|^2}}{3 + 2|\mathbf{u}|^2},\, &j=1,\ldots,d, \\[1em]
           \dfrac{w_{d+1}-p}{w_{d+1}+ 3p} = \dfrac{1 + 2|\mathbf{u}|^2}{3 + 2|\mathbf{u}|^2},\, &j=d+1.
       \end{cases}
\end{align}
The derivatives of the velocity components $u_i, i=1,\dots,d$ can be obtained using \eqref{eq:velo} and \eqref{eq:cons-var}
\begin{align}
    \begin{split}\label{eq:aux-velo-der-1}
      \partial_{w_j} u_i &=
       \dfrac{\delta_{ij}}{\sqrt{4p(w_{d+1} + p)}} + \dfrac{1}{2}\dfrac{u_i}{p}\dfrac{w_j(w_{d+1} + 2p)}{(w_{d+1} + p)(w_{d+1}+ 3p)}\\
       &= \dfrac{\delta_{ij}}{4p\sqrt{1 + |\mathbf{u}|^2}} + \dfrac{u_iu_j}{4p}\dfrac{5 + 4|\mathbf{u}|^2}{\sqrt{1 + |\mathbf{u}|^2}(3 + 2|\mathbf{u}|^2)},\, j=1,\ldots,d, \\[1.5em]
      \partial_{w_{d+1}} u_i &= -\dfrac{1}{2}\dfrac{u_i}{p}\dfrac{w_{d+1} + p}{w_{d+1} + 3p} = -\dfrac{u_i}{p}\dfrac{1 + |\mathbf{u}|^2}{3 + 2|\mathbf{u}|^2}.
    \end{split}
\end{align}
We yield for the derivatives of the entropy with respect to the primitive variables
\begin{align}
	\frac{\partial\eta}{\partial p} &= \frac{3}{4}p^{-1/4}\sqrt{1 + |\mathbf{u}|^2} = \frac{3}{4}\frac{\eta}{p},\label{deriv_eta_p}\\
	\frac{\partial\eta}{\partial u_i} &= p^{3/4}\frac{u_i}{\sqrt{1 + |\mathbf{u}|^2}} = \eta\frac{u_i}{1 + |\mathbf{u}|^2}.\label{deriv_eta_u}
\end{align}
Note that corresponding to the flux potential \eqref{eq:flux_potential} we can also calculate the entropy potential. Sometimes it is also called generating potential; we refer to the recent survey \cite{Warnecke2024} and the references therein.
We obtain for the entropy potential
\begin{align}
    \phi = \nabla_\mathbf{w}\eta\cdot\mathbf{w} - \eta = -\frac{1}{4}\eta.\label{eq:main_field_pot}
\end{align}
Next, we provide the second order derivatives of $\eta$ with respect to the conservative variables, i.e.,
\begin{align}
	\partial_{w_i}^2\eta &= -\frac{\left(1 + u_i^2\right)\left(3 + 2|\mathbf{u}|^2\right) + u_i^2\left(1 + |\mathbf{u}|^2\right)}{16p^{\frac{5}{4}}\sqrt{1 + |\mathbf{u}|^2}\left(3 + 2|\mathbf{u}|^2\right)},\;i=1,\dots,d,\\
	\partial_{w_iw_j}^2\eta &= -\frac{u_iu_j\left(4 + 3|\mathbf{u}|^2\right)}{16p^{\frac{5}{4}}\sqrt{1 + |\mathbf{u}|^2}\left(3 + 2|\mathbf{u}|^2\right)},\;i,j=1,\dots,d,\;i\neq j,\\
	\partial_{w_iw_{d+1}}^2\eta &= -\frac{u_i\left(5 + 6|\mathbf{u}|^2\right)}{16p^{\frac{5}{4}}\left(3 + 2|\mathbf{u}|^2\right)},\;i=1,\dots,d,\\
	\partial_{w_{d+1}}^2\eta &= -\frac{\sqrt{1 + |\mathbf{u}|^2}\left(1 + 6|\mathbf{u}|^2\right)}{16p^{\frac{5}{4}}\left(3 + 2|\mathbf{u}|^2\right)}.
\end{align}
For simplicity, we are going to show the convexity of $\tilde{\eta} = -\eta$ and the Hessian can be written as follows
\begin{align}
	\mathbf{D}^2_\mathbf{w}\tilde{\eta} = \frac{1}{16p^{\frac{5}{4}}\sqrt{1 + |\mathbf{u}|^2}(3 + 2|\mathbf{u}|^2)}(\mathbf{A}_1 + \mathbf{A}_2 + \mathbf{A}_3)
\end{align}
with
\begin{align*}
	\mathbf{A}_1 &= \left(3 + 2|\mathbf{u}|^2\right)
	\begin{pmatrix}
		& 						  	& 0\\
		& \mathbf{Id}_{d\times d}	& \vdots\\
	0	& \dots						& 0
	\end{pmatrix},\quad
	\mathbf{A}_2 = \left(4 + 3|\mathbf{u}|^2\right)
	\begin{pmatrix}
	u_1^2	& u_1u_2	& \dots	 & u_1u_d		& 0\\
	u_1u_2	& u_2^2		& \dots	 & u_2u_d		& 0\\
	\vdots	&			& \ddots & \vdots		& \vdots\\
	u_1u_d	& \dots		& \dots	 & u_d^2			& 0\\
	0		& \dots		& \dots 	 & 0			& 0
	\end{pmatrix}\\
	\text{and}\quad
	\mathbf{A}_3 &=
	\begin{pmatrix}
	0	& \dots	& 0	& u_1\left(5 + 6|\mathbf{u}|^2\right)\sqrt{1 + |\mathbf{u}|^2}\\
	\vdots	& 	& \vdots	 & \vdots\\
	0	& \dots		& 0	 & u_d\left(5 + 6|\mathbf{u}|^2\right)\sqrt{1 + |\mathbf{u}|^2}	\\
	u_1\left(5 + 6|\mathbf{u}|^2\right)\sqrt{1 + |\mathbf{u}|^2}		& \dots		& u_d\left(5 + 6|\mathbf{u}|^2\right)\sqrt{1 + |\mathbf{u}|^2}	 	 & \left(1 + |\mathbf{u}|^2\right)\left(1 + 6|\mathbf{u}|^2\right)
	\end{pmatrix}
\end{align*}
To prove convexity we now want to show that
\[
  \mathbf{v}^T\ \mathbf{D}^2_\mathbf{w}\tilde{\eta}\, \mathbf{v} > 0,\;\forall\mathbf{v}\in\R^{d + 1}\setminus\{\mathbf{0}\}.
\]
Therefore, we first note that all elements $(v_1,\dots,v_d,0)^T\in\R^{d+1}$ belong to $\ker(\mathbf{A}_3)$ and conversely all elements $(0,\dots,0,v_{d+1})^T\in\R^{d+1}$ belong to $\ker(\mathbf{A}_1)$ as well as to $\ker(\mathbf{A}_2)$.
Clearly, $\mathbf{A}_1 > 0$ for all $(v_1,\dots,v_d,0)^T\in\R^{d+1}\setminus\{\mathbf{0}\}$.
For $\mathbf{A}_2$ we first consider the non--zero block matrix and verify that the rows are pairwise linearly dependent, since the $i$-th row is written as $u_i\mathbf{u}$. Thus, we have that $\mathbf{A}_2 \geq 0$ for all $(v_1,\dots,v_d,0)^T\in\R^{d+1}\setminus\{\mathbf{0}\}$. For $\mathbf{A}_3$ we yield when multiplying with elements $\tilde{\mathbf{v}} = (0,\dots,0,v_{d+1})^T\in\R^{d+1}\setminus\{\mathbf{0}\}$
\[
  \tilde{\mathbf{v}}^T\mathbf{A}_3\tilde{\mathbf{v}} = \left(1 + |\mathbf{u}|^2\right)\left(1 + 6|\mathbf{u}|^2\right)v_{d+1}^2 > 0
\]
and thus $\mathbf{A}_3 > 0$ for all $(0,\dots,0,v_{d+1})^T\in\R^{d+1}\setminus\{\mathbf{0}\}$.
Altogether we hence have shown the convexity of $\tilde{\eta} = -\eta$.\\
\newline
We now want to verify relation \eqref{entropy_flux_rel}. Therefore, we first calculate
the derivatives of the entropy flux with respect to the primitive quantities. These  are given by
\begin{align}
	\frac{\partial q_i}{\partial p} &= \frac{3}{4}p^{-1/4}u_i = \frac{3}{4}\frac{q_i}{p},\label{deriv_etaflux_p}\\
	\frac{\partial q_i}{\partial u_j} &= p^{3/4}\delta_{ij}.\label{deriv_etaflux_u}
\end{align}
From these we calculate the derivatives with respect to the conservative variables
\begin{align}
	\label{deriv_etaflux_cons}
      \partial_{w_j} q_i =
      \begin{cases}
       \dfrac{\delta_{ij}(3 + 2|\mathbf{u}|^2) - u_iu_j(1 + 2|\mathbf{u}|^2)}{4p^{\frac{1}{4}}\sqrt{1 + |\mathbf{u}|^2}(3 + 2|\mathbf{u}|^2)},\, &j=1,\ldots,d, \\[1.5em]
           -\dfrac{u_i}{4p^{\frac{1}{4}}}\dfrac{1 - 2|\mathbf{u}|^2}{3 + 2|\mathbf{u}|^2},\, &j=d+1.
       \end{cases}
\end{align}
We further need the Jacobians of the conservative fluxes $\mathbf{F}_k$ in the normal directions w.r.t.\ $x_k$. To this end we refer to \cite{kunik2024radially} where the Jacobians are given for any normal direction, see \cite[(B.5)]{kunik2024radially}.
Here the situation simplifies, since we only consider the normal directions $\mathbf{n} = \mathbf{e}_k, k=1,\dots,d$.
The flux Jacobian is determined by
\begin{align}
    \label{eq:flux-normal-jac}
    \mathbf{D}_\mathbf{w}\mathbf{F}_k =
    \begin{pmatrix}
        \overline{\mathbf{A}}   & \overline{\mathbf{a}}   \\
        \mathbf{e}_k^T & 0
    \end{pmatrix}
\end{align}
with
\begin{align*}
    \overline{\mathbf{A}} &=
    -2\frac{\sqrt{1 + |\mathbf{u}|^2}}{3 + 2|\mathbf{u}|^2}
    \begin{pmatrix}
    	0	& \dots			& 0\\
    		& \mathbf{u}^T	& \\
    	0	& \dots			& 0
    \end{pmatrix}
    + \frac{u_k}{\sqrt{1 + |\mathbf{u}|^2}}\mathbf{Id}_{d\times d}
    + \frac{1}{\sqrt{1 + |\mathbf{u}|^2}}
    \begin{pmatrix}
    	0		& 				& 0\\
    	\vdots	& \mathbf{u}		& \vdots\\
    	0		&				& 0
    \end{pmatrix}\\
    &+ 2\frac{u_k}{\sqrt{1 + |\mathbf{u}|^2}\left(3 + 2|\mathbf{u}|^2\right)}
	\begin{pmatrix}
		u_1^2	& u_1u_2	& \dots	 & u_1u_d\\
		u_1u_2	& u_2^2		& \dots	 & u_2u_d\\
		\vdots	&			& \ddots & \vdots\\
		u_1u_d	& \dots		& \dots	 & u_d^2
	\end{pmatrix},\\
	\overline{\mathbf{a}} &= \frac{1 + 2|\mathbf{u}|^2}{3 + 2|\mathbf{u}|^2}\mathbf{e}_k - \frac{4u_k}{3 + 2|\mathbf{u}|^2}\mathbf{u}.
\end{align*}
If \eqref{entropy_flux_rel} is written component wise we yield
\begin{align*}
	\partial_{w_i}q_k &= \nabla_{\overline{\mathbf{w}}}\eta^T\overline{\mathbf{A}}^{(i,\cdot)} + \partial_{w_{d+1}}\eta\delta_{ik},\;i=1,\dots,d,\\
	\partial_{w_{d+1}}q_k &= \nabla_{\overline{\mathbf{w}}}\eta^T\overline{\mathbf{a}}.
\end{align*}
Considering $i=1,\dots,d, i\neq k$ yield
\begin{align*}
	\overline{\mathbf{A}}^{(i,\cdot)} &= -2\frac{\sqrt{1 + |\mathbf{u}|^2}}{3 + 2|\mathbf{u}|^2}u_i\mathbf{e}_k + \frac{u_k}{\sqrt{1 + |\mathbf{u}|^2}}\mathbf{e}_i + 2\frac{u_iu_k}{\sqrt{1 + |\mathbf{u}|^2}\left(3 + 2|\mathbf{u}|^2\right)}\mathbf{u},\\
	\nabla_{\overline{\mathbf{w}}}\eta^T\overline{\mathbf{A}}^{(i,\cdot)} &= -2\frac{\sqrt{1 + |\mathbf{u}|^2}}{3 + 2|\mathbf{u}|^2}u_i\partial_{w_k}\eta + \frac{u_k}{\sqrt{1 + |\mathbf{u}|^2}}\partial_{w_i}\eta + 2\frac{u_iu_k}{\sqrt{1 + |\mathbf{u}|^2}\left(3 + 2|\mathbf{u}|^2\right)}\nabla_{\overline{\mathbf{w}}}\eta^T\mathbf{u}\\
	&= -\frac{1}{4}\frac{\eta}{p}\frac{u_iu_k\left(1 + 2|\mathbf{u}|^2\right)}{\left(1 + |\mathbf{u}|^2\right)\left(3 + 2|\mathbf{u}|^2\right)} = \partial_{w_i}q_k.
\end{align*}
Now we deal with $i=k=1,\dots,d$ and obtain
\begin{align*}
	\overline{\mathbf{A}}^{(k,\cdot)} &= -2\frac{\sqrt{1 + |\mathbf{u}|^2}}{3 + 2|\mathbf{u}|^2}u_k\mathbf{e}_k + \frac{u_k}{\sqrt{1 + |\mathbf{u}|^2}}\mathbf{e}_k + \frac{1}{\sqrt{1 + |\mathbf{u}|^2}}\mathbf{u} + \frac{2u_k^2}{\sqrt{1 + |\mathbf{u}|^2}\left(3 + 2|\mathbf{u}|^2\right)}\mathbf{u},\\
	\nabla_{\overline{\mathbf{w}}}\eta^T\overline{\mathbf{A}}^{(k,\cdot)} &= u_k\left(\frac{1}{\sqrt{1 + |\mathbf{u}|^2}}-2\frac{\sqrt{1 + |\mathbf{u}|^2}}{3 + 2|\mathbf{u}|^2}\right)\partial_{w_k}\eta + \left(\frac{1}{\sqrt{1 + |\mathbf{u}|^2}} + \frac{2u_k^2}{\sqrt{1 + |\mathbf{u}|^2}\left(3 + 2|\mathbf{u}|^2\right)}\right)\nabla_{\overline{\mathbf{w}}}\eta^T\mathbf{u}\\
	&= -\frac{1}{4}\frac{\eta}{p}\frac{1}{\left(1 + |\mathbf{u}|^2\right)\left(3 + 2|\mathbf{u}|^2\right)}\left(u_k^2 + 3|\mathbf{u}|^2 + 2|\mathbf{u}|^4 + 2u_k^2|\mathbf{u}|^2\right),\\
	\nabla_{\overline{\mathbf{w}}}\eta^T\overline{\mathbf{A}}^{(k,\cdot)} + \partial_{w_{d+1}}\eta &= -\frac{1}{4}\frac{\eta}{p}\frac{u_k^2\left(1 + |\mathbf{u}|^2\right) - \left(3 + 2|\mathbf{u}|^2\right)}{\left(1 + |\mathbf{u}|^2\right)\left(3 + 2|\mathbf{u}|^2\right)} = \partial_{w_k}q_k.
\end{align*}
We now consider the final case and get
\begin{align*}
	\nabla_{\overline{\mathbf{w}}}\eta^T\overline{\mathbf{a}} &= \frac{1 + 2|\mathbf{u}|^2}{3 + 2|\mathbf{u}|^2}\partial_{w_k}\eta - \frac{4u_k}{3 + 2|\mathbf{u}|^2}\nabla_{\overline{\mathbf{w}}}\eta^T\mathbf{u}
	= -\frac{1}{4}\frac{\eta}{p}\frac{u_k\left(1 - 2|\mathbf{u}|^2\right)}{\left(1 + |\mathbf{u}|^2\right)\left(3 + 2|\mathbf{u}|^2\right)} = \partial_{w_{d+1}}q_k.
\end{align*}
Thus, \eqref{entropy_flux_rel} is verified.
\printbibliography
%
\end{document}